\documentclass[psamsfonts,a4paper,reqno]{amsart}
\usepackage[left=1in,right=1in,bottom=1in,top=1in]{geometry}
\usepackage{tikz,amsthm,amsmath,amstext,amssymb,amscd,epsfig,euscript, mathrsfs, dsfont,pspicture,multicol, mathtools,graphpap,graphics,graphicx,times,enumerate,subfig,sidecap,wrapfig,color}
\usepackage{hyperref}
\usepackage{bm}
\usepackage{pgfplots}

\def \rr {{\mathbb R}}
\def \cc {{\mathbb C}}
\def \nn {{\mathbb N}}

\def \zz {{\mathbb Z}}

\DeclareGraphicsExtensions{.png}
 \numberwithin{equation}{section}

\def\var{\varphi}

\def\l{\lambda}

\def\z{\zeta}
\def\ep{\epsilon}

\def\g{\gamma}
\def\t{\tau}

\def\G{\Gamma}

\def\O{\Omega}
\def\o{\omega}
\def\S{\Sigma}

\def \tn {\tilde{n}}
\def \ty {\tilde{y}}
\def \tx {\tilde{x}}

\def \Hm {H_{\text{MIT}}(m)}
\def\sS{{\mathcal{S}}} 
\def\cP{{\mathcal{P}}} 

\def \cc {{\mathbb{C}}}

\newtheorem{definition}{Definition}[section]
\newtheorem{theorem}{Theorem}[section]
\newtheorem{proposition}{Proposition}[section]
\newtheorem{lemma}{Lemma}[section]
\newtheorem{corollary}{Corollary}[section]
\newtheorem{remark}{Remark}[section]
\newtheorem{notation}{Notation}[section]


\title[A Poincar\'e-Steklov map for the MIT bag model]{A Poincar\'e-Steklov map for the MIT bag model}

\author[B. BENHELLAL]{BADREDDINE BENHELLAL \textsuperscript{1}}
\address{\textsuperscript{1}Basque Center for Applied Mathematics, Alameda Mazarredo 14, 48009\\
and Departamento de Matem\'aticas, Universidad del Pa\' is Vasco, Barrio Sarriena s/n 48940 Leioa, Spain\\
 and Institut de Math\'ematiques de Bordeaux, UMR 5251, Universit\' e de Bordeaux 33405 Talence Cedex, France.}
\email{\textsuperscript{1} benhellal.badreddine@ehu.eus and badreddine.benhellal@u-bordeaux.fr}
\author[V. BRUNEAU]{VINCENT BRUNEAU\textsuperscript{2}}
\address{\textsuperscript{2}Institut de Math\'ematiques de Bordeaux, UMR 5251, Universit\' e de Bordeaux,  33405 Talence Cedex,  France.}
\email{\textsuperscript{2}vbruneau@math.u-bordeaux.fr}
\author[M. ZREIK]{MAHDI ZREIK\textsuperscript{3}}
\address{\textsuperscript{3}Institut de Math\'ematiques de Bordeaux, UMR 5251, Universit\' e de Bordeaux 33405 Talence Cedex,  France\\ and Departamento de Matem\'aticas, Universidad del Pa\' is Vasco, Barrio Sarriena s/n 48940 Leioa, Spain.}
	\email{\textsuperscript{3}mahdi.zreik@ehu.eus and mahdi.zreik@math.u-bordeaux.fr}
	\date{\today}
\subjclass[2010]{Primary: 35Q40; Secondary 35P05, 81Q10, 81Q20}
\keywords{Poincar\'e-Steklov operators, the MIT bag Model, $h$-Pseudodifferential operators, Large coupling
limits}
	
\begin{document}
	\begin{abstract} The purpose of this paper is to introduce and study Poincar\'e-Steklov (PS) operators associated to the Dirac operator $D_m$ with the so-called  MIT bag boundary condition. In a domain $\O\subset\rr^3$, for a complex number $z$ and for $U_z$ a solution of $(D_m-z)U_z=0$, the associated PS operator maps the value of $\Gamma_- U_z$, the MIT bag boundary value of $U_z$, to $\Gamma_+ U_z$, where $\Gamma_\pm$ are projections along the boundary $\partial\O$ and $(\Gamma_ - + \Gamma_+) = t_{\partial\O}$ is the trace operator on $\partial\O$.

  In the first part of this paper, we show that the PS operator is a zero-order pseudodifferential operator and give its principal symbol. In the second part, we study the PS operator when the mass $m$ is large, and we prove that it fits into the framework of $1/m$-pseudodifferential operators,  and we derive some important properties, especially its semiclassical principal symbol.  Subsequently, we apply these results to establish a Krein-type resolvent formula for the Dirac operator $H_M= D_m+ M\beta 1_{\rr^3\setminus\overline{\O}}$ for large masses  $M>0$, in terms of the resolvent of the MIT bag operator on $\O$. With its help, the large coupling convergence with a convergence rate of $\mathcal{O}(M^{-1})$ is shown.
	\end{abstract}
\maketitle	
		\section{Introduction}
	\subsection*{Motivation}
	Boundary integral operators have played a key role in the study of many boundary value problems  for partial differential equations arising in various areas of mathematical physics, such as electromagnetism, elasticity,  and potential theory. In particular, they are used as a tool for proving the existence of solutions as well as for their construction by means of integral equation methods, see, e.g., \cite{FJR,JK1,JK2,Verc}.

The study of boundary integral operators has also been the motivation for the development of various tools and branches of mathematics, e.g., Fredholm theory, Singular integral and Pseudodifferential operators. Moreover, it turned out that functional analytic and spectral properties of some of these operators are strongly related to the regularity and geometric properties of surfaces, see for example \cite{HMT,HMMPT}. A typical and well-known example which occurs in many applications is the Dirichlet-to-Neumann (DtN) operator.  In the classical setting of a bounded domain $\O\subset\rr^d$ with a smooth boundary,  the DtN operator, $\mathcal{N}$, is defined by 
  \begin{align*}
 \begin{split}
\mathcal{N}:\mathit{H}^{1/2}(\partial\O)  \longrightarrow  \mathit{H}^{-1/2}(\partial\O),\quad g\longmapsto  \mathcal{N}g =\G_NU(g),
\end{split}
\end{align*}
where $U(g)$ is the harmonic extension of $g$ (i.e., $\Delta U(g)=0$ in $\O$ and $\G_D U=g$ on $\partial\O$). Here $\G_{ D}$ and $\G_{N}$ denote the Dirichlet and the Neumann traces, respectively.  In this setting, it is well known that the DtN operator fits into the framework of pseudodifferential operators, see e.g., \cite{T1}. Moreover, from the viewpoint of the spectral theory,  several geometric properties of the eigenvalue problem for the DtN operator  (such as isoperimetric inequalities, spectral asymptotics and geometric invariants) are closely related to the theory of minimal surfaces \cite{FS}, as well as the problem of determining a complete Riemannian manifold with boundary from the Cauchy data of harmonic functions, see \cite{LTU} (see also the survey \cite{GP} for further details).

The main goal of this paper is to introduce a Poincar\'e-Steklov map for the Dirac operator (i.e., an analogue of the DtN map for the Laplace operator) and to study its (semiclassical) pseudodifferential properties. Our main motivation for considering this operator is that it arises naturally in the study of the well-known Dirac operator with the MIT bag boundary condition, $\Hm$, which will be rigorously defined below.

\subsection*{Description of main results}To give a rigorous definition of the operator we are dealing with in this paper and  go more into details, we need to introduce some notations. Given $m>0$, the free Dirac operator $D_m$ on $\rr^3$ is defined by  $D_m:=- i \alpha \cdot\nabla + m\beta$, where   
\begin{align*}
	\alpha_k&=\begin{pmatrix}
	0 & \sigma_j\\
	\sigma_j & 0
	\end{pmatrix}\quad \text{for } j=1,2,3,
	\quad
	\beta=\begin{pmatrix}
	\mathit{I}_2 & 0\\
	0 & -\mathit{I}_2
	\end{pmatrix},\quad 
\mathit{I}_2 :=\begin{pmatrix}
1& 0 \\
0 & 1
\end{pmatrix}, \\
&\text{and }\,
\sigma_1=\begin{pmatrix}
	0 & 1\\
	1 & 0
	\end{pmatrix},\quad \sigma_2=
	\begin{pmatrix}
	0 & -i\\
	i & 0
	\end{pmatrix} ,
	\quad
	\sigma_3=\begin{pmatrix}
	1 & 0\\
	0 & -1
	\end{pmatrix},
	\end{align*}  
 are the family of Dirac and Pauli matrices. As usual, we use the notation $\alpha \cdot x=\sum_{j=1}^{3}\alpha_j x_j$  for $x=(x_1,x_2,x_3)\in\rr^3$.  We recall that $D_m$ is self-adjoint in $\mathit{L}^2(\rr^3)^4$ with $\mathrm{dom}(D_m)=\mathit{H}^1(\rr^3)^4$ (see, e.g., \cite[ subsection~1.4]{Tha}), and for the spectrum and the continuous spectrum, we have:
	\begin{align*}
	\mathrm{Sp}(D_m) =\mathrm{Sp}_{\mathrm{cont}}(D_m)=(-\infty,-m]\cup [m,+\infty).
	\end{align*}
	Let $\O\subset\rr^3$ be a domain with a compact smooth boundary $\partial\O$, let $n$ be the outward unit normal to $\O$, and let $\G_\pm$ and $P_\pm$ be the trace mappings and the orthogonal projections, respectively,  defined by
 \begin{align*}
	\G_\pm=P_{\pm}\G_D: \mathit{H}^{1}(\O)^4  \longrightarrow  P_\pm\mathit{H}^{1/2}(\partial\O)^4\,\,\text{ and }\,\, P_{\pm}:=\frac{1}{2}\left(I_{4}\mp i\beta(\alpha\cdot n(x))\right), \quad x\in\partial\O.
\end{align*}
In the present paper, we investigate the specific case of the Poincar\'e-Steklov (PS for short) operator, $\mathscr{A}_m$,  defined by 
 \begin{align*}
\mathscr{A}_{m}:P_-\mathit{H}^{1/2}(\partial\O)^4  \longrightarrow  P_+\mathit{H}^{1/2}(\partial\O)^4,\quad g\longmapsto  \mathscr{A}_m(g) =\G_+U_z, 
\end{align*}
where $z$ belongs to the resolvent set of the MIT bag operator on $\O$ (i.e.,  $z\in\rho(\Hm)$),    $U_z\in\mathit{H}^{1}(\O)^4$  is the unique solution to the following elliptic boundary problem:
\begin{equation}\label{Dj}
	\left\{
	\begin{aligned}
	(D_m-z)U_z&=0, \quad \text{ in } \O,\\
	\G_-U_z&= g, \quad \text{ on } \partial\O.
	\end{aligned}
	\right.
	\end{equation}
	

We point out that in the R-matrix theory and the embedding method for the Dirac equation, similar operators linking on $ \partial\O$ values of the upper and lower components of the spinor wavefunctions have been studied in \cite{Sz, Ag, AR, BS}.
It corresponds to a different boundary condition (the trace of the upper/lower components) which is not necessarily  elliptic.
As far as we know, such operators for the MIT bag boundary condition have not been studied yet.


Let us now briefly describe the contents of the present paper. Our results are mainly concerned with the pseudodifferential properties of $\mathscr{A}_m$ and their applications. Thus, our first goal is to show that  $\mathscr{A}_m$ fits into the framework of pseudodifferential operators.  In Section \ref{SEC4},   we show that when the mass $m$ is fixed and $z\in\rho(D_m)$, then the Poincar\'e-Steklov operator $\mathscr{A}_m$ is a classical homogeneous pseudodifferential operators of order $0$, and that 
\begin{align*}
\mathscr{A}_m&=S\cdot\left(\frac{\nabla_{\partial\O}\wedge n}{\sqrt{-\Delta_{\partial\O}}}\right)P_- \quad\mathrm{mod}\, \mathit{Op} \sS^{-1}(\partial\O),
\end{align*}
where $S=i(\alpha\wedge\alpha)/2$ is the spin angular momentum,   $\nabla_{\partial\O}$ and $\Delta_{\partial\O}$  are, respectively,   the surface gradient and the Laplace-Beltrami operator on $\partial\O$ (equipped with the Riemann metric induced by the euclidian one in $\rr^3$) and $\mathit{Op} \sS^{-1}$ is the classical class of pseudodifferential operators of order $-1$ (see Theorem \ref{Th SymClas} for details). 
For $D_{\partial\O}$, the extrinsically defined Dirac operator introduced in Section \ref{sDirac}, we also have:
\begin{align*}
\mathscr{A}_m&=D_{\partial\O} \, \left(-\Delta_{\partial\O} \right)^{-\frac12} P_- \quad\mathrm{mod}\, \mathit{Op} \sS^{-1}(\partial\O).
\end{align*}
 The proof of the above result is based on the fact that we have an explicit solution of the system \eqref{Dj} for any $z\in\rho(D_m)$, and in this case  the PS operator takes the following layer potential form:
 \begin{align}\label{def2PSS}
\mathscr{A}_m= -P_+\beta\left(\beta/2 +\mathscr{C}_{z,m}\right)^{-1}P_-,
\end{align}
where $\mathscr{C}_{z,m}$ is the Cauchy operator associated with $(D_m-z)$ defined on $\partial\O$ in the principal value sense (see Subsection \ref{Subintope} for the precise definition). So the starting point of the proof is to analyze the pseudodifferential properties of the Cauchy operator. In this sense, we show that $2\mathscr{C}_{z,m}$ is equal, modulo $\mathit{Op} \sS^{-1}(\partial\O)$, 
to $\alpha\cdot(\nabla_{\partial\O}(-\Delta_{\partial\O})^{-1/2})$. Using this,  the explicit layer potential description of $\mathscr{A}_m$, and the symbol calculus, we then prove that $\mathscr{A}_m$ is a pseudodifferential operator and catch its principal symbol (see Theorem \ref{Th SymClas}).

While the above strategy allows us to capture the pseudodifferential character of $\mathscr{A}_m$, but unfortunately it does not allow us to trace the dependence on the parameter $m$,  and it also imposes a restriction on the spectral parameter $z$ (i.e., $z\in\rho(D_m)$), whereas $\mathscr{A}_m$ is well-defined for any $z\in\rho(\Hm)$. In Section \ref{SEC5}, we address the $m$-dependence of the pseudodifferential properties of $\mathscr{A}_m$ for any $z\in\rho(\Hm)$.  Since we are mainly concerned with large masses $m$ in our application, we  treat this problem from the semiclassical point of view, where  $h=1/m \in(0,1]$ is the semiclassical parameter.  In fact, we show in Theorem \ref{pseudo de Am}  that $\mathscr{A}_{1/h}$  admits a semiclassical approximation, and that
 \begin{equation*}
\mathscr{A}_{1/h} = \frac{hD_{\partial\O}}{\sqrt{-h^2\Delta_{\partial\O}+I}+I} \, P_- \quad\mathrm{mod}\,  \,  h \mathit{Op}^h \sS^{-1}(\partial\O).
 \end{equation*}
The main idea of the proof is to use the system  \eqref{Dj} instead of the explicit formula \eqref{def2PSS}, and it is based on the following two steps. The first step is to construct a local approximate solution for the pushforward of the system \eqref{Dj} of the form
\begin{equation*}
U^h(\tx,x_3) = \mathit{Op}^h (A^h(\cdot, \cdot, x_3) ) g = \frac{1}{2\pi}  \int_{\rr^2} A^h(\tx, h \xi , x_3) e^{i y \cdot \xi} \hat{g}(\xi) \mathrm{d} \xi ,\quad  (\tx,x_3)\in\rr^2\times[0,\infty),
\end{equation*}
where $A^h$  belongs to a specific symbol class and has the following asymptotic expansion 
\begin{align*}A^h(\tx,  \xi , x_3) \sim \sum_{j \geq 0} h^j A_j(\tx, \xi , x_3).
\end{align*}
The second step is to show that when applying  the trace mapping  $\G_+$ to the pull-back of $U^h(\cdot , 0)$ it coincides locally with $\mathscr{A}_{1/h}$ modulo a regularizing and negligible operator.  At this point, the properties of the MIT bag operator become crucial, in particular,  the regularization property of its resolvent which allows us to achieve this second step, as we will see in Section \ref{SEC5}. The MIT bag operator on $\O$ is the Dirac operator on $\mathit{L}^{2}(\O )^4$ defined by 
\begin{align*}
\Hm \psi=D_m \psi, \quad \forall\psi\in\mathrm{dom}(\Hm):=  \left\{ \psi \in\mathit{H}^{1}(\O )^4 :   \G_- \psi=0 \text{ on }\partial\O \right\}.
\end{align*}
It is well-known that $(\Hm, \mathrm{dom}(\Hm)$ is self-adjoint when $\O$ is smooth, see, e.g., \cite{OBV}. In Section \ref{SEC3}, we briefly discuss the basic spectral properties of $\Hm$ when $\O$ is a domain with compact Lipschitz boundary (see Theorem \ref{MITLipshictz}). Moreover, in Theorem \ref{the de regu} we establish regularity results concerning  the regularization property of the resolvent and the Sobolev regularity of the eigenfunctions of $H_{\text{MIT}}$. In particular, we prove that $( \Hm -z)^{-1}$ is bounded from $ \mathit{H}^{n}(\O )^4$ into $\mathit{H}^{n+1}(\O )^4\cap\mathrm{dom}(\Hm)$, for all $n\geqslant1$. 

Motivated by the natural way in which the PS operator is related to the MIT bag operator, and to illustrate its usefulness,  we consider in Section \ref{SEC6} the large mass problem for the self-adjoint Dirac operator $H_M= D_m +M\beta 1_{\mathcal{U}}$, where $\mathcal{U}=\rr^3\setminus \overline{\O}$.  Indeed, it is known that, in the limit $M\rightarrow\infty$, every eigenvalue of $\Hm$ is a limit of eigenvalues of $H_M$, cf. \cite{ALMR,MOP} (see also  \cite{BCLS,BB1,SV} for the two-dimensional setting). Moreover, it is shown in \cite{BCLS,BB1} that the two-dimensional analogue of $H_M$ convergences to the two-dimensional analogue of $\Hm$ in the norm resolvent sense with a convergence rate of $\mathcal{O}(M^{-1/2})$. 

The main goal of Section \ref{SEC6} is to address the following question:
Let $M_0>0$ be large enough and fix $M\geqslant M_0$ and $z\in\rho(\Hm)\cap\rho(H_M)$. Given $f\in\mathit{L}^2(\rr^3)^4$  such that $f=0$ in $\rr^3\setminus\overline{\O}$,  and $U\in \mathit{H}^1(\rr^3)^4$, what is the boundary value problem on $\O$ whose solutions closely approximate those of $(D_m +M\beta 1_{\rr^3\setminus\overline{\O}} -z)U= f$?

It is worth noting that the answer to this question becomes trivial if one establishes an explicit formula for the resolvent of $H_M$. Having in mind the connection between the Dirac operators  $H_M$ and $\Hm$, this leads us to address the following question: \textit{for $M$ sufficiently large, is it possible to relate the resolvents of  $H_M$ and $H_{\mathrm{MIT}}$ via a Krein-type resolvent formula?} In Theorem \ref{Formulederes}, which is the main result of Section \ref{SEC6}, we establish a Krein-type resolvent formula for $H_M$ in terms of the resolvent of $\Hm$.  
The key point to establish this result is to treat the elliptic problem $(H_M-z)U=f\in\mathit{L}^2(\rr^3)^4$ as a transmission problem (where $\G_\pm U_{| \O}=\G_\pm U_{| \rr^3\setminus{\O}}$ are the transmission conditions) and to use the semiclassical properties of the Poincar\'e-Steklov operators in order to invert the auxiliary operator $\Psi_M(z)$ acting on the boundary $\partial\O$ (see Theorem \ref{Formulederes} for the precise definition). In addition, we prove an adapted Birman-Schwinger principle relating the eigenvalues of $H_M$ in the gap $(-(m+M),m+M)$ with a spectral property of $\Psi_M(z)$. With their help, we show in Corollary \ref{une autre est}  that the restriction of $U$ on $\O$ satisfies the elliptic problem
\begin{equation*}
\left\{
\begin{aligned}
	(D_m- z) U_{| \O} &=f \quad&\text{ in }\,\Omega,\\
	\G_-U_{| \O}&= \mathscr{B}_{M} \G_+R_{MIT}(z)f  \,\,&\text{ on }\, \partial\O,\\
	\G_+U_{| \O}&= \G_+R_{MIT}(z)f  +\mathscr{A}_{m}\G_-v \,\,&\text{ on }\, \partial\O,
\end{aligned}
\right.
\end{equation*}
where $\mathscr{B}_{M}$ is  a semiclassical pseudodifferential operators of order $0$. Here, the semiclassical parameter is $1/M$. Moreover, we show that the convergence of $H_M$ to $H_{\text{MIT}}$ in the norm-resolvent sense indeed holds with a convergence rate of $\mathcal{O}(M^{-1})$, which improves previous works, see Proposition \ref{Resolvent convergence}.  The most important ingredient in proving these results is the use of the Krein formula relating the resolvents of $H_M$ and $\Hm$, as well as regularity estimates for the PS operators (see Theorem \ref{PseudESt}) and layer potential operators (see Lemma \ref{Estimations des op} for details).

\subsection*{Organization of the paper} The paper is organized as follows.  Sections  \ref{SEC2} and  \ref{SEC3} are devoted to preliminaries for the sake of completeness and self-containedness of the paper. In Section \ref{SEC2} we set up some notations and we recall some basic properties of boundary integral operator associated with $(D_m-z)$. Section \ref{SEC3} is devoted to the study of the MIT bag operator, where we gather its basic properties in Theorem  \ref{MITLipshictz} and we establish the regularization property of its resolvent in  Theorem \ref{the de regu}. In Section \ref{SEC4} we establish Theorem \ref{Th SymClas}, proving that  the PS operator is a classical pseudodifferantial operator. Then, in Section \ref{SEC5}  we study the PS operator from viewpoint of semiclassical pseudodifferantial operators, the main result being  Theorem \ref{pseudo de Am}. Finally,  Section \ref{SEC6}  is devoted to the study of the large mass problem for the operator $H_M$. There, we prove Theorem \ref{Formulederes} regarding the Krein-type resolvent formula and we solve the large mass problem, and  Proposition \ref{Resolvent convergence} on  the resolvent convergence.
	
	\section{ Preliminaries}\label{SEC2}
	In this section we gather some well-known results about boundary integral operators. 
	We also recall some properties of symbol classes and their associated pseudodifferential operators.  Before proceeding further, however, we need to introduce some notations that we will use in what follows. 
\subsection{Notations}\label{SubNot} Throughout this paper we will write $a\lesssim b$ if there is $C>0$ so that  $a\leqslant C b$ and $a\lesssim_h b$  if the constant $C$ depends on the parameter $h$. As usual, the letter $C$ stands for some constant which may change its value at different occurrences.

For a bounded or unbounded Lipschitz domain  $\O\subset\rr^3$, we write $\partial\O$ for its boundary and we denote by $n$ and $\mathrm{\sigma}$ the outward pointing normal to $\O$ and the surface measure on $\partial\O$, respectively.  By $\mathit{L}^2(\rr^{3})^4:=\mathit{L}^2(\rr^{3};\cc^4)$ (resp. $\mathit{L}^2(\O)^4:=\mathit{L}^2(\O,\cc^4)$) we denote the usual $\mathit{L}^2$-space over $\rr^3$ (resp. $\O$), and we let  $r_{\O}: \mathit{L}^{2}(\rr^3)^4\rightarrow \mathit{L}^{2}(\O)^4$ be the restriction operator on $\O$ and $e_{\O}:\mathit{L}^{2}(\O)^4 \rightarrow \mathit{L}^{2}(\rr^3)^4$ its adjoint operator, i.e., the extension by $0$ outside of $\O$.

For $s\in[0,1]$, we define the usual  Sobolev space $\mathit{H}^s(\rr^{d})^4$  as
\begin{align*}
 \mathit{H}^s(\rr^{d})^4:=\{ u\in\mathit{L}^2(\rr^{d})^4: \int_{\rr^d}(1+|\xi|^2)^{s} \left|\mathcal{F}[u](\xi)\right|^2\mathrm{d}\xi<\infty\},
\end{align*}
where $\mathcal{F}:\mathit{L}^2(\rr^{d}) \rightarrow \mathit{L}^2(\rr^{d})$ is the unitary Fourier-Plancherel operator, and we let $\mathit{H}^{s}(\O)^4$ to be standard $\mathit{L}^2$-based  Sobolev space of order $s$. By $\mathit{L}^2(\partial\O)^4:=\mathit{L}^2(\partial\O, \mathrm{d\sigma})^4$ we denote the usual $\mathit{L}^2$-space over $\partial\O$.  If $\O$ is a $\mathit{C}^2$-smooth domain with a compact boundary $\partial\O$, then the Sobolev space of order $s\in(0,1]$ along the boundary, $\mathit{H}^{s}(\partial\O)^4$,  is defined using local coordinates representation on the surface $\partial\O$. As usual, we use the symbol $\mathit{H}^{-s}(\partial\O)^4$ to denote the dual space of $\mathit{H}^{s}(\partial\O)^4$. We denote by  $t_{\partial\O}:\mathit{H}^1(\O)^4\rightarrow \mathit{H}^{1/2}(\partial\O)^4$ the classical trace operator, and by  $\mathcal{E}_{\O}:  \mathit{H}^{1/2}(\partial\O)^4\rightarrow  \mathit{H}^1(\O)^4$ the extension operator, that is 
$$t_{\partial\O}\mathcal{E}_{\O}[f]=f,\quad \forall f\in \mathit{H}^{1/2}(\partial\O)^4.$$
Throughout the current paper, we denote by $P_\pm$ the orthogonal projections defined by 
 \begin{align}\label{Projections}
	P_{\pm}:=\frac{1}{2}\left(I_{4}\mp i\beta(\alpha\cdot n(x))\right), \quad x\in\partial\O.
\end{align}
We use the symbol $\mathit{H}(\alpha,\O)$ for the Dirac-Sobolev space on a smooth domain $\O$ defined as
 \begin{align}\label{Sobolev Dirac}
\mathit{H}(\alpha,\O)=\{ \varphi\in\mathit{L}^2(\O)^4: (\alpha\cdot\nabla)\varphi\in\mathit{L}^2(\O)^4 \},
\end{align}
which is a Hilbert space (see \cite[Section 2.3]{OBV}) endowed with the following scalar product
\begin{align*}
\langle \varphi, \psi \rangle_{\mathit{H}(\alpha,\O)}=\langle \varphi, \psi\rangle_{\mathit{L}^{2}(\O)^4}+\langle (\alpha\cdot\nabla)\varphi, (\alpha\cdot\nabla)\psi\rangle_{\mathit{L}^{2}(\O)^4},\quad \varphi,\psi\in\mathit{H}(\alpha,\O).
\end{align*} 
 We also recall that the trace operator $t_{\partial\O}$ extends into a continuous map $t_{\partial\O}:\mathit{H}(\alpha,\O)\rightarrow \mathit{H}^{-1/2}(\partial\O)^4$. Moreover, if $v\in\mathit{H}(\alpha,\O)$  and $t_{\partial\O} v\in \mathit{H}^{1/2}(\partial\O)^{4}$, then $v\in \mathit{H}^1(\O)^4$, cf.  \cite[Proposition 2.1 \& Proposition 2.16]{OBV}.
\subsection{Boundary integral operators}\label{Subintope} The aim of this part is to introduce boundary integral operators associated with the fundamental solution of the free Dirac operator $D_m$  and to summarize some of their well-known properties.  

For $z\in\rho(D_m)$, with the convention that $\mathrm{Im}\sqrt{z^2-m^2}>0$,  the fundamental solution of $(D_m-z)$ is given by
\begin{align}\label{defsolo}
	\phi^{z}_{m}(x)=\frac{e^{i\sqrt{z^2-m^2}|x|}}{4\pi|x|}\left(z +m\beta+( 1-i\sqrt{z^2-m^2}|x|)i\alpha\cdot\frac{x}{|x|^2}\right), \quad \forall  x\in\rr^3\setminus\{0\}.
\end{align}
 We define the potential operator $\Phi^{\O}_{z,m}: \mathit{L}^2(\partial\O)^4  \longrightarrow \mathit{L}^2(\O)^4$ by
\begin{align}\label{Phi}
  \Phi^{\O}_{z,m}[g](x) =\int_{\partial\O} \phi^{z}_{m}(x-y)g(y)\mathrm{d\sigma}(y), \quad \text{for all } x\in\O,
\end{align}
 and the Cauchy operators $\mathscr{C}_{z,m}: \mathit{L}^2(\partial\O)^4  \longrightarrow \mathit{L}^2(\partial\O)^4$ as the singular integral operator acting as
\begin{align}\label{CauchyOpe}
\mathscr{C}_{z,m}[f](x)&=  \lim\limits_{\rho\searrow 0}\int_{|x-y|>\rho}\phi^{z}_{m}(x-y)f(y)\mathrm{d\sigma}(y), \quad \text{for } \mathrm{\sigma}\text{-a.e. }x\in\partial\O,\,f\in\mathit{L}^2(\partial\O)^4.
\end{align}
It is well known that $\Phi^{\O}_{z,m}$ and $\mathscr{C}_{z,m}$ are bounded and everywhere defined (see, for instance,  \cite[Section. 2]{AMV1}), and that 
\begin{align}\label{PRI}
((\alpha\cdot n)\mathscr{C}_{z,m})^2=\left(\mathscr{C}_{z,m}(\alpha\cdot n)\right)^2=-4\mathit{I}_4, \quad \forall z\in\rho(D_m),
\end{align}
holds in $ \mathit{L}^2(\partial\O)^4$,  cf. \cite[Lemma 2.2]{AMV2}. In particular, the inverse $\mathscr{C}_{z,m}^{-1}=-4(\alpha\cdot n)\mathscr{C}_{z,m}(\alpha\cdot n)$ exists and is bounded and everywhere defined. Since we have $\phi^{z}_m(y-x)^{\ast}=\phi^{\overline{z}}_m(x-y)$ for all $z\in\rho(D_m)$, it follows that  $\mathscr{C}^{\ast}_{z,m}=\mathscr{C}_{\overline{z},m}$ as operators in $ \mathit{L}^2(\partial\O)^4$. In particular, $\mathscr{C}_{z,m}$ is self-adjoint in $\mathit{L}^2(\partial\O)^4$ for all $z\in(-m,m)$.

Next,  recall that the trace of the single layer operator, $S_z$,  associated with the  Helmholtz operator $(-\Delta+m^2-z^2)\mathit{I}_4$ is defined, for every $f\in\mathit{L}^{2}(\partial\O)^4$ and $z\in\rho(D_m)$,  by
\begin{align*}
S_z[f](x):=\int_{\partial\O}\frac{e^{i\sqrt{z^2-m^2}|x-y|}}{4\pi|x-y|}\mathit{I}_4f(y)\mathrm{d\sigma}(y), \quad \text{for } x\in\partial\O. 
\end{align*}
It is well-known that $S_z$ is bounded from $\mathit{L}^{2}(\partial\O)^4$ into $\mathit{H}^{1/2}(\partial\O)^4$, and it is a positive operator in $\mathit{L}^{2}(\partial\O)^4$ for all $z\in(-m,m)$, cf. \cite[Lemma 4.2]{AMV2}. Now we define the operator  $\Lambda^{z}_{m}$ by
\begin{align*}
\Lambda^{z}_ {m}=\frac{1}{2}\beta + \mathscr{C}_{z,m},\quad\text{ for all }  z\in\rho(D_m),
\end{align*}
which is clearly a bounded operator from $\mathit{L}^{2}(\partial\O)^4$ into itself. \\

In the next lemma  we collect  the main  properties of the operators $\Phi^{\O}_{z,m}$, $\mathscr{C}_{z,m}$ and $\Lambda^{z}_{m}$.
\begin{lemma}\label{prop of C} Assume that $\O$ is $\mathit{C}^2$-smooth. Given $z\in\rho(D_m)$ and let $\Phi^{\O}_{z,m}$, $\mathscr{C}_{z,m}$ and $\Lambda^{z}_{ m}$ be as above. Then the following hold true:
\begin{itemize}
   \item[($\mathrm{i}$)] The operator $\Phi^{\O}_{z,m}$ is bounded from $\mathit{H}^{1/2}(\partial\O)^4$ to $\mathit{H}^{1}(\O)^4$, and  extends into a bounded operator from $\mathit{H}^{-1/2}(\partial\O)^4$ to $\mathit{H}(\alpha,\O)$. Moreover,  it holds that 
   \begin{align}\label{traceof Phi}
   t_{\partial\O} \Phi^{\O}_{z,m}[f]=  \left(-\frac{i}{2}(\alpha\cdot n) + \mathscr{C}_{z,m}\right)[f], \quad \forall f\in \mathit{H}^{1/2}(\partial\O)^4.
   \end{align}
  \item[($\mathrm{ii}$)] The operator $\mathscr{C}_{z,m}$ gives rise to a bounded operator $\mathscr{C}_{z,m} : \mathit{H}^{1/2}(\partial\O)^4  \longrightarrow \mathit{H}^{1/2}(\partial\O)^4$. 
  \item[($\mathrm{iii}$)]  The operator $\Lambda^{z}_{ m} : \mathit{H}^{1/2}(\partial\O)^4  \longrightarrow \mathit{H}^{1/2}(\partial\O)^4$ is bounded invertible for all  $z\in\rho(D_m)$. 
\end{itemize}
\end{lemma}
\textbf{Proof.} $(\mathrm{i})$ The proof of  the boundedness of   $\Phi^{\O}_{z,m}$ from $\mathit{H}^{1/2}(\partial\O)^4$ into $\mathit{H}^{1}(\O)^4$ is contained in \cite[Proposition 4.2]{BH}, and the jump formula \eqref{traceof Phi} is proved in \cite[Lemma 3.3]{AMV1} in terms of non-tangential limit which coincides (almost everywhere in $\partial\O$) with the trace operator for functions in $\mathit{H}^1(\O)^4$. The boundedness of $\Phi^{\O}_{z,m}$ from $\mathit{H}^{-1/2}(\partial\O)^4$ to $\mathit{H}(\alpha,\O)$ is established in \cite[Theorem 2.2]{OBV}.

Since $n$ is smooth, it is clear from $(\mathrm{i})$ that $\mathscr{C}_{z,m}$ is bounded from $ \mathit{H}^{1/2}(\partial\O)^4 $ into itself, which proves  $(\mathrm{ii})$.  As consequence  we also obtain that $\Lambda^{z}_{ m}$ is bounded from $ \mathit{H}^{1/2}(\partial\O)^4 $ into itself.  Now, the invertibility of $\Lambda^{z}_{ m}$ in $ \mathit{H}^{1/2}(\partial\O)^4 $ for $z\in\cc\setminus\rr$ is shown in \cite[Lemma 3.3 (iii)]{BEHL2}, see also \cite[Lemma 3.12]{BHM}. To complete the proof of $(\mathrm{iii})$, note that if $f\in\mathit{L}^{2}(\partial\O)^4$ is such that $\Lambda^{z}_{ m}[f]\in\mathit{H}^{1/2}(\partial\O)^4 $, then a simple computation shows that 
\begin{align*}
\mathit{H}^{1/2}(\partial\O)^4\ni(\Lambda^{z}_{m})^2[f]= \left(1/4+(\mathscr{C}_{z,m})^2+ (m\mathit{I}_{4}+z\beta) S_z\right)[f],
\end{align*}
which means that $f\in\mathit{H}^{1/2}(\partial\O)^4$. From the above computation we see that $\Lambda^{z}_{ m}$ is invertible from $\mathit{H}^{1/2}(\partial\O)^4$ into itself for all $z\in(-m,m)$, since $((\mathscr{C}_{z,m})^2+ (m\mathit{I}_{4}+z\beta)S_z)$ is a positive operator.  This completes the proof of the lemma. \qed
\\
\begin{remark}\label{remLip} Note that if $\O$ is a Lipschitz domain with a compact boundary, then for all $z\in\rho(D_m)$ the operators $\mathscr{C}_{z,m}$ and $\Lambda^{z}_{ m}$ are bounded  from $\mathit{L}^2(\partial\O)^4$ into itself (see, e.g, \cite[Lemma 3.3]{AMV1}), and since $\Lambda^{z}_{ m}$ is an injective Fredholm operator (see the proof of \cite[Theorem 4.5]{BB3}) it follows that it is also invertible in $\mathit{L}^2(\partial\O)^4$.  Note also that, thanks to \cite[Lemma 5.1 and Lemma 5.2]{BHSS},  we know that the mapping $\Phi^{\O}_{z,m}$ defined by \eqref{Phi}  is bounded from  $\mathit{L}^2(\partial\O)^4$ to  $\mathit{H}^{1/2}(\O)^4$,  $ t_\S \Phi^{\O}_{z,m}[g]\in \mathit{L}^{2}(\partial\O)^4$ and the formula \eqref{traceof Phi}  still holds true for all $g\in \mathit{L}^{2}(\partial\O)^4$. 
\end{remark}	
\subsection{Symbol classes  and Pseudodifferential operators}\label{SymPseuOpSur} We recall here the basic facts concerning the classes of pseudodifferential operators that will serve in the rest of the paper. 

Let $\mathscr{M}_{4}(\cc)$ be the set of $4\times 4$ matrices over $\cc$.  For $d\in\mathbb{N}^{\ast}$ we let $\sS^{m}(\rr^d\times\rr^d)$ be the standard symbol class  of order $m\in\rr$ whose elements are matrix-valued functions $a$ in the space $\mathit{C}^{\infty}(\rr^d\times\rr^d;\mathscr{M}_{4}(\cc))$ such that 
\begin{align*}
|\partial^{\alpha}_{x}\partial^{\beta}_{\xi}a(x,\xi) |\leqslant C_{\alpha\beta}(1+|\xi|^2)^{m-|\beta|}, \quad \forall (x,\xi)\in\rr^d\times\rr^d, \,\,\forall \alpha\in\nn^{d},\,\,\forall \beta\in\nn^{d}.
\end{align*}
Let $\mathscr{S}(\rr^d)$ be the Schwartz class of functions. Then,  for each  $a\in  \sS^{m}(\rr^d\times\rr^d)$ and any  $h \in (0,1]$,  we associate a semiclassical pseudodifferential operator $Op^h(a):\mathscr{S}(\rr^d)^4\rightarrow \mathscr{S}(\rr^d)^4$ via the standard formula
$$Op^h(a)u(x)=\frac{1}{(2\pi)^{d/2}}\int_{\rr^d}e^{i\xi\cdot x}a(x,h\xi) \mathcal{F}[u](\xi)\mathrm{d}\xi, \quad \forall u\in \mathscr{S}(\rr^d)^4.$$
If $a\in \sS^0 (\rr^d\times\rr^d)$, then Calder\'on-Vaillancourt theorem's (see, e.g., \cite{CAVA}) yields that $Op^h(a)$ extends to a bounded operator from $\mathit{L}^{2}(\rr^d)^4$ into itself, and there exists $C, N_C>0$  such that 
\begin{align}\label{Calderon-Vaillancourt}
\left| \left| Op^{h}(a)\right|\right|_{\mathit{L}^{2}\rightarrow \mathit{L}^{2}}\leqslant C \max_{ |\alpha+\beta|\leqslant N_C}\left| \left| \partial_{x}^{\alpha}\partial_{\xi}^{\beta}a \right|\right|_{\mathit{L}^\infty}.
\end{align}

Given a $\mathit{C}^{\infty}$-smooth domain  $\O\subset\rr^3$ with a compact boundary $\S=\partial\O$. Then $\S$ is a $2$-dimensional parameterized surface, which in the sense of differential geometry,   can also be viewed as a smooth $2$-dimensional manifold immersed into $\rr^3$.  Thus, $\S$ can be  covered by an atlas $\mathbb{A}=\{ (U_j, V_j,\var_j)| j\in\{1,\cdots, N\} \}$ (i.e., a collection of smooth charts)  where $N\in \mathbb{N}^{\ast}$.  That is
\begin{align*}
\S= \bigcup\limits_{j= 1}^{N}U_j ,
\end{align*}
and for each  $j\in\{1,\cdots, N\}$,  $U_j$ is an open set of $\S$, $ V_j\subset\rr^2$ is an open set of the parametric space $\rr^2$, and $\var_j:U_j\rightarrow V_j$ is a $\mathit{C}^\infty$- diffeormorphism. Moreover, by definition of a smooth manifold,  if $U_j\cap U_k\neq\emptyset$  then 
$$\varphi_k \circ (\varphi_j)^{-1} \in C^\infty \Big( \varphi_j (U_{j} \cap U_{k}) ; \, \varphi_k (U_{j} \cap U_{k}) \Big).$$
As usual,  the pull-back $(\var_j^{-1})^{\ast}$  and the  pushforward $\var_j^{\ast}$ are defined by 
 $$(\var_j^{-1})^{\ast}u= u\circ \var_j^{-1} \quad \text{ and } \quad  \var_j^{\ast}v= v \circ \var_j ,$$ 
for  $u$ and $v$  functions on $U_j$ and  $V_j$, respectively. We also recall that a function $ u$ on $\S$ is said to be in the class $\mathit{C}^k(\S)$ if for every chart the pushforward has the property  $(\var_j^{-1})^{\ast}u \in\mathit{C}^k(V_j)$.

Following Zworski \cite[Part 4.]{MZ},  we define pseudodifferential operators on the boundary $\S$  as follows:
\begin{definition} Let $\mathscr{A}: \mathit{C}^\infty(\S)^4 \rightarrow \mathit{C}^\infty (\S)^4$ be a continuous linear operator. Then $\mathscr{A}$ is said to be a $h$-pseudodifferential operator of order $m\in\rr$ on $\S$, and we write  $\mathscr{A}\in \mathit{Op} ^h \sS^m(\S)$,
if 
\begin{itemize}
  \item[(1)] for every chart $ (U_j, V_j,\var_j)$  there exists a symbol $a\in  \sS^m$ such that 
  \begin{align*}
\psi_1\mathscr{A}(\psi_2 u)= \psi_1\var_j^{\ast} \mathit{Op}^h(a) (\var_j^{-1})^{\ast}(\psi_2 u),
\end{align*}
for any $\psi_1,\psi_2\in\mathit{C}^{\infty}_{0}(U_j)$ and $u\in \mathit{C}^\infty(\S)^4$.
  \item [(2)] for all $\psi_1,\psi_2\in\mathit{C}^\infty (\S)$ such that  $\mathrm{supp}(\psi_1)\cap \mathrm{supp}(\psi_2)=\emptyset$ and for all $N\in\mathbb{N}$ we have 
 \begin{align*}
 \left\|  \psi_1\mathscr{A}\psi_2 \right\|_{\mathit{H}^{-N}(\S)^4\rightarrow \mathit{H}^{N}(\S)^4}= \mathcal{O}(h^\infty).
\end{align*}
\end{itemize}
For $h$ fixed (for example $h=1$), $\mathscr{A}$ is called a pseudodifferential operator.
\end{definition}

Since the study of a given pseudodifferential operator on $\S$ reduces to local study on local charts, in what follows, we will recall below the specific local coordinates and surface geometry notations we will use in the rest of the paper. 

We always fix an open set $U\subset\S$, and we let  $\chi:V\rightarrow \rr$  to be a $\mathit{C}^{\infty}$-function  (where $V\subset\rr^2$ is open) such that its graph coincides with $U$. Set $\varphi(\tx)=(\tx,\chi(\tx))$, then for $x\in U$ we write $x=\varphi(\tx)$ with $\tx\in V$.  Here and also in what follows, $\partial_1\chi$ and $\partial_2\chi$ stand for the partial derivatives $\partial_{\tx_1}\chi$ and $\partial_{\tx_2}\chi$, respectively. Recall that the first fundamental form, $\mathrm{I}$, and  the metric tensor $G(\tx)=(g_{jk}(\tx))$, have the following forms:
\begin{align*}
\mathrm{I}&=g_{11}d\tx_1^2+2g_{12}d\tx_1d\tx_2 +g_{22}d\tx_2^2,\\
G(\tx)=\left(g_{jk}(\tx)\right)&=\begin{pmatrix}
	g_{11} & g_{12}\\
	g_{21} & g_{22}
	\end{pmatrix}(\tx):=\begin{pmatrix}
	1+|\partial_1\chi|^2 &\partial_1\chi\partial_2\chi\\
	\partial_1\chi\partial_2\chi &1+|\partial_2\chi|^2
	\end{pmatrix}(\tx).
\end{align*}
As $G(\tx)$ is symmetric, it follows that it is diagonalizable  by an orthogonal matrix. Indeed, let 
\begin{align}\label{defPmatrix}
Q(\tx):=\begin{pmatrix}
	\frac{|\partial_2\chi|}{|\nabla\chi|} &\frac{\partial_1\chi\partial_2\chi}{|\partial_2\chi||\nabla\chi|}\\
	-\frac{\partial_1\chi\partial_2\chi}{|\partial_2\chi||\nabla\chi|}&	\frac{|\partial_2\chi|}{|\nabla\chi|},
	\end{pmatrix}\begin{pmatrix}
        1 &0\\
	0&g^{-1/2}
	\end{pmatrix}(\tx).
\end{align}
where $g$ stands for  the determinant of  $G$. Then, it is straightforward to  check that 
\begin{align}\label{lesrelations}
Q^{\mathrm{t}}GQ(\tx)=\mathit{I}_2,\quad  QQ^{\mathrm{t}}(\tx)=G(\tx)^{-1}=:\left(g^{jk}(\tx)\right),\quad\mathrm{det}(Q)=\mathrm{det}(Q^{\mathrm{t}})=g^{-1/2}.
\end{align}

\subsection{Operators on the boundary $\Sigma = \partial \O$} \label{sDirac}
As above, we consider $\Sigma = \partial \O$ the boundary of a smooth bounded domain $\O$. On $\Sigma$ equipped with the Riemann metric induced by the euclidian one in $\rr^3$, we consider the Laplace-Beltrami operator $-\Delta_\Sigma$ and  the  surface gradient $\nabla_\Sigma = \nabla - n (n \cdot \nabla) $ where $n$ is the unit normal to the surface pointing outside $\O$. With the notation of the previous section, in local coordinates, these operators are pseudodifferential operators with respective principal symbols
\begin{equation}\label{symbolLG}
p_{\Delta_\Sigma} (\tx , \xi) = \langle G(\tx)^{-1}\xi,\xi\rangle, \qquad 
p_{\nabla_\Sigma} (\tx , \xi) = \xi_G:= \begin{pmatrix}
	G(\tx)^{-1}\xi \\
	\langle \nabla\chi(\tx),G(\tx)^{-1}\xi\rangle
	\end{pmatrix}.
\end{equation}

Let us now introduce $D_\Sigma$,  the extrinsically defined Dirac operator. To any $x \in \rr^3$ we associate the matrix $\alpha (x) = \alpha \cdot x $, where $\alpha = (\alpha_1, \alpha_2, \alpha_3)$. For $H_1$ the mean curvature of $\Sigma$,  $D_\Sigma$ is given by (for more details see Appendix B of \cite{MOP}):
\begin{align*}
D_\Sigma =  - \alpha (n) \, \alpha (\nabla_\Sigma)  + \frac{H_1}2.
\end{align*}
It is a  pseudodifferential operator with principal symbol:
$$ p_{D_\Sigma} (\tx , \xi) = - i \alpha (n^\varphi (\tx)) \, \alpha ( \xi_G),$$
where $n^\varphi = \varphi^* n$. Using the anticommutation relations of the Dirac's matrices (see also (\ref{AlphaSigma})) and that $n \cdot \xi_G=0$ we have:
$$
p_{D_\Sigma} (\tx , \xi) = - i \alpha \cdot n^\varphi (\tx) \, \alpha \cdot \xi_G = S \cdot (\xi_G \wedge n^\varphi (\tx) ).$$
Moreover for $\overline{\xi} := \begin{pmatrix}
	\xi \\
	0
	\end{pmatrix}$, we have:
	$ \overline{\xi} = \xi_G + (\overline{\xi} \cdot n^\varphi ) n^\varphi $.
	Thus,  in local coordinates, the principal symbol of $D_\Sigma$ is also:
	\begin{equation}\label{symbolD}
p_{D_\Sigma} (\tx , \xi) = S \cdot ( \overline{\xi} \wedge n^\varphi (\tx) ).
\end{equation}
Let us also point out the relationship between the principal symbols of $\Delta_\Sigma$ and $ D_\Sigma$:
\begin{equation}\label{symbolLD}
 | \overline{\xi} \wedge n^\varphi (\tx) |^2 =    \langle G(\tx)^{-1}\xi,\xi\rangle.  
\end{equation}

\section{ Basic properties of the MIT bag model}\label{SEC3}
	In this section, we give a brief review of the basic spectral properties of the Dirac operator with the MIT bag boundary condition  on Lipschitz domains.  Then, we establish some results concerning the regularization properties of the resolvent and the Sobolev regularity of the eigenfunctions in the case of smooth domains. 
	
Let $\mathcal{U}\subset\rr^3$ be a Lipschitz domain with a compact boundary $\partial\mathcal{U}$. Then, for $m>0$,  the Dirac operator with the MIT bag boundary condition on $\mathcal{U}$, $(\Hm,\mathrm{dom}(\Hm))  $, or simply the MIT bag  operator, is defined on the domain 
\begin{align*}
\mathrm{dom}(\Hm):=  \left\{ \psi \in\mathit{H}^{1/2}(\mathcal{U} )^4 : (\alpha\cdot\nabla)\psi\in\mathit{L}^{2}(\mathcal{U} )^4  \text{ and }  P_-t_{\partial\mathcal{U}} \psi=0 \text{ on }\partial\mathcal{U} \right\},
\end{align*}
by $\Hm\psi=D_m \psi$, for all $\psi\in \mathrm{dom}(\Hm)$, and where the boundary condition holds in $\mathit{L}^{2}(\partial\mathcal{U} )^4$. Here $P_\pm$ are the orthogonal projections defined by \eqref{Projections}.

The following theorem gathers the basic properties of  the MIT bag  operator. We mention that some of theses properties  are well-known in the case of smooth domains, see, e.g., \cite{ALMR,ALR,AMSV,BHM,OBV}.
\begin{theorem}\label{MITLipshictz} The operator $(\Hm,\mathrm{dom}(\Hm))$ is self-adjoint and we have 
\begin{align}\label{KreinMIT}
(\Hm-z)^{-1}= r_{\mathcal{U}}(D_m -z)^{-1}e_{\mathcal{U}} - \Phi_{z,m}^{\mathcal{U}}(\Lambda^z_m)^{-1}t_{\partial\mathcal{U}}(D_m -z)^{-1}e_{\mathcal{U}}, \quad \forall z\in\rho(D_m).
	\end{align} 
Moreover,  the following statements hold true:
\begin{itemize}
\item[(i)] If  $\mathcal{U}$ is bounded, then $\mathrm{Sp}(\Hm)=\mathrm{Sp}_{\mathrm{disc}}(\Hm)\subset\rr\setminus[-m,m]$.
\item[(ii)] If  $\mathcal{U}$ is unbounded, then $\mathrm{Sp}(\Hm)=\mathrm{Sp}_{\mathrm{ess}}(\Hm)=(-\infty,-m]\cup [m,+\infty)$. Moreover, if $\mathcal{U}$ is connected then $\mathrm{Sp}(\Hm)$ is purely continuous.
\item[(iii)]	Let $z\in\rho(\Hm)$ be such that $2|z|<m$, then for all $f\in\mathit{L}^2(\mathcal{U})^4$ it holds that
 $$\left\|(\Hm-z)^{-1}f\right\|_{\mathit{L}^{2}(\mathcal{U} )^4}\lesssim \frac{1}{m} \left\|f\right\|_{\mathit{L}^2(\mathcal{U})^4}.$$
\end{itemize}
\end{theorem}

\textbf{Proof.} Let $\var,\psi \in\mathrm{dom}(\Hm)$, then by density arguments we get the Green's formula 
 \begin{align}\label{2Geen formula}
\langle (-i\alpha\cdot\nabla)\varphi, \psi \rangle_{\mathit{L}^{2}(\mathcal{U})^4}-\langle \varphi,  (-i\alpha\cdot\nabla)\psi\rangle_{\mathit{L}^{2}(\mathcal{U})^4}=\langle (-i\alpha\cdot n)t_{\partial\mathcal{U}}\varphi, t_{\partial\mathcal{U}}\psi\rangle_{\mathit{L}^{2}(\partial\mathcal{U})^4}.
\end{align} 
Since $P_-t_{\partial\mathcal{U}}\var=P_-t_{\partial\mathcal{U}}\psi=0$ and $P_\pm(\alpha\cdot n)=(\alpha\cdot n)P_\mp$, it follows that   
 \begin{align*}
\langle (-i\alpha\cdot\nabla)\varphi, \psi \rangle_{\mathit{L}^{2}(\mathcal{U})^4}-\langle \varphi,  (-i\alpha\cdot\nabla)\psi\rangle_{\mathit{L}^{2}(\mathcal{U})^4}=\langle P_+(-i\alpha\cdot n)P_+t_{\partial\mathcal{U}}\varphi, P_+t_{\partial\mathcal{U}}\psi\rangle_{\mathit{L}^{2}(\partial\mathcal{U})^4}=0.
\end{align*} 
Consequently, we obtain 
\begin{align*}\label{}
  \langle \Hm\varphi,\psi\rangle_{\mathit{L}^2(\mathcal{U})^4}- \langle \varphi,\Hm\psi \rangle_{\mathit{L}^2(\mathcal{U})^4}&=  \langle D_m\var,\psi\rangle_{\mathit{L}^2(\mathcal{U})^4}- \langle \var,D_m\psi \rangle_{\mathit{L}^2(\mathcal{U})^4}\\ 
& = \langle (-i\alpha\cdot\nabla)\varphi, \psi \rangle_{\mathit{L}^{2}(\mathcal{U})^4}-\langle \varphi,  (-i\alpha\cdot\nabla)\psi\rangle_{\mathit{L}^{2}(\mathcal{U})^4}=0.
\end{align*}
Therefore $(\Hm,\mathrm{dom}(\Hm))$ is symmetric. Now, thanks to \cite[Proposition 4.3]{BB3} we know that the MIT bag operator defined on the domain 
\begin{align}\label{domlipmit}
\mathscr{D}=\left\{ \psi =u+\Phi_{m}^{\mathcal{U}}[g], u\in\mathit{H}^1(\mathcal{U})^4, g\in\mathit{L}^2(\partial\mathcal{U})^4: P_- t_{\partial\mathcal{U}}\psi=0 \text{ on } \partial\mathcal{U} \right\},
\end{align}
by $\Hm(u+\Phi_{m}^{\mathcal{U}}[g])=D_m u$, for all $(u+\Phi_{m}^{\mathcal{U}}[g])\in \mathscr{D}$,  is a self-adjoint operator.  As $\Hm$ is symmetric on $\mathrm{dom}(\Hm)$ we deduce that  $\mathrm{dom}(\Hm)\subset\mathscr{D}$.  Now, by Remark \ref{remLip} we also get that $\mathscr{D}\subset \mathrm{dom}(\Hm)$ which proves the equality $\mathscr{D}= \mathrm{dom}(\Hm)$, and thus  $(\Hm,\mathrm{dom}(\Hm))$ is self-adjoint.  Next, we check the resolvent formula \eqref{KreinMIT}. So let $f\in\mathit{L}^2(\mathcal{U})^4$, $z\in\rho(D_m)$ and set 
\begin{align*}
\psi= r_{\mathcal{U}}(D_m -z)^{-1}e_{\mathcal{U}}f - \Phi_{z,m}^{\mathcal{U}}(\Lambda^z_m)^{-1}t_{\partial\mathcal{U}}(D_m -z)^{-1}e_{\mathcal{U}}f .
	\end{align*}
Since  $(D_m -z)^{-1}e_{\mathcal{U}}$ is bounded from $ \mathit{L}^2(\mathcal{U})^4 $ into $ \mathit{H}^1(\rr^3)^4$ and $(\Lambda^z_m)^{-1}$ is well-defined by Remark \ref{remLip}, it follows that  
	$$u:=r_{\mathcal{U}}(D_m -z)^{-1}e_{\mathcal{U}}f\in \mathit{H}^1(\mathcal{U})^4 \,\,\text{ and } \,\,g:=-(\Lambda^z_m)^{-1}t_{\partial\mathcal{U}}(D_m -z)^{-1}e_{\mathcal{U}}f\in \mathit{L}^2(\partial\mathcal{U})^4,$$ 
which entails that $\psi \in\mathit{H}^{1/2}(\mathcal{U} )^4$ and that $(\alpha\cdot\nabla)\psi\in\mathit{L}^2(\mathcal{U} )^4$. Next, using Lemma \ref{prop of C}-$(\mathrm{i})$ and Remark \ref{remLip} we easily get 
\begin{align*}
t_{\partial\mathcal{U}} \psi&=P_+\beta(\Lambda^z_m)^{-1}t_{\partial\mathcal{U}}(D_m -z)^{-1}e_{\mathcal{U}}f,
\end{align*}
 thus $P_-t_{\partial\mathcal{U}} \psi=0$ on $\partial\mathcal{U}$, which means that $\psi\in\mathrm{dom}(\Hm)$. Since $(D_m-z)\Phi_{z,m}^{\mathcal{U}}[g]=0$ holds in $\mathcal{U}$, it follows that   $(\Hm-z)\psi=f$ and the formula  \eqref{KreinMIT} is proved. 

Now, we are going to prove assertions  $(\mathrm{i})$ and  $(\mathrm{ii})$. First, note that for $\psi\in\mathrm{dom}(\Hm)$ a straightforward application of the Green formula \eqref{2Geen formula} yields that
\begin{align}\label{Quadraticform}
		\left\|\Hm\psi\right\| ^{2}_{\mathit{L}^2(\mathcal{U})^4} =& \left\| (\alpha\cdot\nabla)\psi \right\|^2_{\mathit{L}^2(\mathcal{U})^4}+ m^2 \left\| \psi \right\|^2_{\mathit{L}^2(\mathcal{U})^4}+ m\left\| P_+t_{\partial\mathcal{U}}\psi \right\|^2_{\mathit{L}^2(\partial\mathcal{U})^4}.
\end{align}
Thus $\left\|\Hm\psi\right\| ^{2}_{\mathit{L}^2(\mathcal{U})^4}\geqslant m^2\left\| \psi \right\|^2_{\mathit{L}^2(\mathcal{U})^4}$ which yields that $\mathrm{Sp}(\Hm)\subset(-\infty,-m]\cup [m,+\infty)$.  Note that this fact can be seen immediately from the formula \eqref{KreinMIT}. Next, we show that $\{ -m, m\}\notin\mathrm{Sp}_{\mathrm{disc}}(\Hm)$. Assume that there is $0\neq \psi \in\mathrm{dom}(\Hm)$ such that $(\Hm-m)\psi=0$ in $\mathcal{U}$. Then, from \eqref{Quadraticform} we have that 
\begin{align*}		
 \left\| (-i\alpha\cdot\nabla)\psi \right\|^2_{\mathit{L}^2(\mathcal{U})^4} + m\left\| P_+t_{\partial\mathcal{U}}\psi \right\|^2_{\mathit{L}^2(\partial\mathcal{U})^4}=0.
\end{align*}
Since $m>0$ it follows that $ P_+t_{\partial\mathcal{U}}\psi=0$, and thus $t_{\partial\mathcal{U}}\psi=0$. Using this and the above equation, an integration by parts (using density arguments)  gives 
\begin{align*}		
 \left\| \nabla\psi \right\|_{\mathit{L}^2(\partial\mathcal{U})^4}= \left\| (-i\alpha\cdot\nabla)\psi \right\|_{\mathit{L}^2(\mathcal{U})^4}=0.
\end{align*}
From this we conclude that $\psi$ vanishes identically, which contradicts the fact that $\psi\neq0$, and thus $m\notin\mathrm{Sp}_{\mathrm{disc}}(\Hm)$. Following the same lines as above we also get that  $-m\notin\mathrm{Sp}_{\mathrm{disc}}(\Hm)$. Thus, if  $\mathcal{U}$ is bounded, then the above considerations and the fact that $\mathrm{dom}(\Hm) \subset\mathit{H}^{1/2}(\mathcal{U} )^4$  is compactly embedded in $\mathit{L}^{2}(\mathcal{U} )^4$ yield that $\mathrm{Sp}(\Hm)=\mathrm{Sp}_{\mathrm{disc}}(\Hm)\subset\rr\setminus[-m,m]$, which shows the assertion $(\mathrm{i})$.\\
Lest us now complete the proof  of $(\mathrm{ii})$, so  suppose that $\mathcal{U}$ is unbounded.   We first show that $(-\infty,-m]\cup [m,+\infty)\subset \mathrm{Sp}_{\mathrm{ess}}(\Hm)$ by constructing Weyl sequences as in the case of half-space, see \cite[Theorem 4.1]{BB2}.  As $\mathcal{U}$ is unbounded it follows that there is $R_1>0$ such that the half-space $\{ x=(x_1,x_2,x_3)\in\rr^3: x_3>R_1\}$ is strictly contained in $\mathcal{U}$ and $ \rr^3\setminus\overline{\mathcal{U}}\subset B(0,R_1)$. Fix $\l\in(-\infty,-m)\cup (m,+\infty)$ and let $\xi= (\xi_1,\xi_2)$ be such that $|\xi|^2=\lambda^2-m^2$. We define the function  $\var:\rr^3 \longrightarrow \cc^4$ by
  \begin{align*}
\var( \overline{x},x_3)= \left(\frac{\xi_1-i\xi_2}{\lambda-m},0,0,1\right)^t e^{ i\xi\cdot\overline{x}}, \,\,\text{ with }\, \overline{x}=(x_1,x_2).
\end{align*}
Clearly we have $(D_m-\lambda)\var=0$. Now, fix $R_2>R_1$ and let $\eta\in \mathit{C}^{\infty}_{0}(\rr^2,\rr)$ and $\chi\in \mathit{C}^{\infty}_{0}(\rr,\rr)$ be such that $\mathrm{supp}(\chi)\subset[R_1, R_2]$. For $n\in\nn^\star$, we define the sequences of functions
\begin{align*}
 \varphi_{n}(\overline{x},x_3)=  n^{-\frac{3}{2}} \varphi(\overline{x},x_3)\eta(\overline{x}/n)\chi(x_3/n), \,\quad \text{for }(\overline{x},x_3) \in\mathcal{U}.
\end{align*}
Then, it is easy to check that $ \varphi_{n}\in\mathit{H}^1_0(\mathcal{U})\subset\mathrm{dom}(\Hm)$, $ (\varphi_{n})_{n\in\mathbb{N}^\star}$  converges weakly to zero, and that 
\begin{align*}\label{}
    \Vert \varphi_{n}\Vert_{\mathit{L}^2(\mathcal{U})^4}^2=\frac{2\lambda}{\l-m}\Vert \eta\Vert^{2}_{\mathit{L}^2(\rr^2)}\Vert \chi\Vert^{2}_{\mathit{L}^2(\rr)}>0, \quad \frac{  \Vert\left(D_m-\l\right) \varphi_{n}\Vert_{\mathit{L}^2(\mathcal{U})^4}}{ \Vert \varphi_{n}\Vert_{\mathit{L}^2(\mathcal{U})^4}}\xrightarrow[n\to\infty]{}0,
\end{align*}
for more details see the proof of \cite[Theorem 4.1]{BB2}. Therefore, Weyl's criterion yields that 
$$(-\infty,-m)\cup (m,+\infty)\subset \mathrm{Sp}_{\mathrm{ess}}(\Hm).$$
Since the spectrum of a self-adjoint operator is closed, we then get the first statement of $(\mathrm{ii})$. Now, if we assume in addition that $\mathcal{U}$ is connected, then using the same arguments as in the proof of  \cite[Theorem 3.7]{AMV2} (i.e., using Rellich's lemma and the unique continuation property)  one can verifies that $\Hm$ has no eigenvalues in $\rr\setminus[-m,m]$. As  $\{ -m, m\}\notin\mathrm{Sp}_{\mathrm{disc}}(\Hm)$ it follows that $\Hm$ has a purely continuous spectrum.

Now we prove $(\mathrm{iii})$. Let $\psi\in\mathrm{dom}(\Hm)$, then \eqref{Quadraticform} yields that  $\left\|\Hm\psi\right\| ^{2}_{\mathit{L}^2(\O)^4}\geqslant m^2\left\| \psi \right\|^2_{\mathit{L}^2(\O)^4}$, and thus 
\begin{align*}
m \left\| \psi \right\|_{\mathit{L}^2(\mathcal{U})^4} \leqslant  \left\|\Hm\psi\right\|_{\mathit{L}^2(\mathcal{U})^4}  \leqslant  \left\|(\Hm-z)\psi\right\|_{\mathit{L}^2(\mathcal{U})^4} +|z| \left\| \psi \right\|_{\mathit{L}^2(\mathcal{U})^4} 
\end{align*}
Therefore, for $2|z|<m$ with $z\in\rho(\Hm)$,  we get that $ \left\| \psi \right\|_{\mathit{L}^2(\mathcal{U})^4}\leqslant 2 m^{-1}  \left\|(\Hm-z)\psi\right\|_{\mathit{L}^2(\mathcal{U})^4} $. Thus, $(\mathrm{iii})$ follows by taking $\psi=(\Hm-z)^{-1}f$.
 \qed

\begin{remark}\label{RemMITLipshictz} We mention that the above statement on the self-adjointness can also be deduced from \cite[Theorem 5.4]{BHSS}. We also mention that the MIT bag operator defined on the domain $\mathscr{D}$ given by \eqref{domlipmit} is still self-adjoint for less regular domains, cf. \cite{BB3} for more details. 
\end{remark}	
\begin{remark}\label{Sobolevremar} Note that if $\mathcal{U}$ is in the class of H\"{o}lder's domains $\mathit{C}^{1,\o}$, with $\o\in(1/2,1)$, then  $\Hm$ is self-adjoint and $\mathrm{dom}(\Hm):=  \left\{ \psi \in\mathit{H}^{1}(\mathcal{U} )^4 :  P_-t_{\partial\mathcal{U}} \psi=0 \text{ on }\partial\mathcal{U} \right\}$, see \cite[Theorem 4.3]{BB3}  for example. 
\end{remark}

Now we establish regularity results  which concerns the regularization property of the resolvent and the Sobolev regularity of the eigenfunctions of $\Hm$. The first statement of the following theorem will be crucial in Section \ref{SEC5} when studying the semiclassical pseudodifferential properties of the Poincar\'e-Steklov operator.

\begin{theorem}\label{the de regu} Let $k\geqslant1$ be an integer and assume that $\mathcal{U}$ is $\mathit{C}^{2+k}$-smooth. Then the following statements hold true:  
\begin{itemize}
\item[(i)]  The mapping $(\Hm-z)^{-1}:\mathit{H}^k(\mathcal{U})^4  \longrightarrow  \mathit{H}^{k+1}(\mathcal{U})^4\cap \mathrm{dom}(\Hm)$ is well-defined and bounded for all $m>0$ and all  $z\in\rho(\Hm)$. In particular, for $m_0>0$ and all $z\in\rho(H_{\text{MIT}}(m_0))\cap\rho(\Hm)$ we have 
$$\| (\Hm-z)^{-1}\|_{\mathit{H}^k(\mathcal{U})^4  \longrightarrow  \mathit{H}^{k+1}(\mathcal{U})^4}\lesssim 1,$$
uniformly on $m\geqslant m_0$.
\item[(ii)] If $\phi$ is an eigenfunction associated with an eigenvalue $z\in\mathrm{Sp}(\Hm)$, i.e.,  $(\Hm-z)\phi=0$,   then $\phi\in \mathit{H}^{1+k}(\mathcal{U})^4$. In particular, if  $\mathcal{U}$  is $\mathit{C}^{\infty}$-smooth,  then $\phi\in\mathit{C}^{\infty}(\mathcal{U})^4$.
\end{itemize}
\end{theorem}

To prove this theorem we need the following classical regularity result.
\begin{proposition}\label{Neumann} Let $k$ be a nonnegative integer. Assume that $\mathcal{U}$ is $\mathit{C}^{3+k}$-smooth and $u\in \mathit{H}^1(\mathcal{U})$. If $u$ solves the Neumann problem 
\begin{align*}
-\Delta u= f\in \mathit{H}^k(\mathcal{U})\quad \text{and}\quad \partial_n u=g\in\mathit{H}^{1/2+k}(\partial\mathcal{U}),
\end{align*} 
then $u\in \mathit{H}^{2+k}(\mathcal{U})$.
\end{proposition}
\textbf{Proof.} First, assume that $k=0$. As $\mathcal{U}$ is $\mathit{C}^{3}$-smooth we know that the Neumann trace $ \partial_n : \mathit{H}^2(\mathcal{U}) \rightarrow\mathit{H}^{1/2}(\partial\mathcal{U})$ is surjective. Thus,  there is $G\in  \mathit{H}^2(\mathcal{U})$ such that $\partial_n G=g$ in $\partial\mathcal{U}$. Note that the function $\tilde{u}=u-G$ satisfies the homogeneous Neumann problem 
\begin{align*}
-\Delta \tilde{u}= f +\Delta G\, \text{ in } \mathcal{U}\quad \text{and} \quad\partial_n \tilde{u}=0 \,\text{ on } \partial\mathcal{U}.
\end{align*} 
Therefore, $\tilde{u}\in \mathit{H}^2(\mathcal{U})$ by \cite[Theorem 5, p. 217]{Mik},  which implies that $u\in  \mathit{H}^2(\mathcal{U})$ and this proves the result for $k=0$. If $k\geqslant1$, then the result follows by  \cite[Theorem 2.5.1.1]{Gri}.\qed\\

\textbf{Proof of Theorem \ref{the de regu}.}  The theorem will be proved by induction on $k$. First, we show  $(\mathrm{i})$, so fix $z\in\rho(\Hm)$ and assume that $k=1$.  Let $\phi=(\phi_1,\phi_2)^{\top}\in \mathrm{dom}(\Hm)$ be such that $(D_m-z)\phi=f$ in $\mathcal{U}$, with $f=(f_1,f_2)^{\top}\in \mathit{H}^1(\mathcal{U})^4$. By assumption we have  $(\Delta+m^2-z^2)\phi=(D_m-z)f$  in $\mathcal{D}^{\prime}(\mathcal{U})^4$, and then in $\mathit{L}^2(\mathcal{U})^4$. We next  prove that $\partial_{n}\phi\in \mathit{H}^{1/2}(\partial\mathcal{U})^4$. To this end, consider $\mathcal{U}_{\ep}:=\{x\in\rr^3: \mathrm{dist}(x,\partial\mathcal{U})<\ep\}$ for $\ep>0$. Then, for $\delta>0$ small enough and $0<\ep \leqslant\delta$  the mapping $\Psi: \S\times(-\ep,\ep)\rightarrow \mathcal{U}_{\ep}$, defined by 
\begin{align}\label{the mapping Phi}
 \Psi(x_{\partial\mathcal{U}},t)= x_{\partial\mathcal{U}}+tn(x_{\partial\mathcal{U}}),\quad x_{\partial\mathcal{U}}\in\partial\mathcal{U},\, t\in(-\ep,\ep)
\end{align}
is a $\mathit{C}^{2}$-diffeomorphism  and  $\mathcal{U}_{\ep}:=\{x+tn(x) : x\in\partial\mathcal{U},\, t\in(-\ep,\ep) \}$.

 Let  $\widetilde{P_-}: \mathit{L}^2(\mathcal{U}_\ep\cap\mathcal{U})^4\rightarrow\mathit{L}^2(\mathcal{U}_\ep\cap\mathcal{U})^4$ be the bounded operator defined by  
$$\widetilde{P_-}\varphi( \Psi(x,t))=\frac{1}{2}(\mathit{I}_4+i\beta(\alpha\cdot n(x)))\varphi( \Psi(x,t)), \quad    \Psi(x,t)\in \mathcal{U}_\ep\cap\mathcal{U}.$$ 
Let $x_{\partial\mathcal{U}}^0$ be an arbitrary point on the boundary $\S$, fix $0<r<\ep/2$,  and let $\z: \rr^3\rightarrow [0,1]$ be a $\mathit{C}^{\infty}$-smooth and compactly supported function such that  $\z=1$  on $B(x_{\partial\mathcal{U}}^0,r)$ and   $\z=0$ on $\rr^3\setminus B(x^0_{\partial\mathcal{U}},2r)$. We claim that $\widetilde{P_-}\z\phi$ satisfies the elliptic problem
\begin{eqnarray*}
 \left\{
    \begin{split}\
    \displaystyle -\Delta(\widetilde{P_-}\z\phi) =  g \quad\text{in } \mathcal{U},\\
t_{\partial\mathcal{U}}(\widetilde{P_-}\z \phi) =  0\quad\text{on } \partial\mathcal{U},   
   \end{split}
  \right.
\end{eqnarray*}
with $g\in  \mathit{L}^2(\mathcal{U})^4$. Indeed, set $\mathcal{B}(x)= i\beta(\alpha\cdot n(x))$ for $x\in\partial\mathcal{U}$, and  observe that
 \begin{align*}\label{}
(D_m-z)(\widetilde{P_-}\z\phi)= \left( \widetilde{P_-}\z f +\frac{1}{2} [D_m,\z]\phi\right)+ \frac{1}{2}[D_m,\z \mathcal{B}]\phi=: I(\phi,f)+ \frac{1}{2}[D_m,\z \mathcal{B}]\phi.
\end{align*}
Since $n$ is $\mathit{C}^2$-smooth, $\z$ is infinitely differentiable and $\psi, f\in \mathit{H}^1(\mathcal{U})^4$, it is clear that $I(\phi,f)\in\mathit{H}^1(\mathcal{U})^4$ and $[D_m,\z \mathcal{B}]\phi\in\mathit{L}^2(\mathcal{U})^4$. Now, 
applying $(D_m+z)$ to the above equation yields that $-\Delta(\widetilde{P_-}\z\phi)=g$ with 
 \begin{align*}\label{}
g:=(z^2- m^2)\widetilde{P_-}\z \phi +( D_m+z)I(\phi,f)+\frac{z}{2} [D_m,\z \mathcal{B}]\phi +  \frac{1}{2}D_m[D_m,\z \mathcal{B}]\phi.
\end{align*}
As before, it is clear that the first three terms are square integrable. Next, observe that 
 \begin{align*}
 D_0[D_0,\z \mathcal{B}]\phi&= \{ D_0,  [D_0,\z \mathcal{B}] \} \phi -[D_0,\z \mathcal{B}](f- (m\beta-z)\phi)\\
 &=   [-\Delta,\z \mathcal{B}] \phi -[D_0,\z \mathcal{B}](f- (m\beta-z)\phi).
  \end{align*}
Using this together with the smoothness assumption on  $n$ and the fact  $(D_m-z)\phi=f\in \mathit{H}^1(\mathcal{U})^4$,  we easily see that $D_0[D_0,\z \mathcal{B}]\phi\in\mathit{L}^2(\mathcal{U})^4$. Hence,   $D_m[D_m,\z \mathcal{B}]\phi$ is square integrable, which means that  $g\in\mathit{L}^2(\mathcal{U})^4$.  As $P_-t_{\partial\mathcal{U}}\phi=0$ and $t_{\partial\mathcal{U}}(\widetilde{P_-}\z \phi)=t_{\partial\mathcal{U}}\z P_-t_{\partial\mathcal{U}}\phi=0$ on $\partial\mathcal{U}$,   by  \cite[Theorem 8.12 ]{GT} it follows that $\widetilde{P_-}\z\phi \in  \mathit{H}^2(\mathcal{U}_\ep\cap\mathcal{U})^4$, which implies that
 \begin{align*}
	\z(\phi_1 + i(\sigma\cdot n)\phi_2) \in \mathit{H}^2(B(x_{\partial\mathcal{U}}^0,2r)\cap\mathcal{U})^2 \quad\text{and}\quad \z(-i(\sigma\cdot n)\phi_1 + \phi_2) \in \mathit{H}^2(B(x_{\partial\mathcal{U}}^0,2r)\cap\mathcal{U})^2.
\end{align*}
  Consequently,  we get  
 \begin{align}
	\phi_1 + i(\sigma\cdot n)\phi_2 \in \mathit{H}^2(B(x_{\partial\mathcal{U}}^0,r)\cap\mathcal{U})^2  \quad\text{and}\quad -i(\sigma\cdot n)\phi_1 + \phi_2 \in \mathit{H}^2(B(x_{\partial\mathcal{U}}^0,r)\cap\mathcal{U})^2  \label{eqA1}.
\end{align}
Since $-i(\sigma\cdot\nabla)\phi_2= (z-m)\phi_1 +f_1$ and $-i(\sigma\cdot\nabla)\phi_1= (z+m)\phi_2 +f_2$ hold in $\mathit{H}^1(\mathcal{U})^2$,  it follows from  \eqref{eqA1} that 
\begin{align*}
 (\sigma\cdot\nabla)\phi_j\in   \mathit{H}^1(B(x_{\partial\mathcal{U}}^0,r))^2\quad\text{and}\,\,\, (\sigma\cdot\nabla)(\sigma\cdot n)\phi_j\in   \mathit{H}^1(B(x_{\partial\mathcal{U}}^0,r))^2, \quad j=1,2.
\end{align*}
Using this and the fact that  $n$ is $\mathit{C}^2$-smooth,  we easily get that 
\begin{align*}
(\sigma\cdot n) (\sigma\cdot\nabla)\phi_j + (\sigma\cdot\nabla)(\sigma\cdot n)\phi_j= \langle n,\nabla \rangle \mathit{I}_2\phi_j + F_j  \in   \mathit{H}^1(B(x_{\partial\mathcal{U}}^0,r))^2, 
\end{align*}
with $F_j  \in   \mathit{H}^1(B(x_{\partial\mathcal{U}}^0,r)\cap\mathcal{U})^2$. As a consequence, we get that $ \langle n,\nabla \rangle \mathit{I}_2\phi_j \in   \mathit{H}^1(B(x_{\partial\mathcal{U}}^0,r)\cap\mathcal{U})^2$. Since this holds true for all  $x_{\partial\mathcal{U}}^0\in\partial\mathcal{U}$, using the compactness of $\partial\mathcal{U}$  it follows that $\partial_{n}\phi\in \mathit{H}^{1/2}(\partial\mathcal{U})^4$.  Therefore, Propositions \ref{Neumann} yields that $\phi\in \mathit{H}^{2}(\mathcal{U})^4$. 

Next, assume $k\geqslant 2$,  $\mathcal{U}$ is $\mathit{C}^{2+k}$-smooth and $\phi, f\in \mathit{H}^k(\mathcal{U})^4$. Since $n$  is $\mathit{C}^{1+k}$-smooth and $\Psi$ defined by \eqref{the mapping Phi} is a  $\mathit{C}^{1+k}$-diffeomorphism,   following the same arguments as above we then conclude that $\partial_{n}\phi\in \mathit{H}^{k-1/2}(\S)^4$. Note also that $-\Delta\phi=(z^2-m^2)\phi+ (D_m-z)f \in \mathit{H}^{k-1}(\mathcal{U})^4$. Therefore, thanks to Propositions \ref{Neumann},  we conclude that  $\phi\in \mathit{H}^{k+1}(\mathcal{U})^4$, which proves the first statement of $(\mathrm{i})$.

Now,  the second statement of $(\mathrm{i})$ is a direct consequence of the first one, and this completes the proof of $(\mathrm{i})$.

Finally, the proof of the first statement of $(\mathrm{ii})$ follows the same lines as the one of $(\mathrm{i})$. In particular, if $\mathcal{U}$ is $\mathit{C}^{\infty}$-smooth, we then get  $\phi\in \mathit{H}^{k+1}(\mathcal{U})^4$ for any $k\geqslant 0$, which implies that $\phi$ is infinitely  differentiable in $\mathcal{U}$, and the theorem is proved.\qed
\section{Poincar\'e-Steklov operators as pseudodifferential operators}\label{SEC4}
The main purpose of this section is to define the Poincar\'e-Steklov operator $\mathscr{A}_m$ associated with the Dirac operator and to prove that it fits  into the framework of pseudodifferential operators.

 Throughout this section, let  $\O$ be a smooth domain with a compact boundary $\S$, let  $P_\pm$ be as in \eqref{Projections} and set
 \begin{align}\label{spinangular}
 \gamma_5 :=-i\alpha_1\alpha_2\alpha_3=\begin{pmatrix}
0& \mathit{I}_2 \\
\mathit{I}_2  & 0
\end{pmatrix}\quad\text{and}\quad S\cdot X= -\gamma_5 (\alpha\cdot X), \quad \forall X\in\rr^3.
\end{align}
Using the anticommutation relations of the Dirac's matrices we easily get the following identities
	\begin{align}\label{AlphaSigma}
	\begin{split}
	i (\alpha \cdot X ) (\alpha \cdot Y  ) &= i X \cdot Y + S \cdot (X \wedge Y),\\
	\left\{S \cdot X,   \alpha \cdot Y\right\} &= -(X \cdot Y) \g_5,\quad [S \cdot X,\beta]=0, \quad \forall X,Y\in\rr^3.
	\end{split}
	\end{align}

  Next, we give the rigorous definition of the Poincar\'e-Steklov operator $\mathscr{A}_m$, which is the main subject of this paper.
\begin{definition}(PS operator)\label{PSdef} Let $ z\in\rho(\Hm)$  and $g\in P_-\mathit{H}^{1/2}(\S)^4$. We denote by $E^{\O}_{m}(z): P_-\mathit{H}^{1/2}(\S)^4\rightarrow \mathit{H}^{1}(\O)^4$ the lifting operator associated with the elliptic problem
\begin{equation}\label{LL1}
\left\{
\begin{aligned}
(D_m-z)U_z&=0 \quad \text{ in } \O,\\
P_-t_{\S}U_z&= g \quad \text{ on } \S.
\end{aligned}
\right.
\end{equation}
That is,  $E^{\O}_{m}(z)g$ is the unique function in $\mathit{H}^{1}(\O)^4$ satisfying  $(D_m-z)E^{\O}_{m}(z)g=0$ in $\O$, and $P_-t_{\S}E^{\O}_{m}(z)g=g$ on $\S$. Then the Poincar\'e-Steklov (PS) operator $\mathscr{A}_m:P_-\mathit{H}^{1/2}(\partial\O)^4 \longrightarrow P_+\mathit{H}^{1/2}(\partial\O)^4$ associated with the system \eqref{LL1} is defined by 
\begin{align*}
\mathscr{A}_m (g) = P_+t_{\S}E^{\O}_{m}(z)g, 
\end{align*}
\end{definition}	

Recall the definitions of $\Phi^{\O}_{z,m}$ and $\Lambda^z_m$ from Subsection \ref{Subintope}. Then, the following proposition justifies the existence and the unicity of the solution to the elliptic problem \eqref{LL1}, and gives in particular the explicit formula of the PS operator in terms the operator $(\Lambda^z_m)^{-1}$ when $z\in\rho(D_m)$. The third assertion of the proposition will be particularly important in Section \ref{SEC5} when studying the PS operator from the semiclassical point of view. In the last statement, we use the notations $\mathscr{A}_m(z)$ to highlight the dependence on the parameter $ z\in\rho(\Hm)$. 

\begin{proposition}\label{definitionPSOP} For any $ z\in\rho(\Hm)$ and $g\in P_-\mathit{H}^{1/2}(\S)^4$, the elliptic problem \eqref{LL1} has a unique solution $E^{\O}_{m}(z)[g]\in\mathit{H}^1(\O)^4$. Moreover, the following hold true: 
\begin{itemize}
  \item[(i)] $\left(  E^{\O}_{m}(z)\right)^{\ast}=-\beta P_+t_{\S}(\Hm-\overline{z})^{-1}$.
  \item[(ii)] For  any compact set $K\subset\cc$, there is $m_0>0$ such that for all $m\geqslant m_0$  it holds that   $K\subset\rho(\Hm)$, and for all $z\in K$ we have
\begin{align*}
 \left|\left| E^{\O}_{m}(z)g \right|\right|_{\mathit{L}^2(\O)^4}&\lesssim \frac{1}{\sqrt{m}}\left|\left| g \right|\right|_{\mathit{L}^{2}(\S)^4},\quad\forall g\in P_-\mathit{H}^{1/2}(\S)^4.
 \end{align*}
  \item [(iii)] If $z\in\rho(D_m)$, then $E^{\O}_{m}(z)$ and $\mathscr{A}_{m}$ are explicitly given by
\begin{align}\label{def2PS} 
E^{\O}_{m}(z)= \Phi^{\O}_{z}(\Lambda^{z}_{m})^{-1}P_- \quad\text{and}\quad\mathscr{A}_{m}=-P_+\beta(\Lambda^z_{m})^{-1}P_{-}.
\end{align}
 \item[(iv)] Let $z\in\rho(\Hm)$ and let $E^{\O}_{m}(z)$ be as above. Then,  for any $\xi \in\rho(\Hm)$, the operator $E^{\O}_{m}(\xi)$ has the following representation  
 \begin{align}\label{properiesofPS1}
E^{\O}_{m}(\xi)= (\mathit{I}_{4}+(\xi- z)(\Hm -\xi)^{-1})  E^{\O}_{m}(z).
 \end{align}
 In particular, we have
  \begin{align}\label{properiesofPS2}
\mathscr{A}_m(\xi)-\mathscr{A}_m(z)= (z - \xi)\beta \left(  E^{\O}_{m}(\overline{\xi})\right)^{\ast} E^{\O}_{m}(z).
 \end{align}
  \item[(v)] For any  $z\in\rho(\Hm)$ the operator $E^{\O}_{m}(z)$ extends into a bounded operator from  $P_-\mathit{H}^{-1/2}(\S)^4$ to $\mathit{H}(\alpha,\O)$.
\end{itemize} 
 \end{proposition}
\textbf{Proof.}  We first show that the boundary value problem \eqref{LL1} has a unique solution. For this, assume that $u_1$ and $u_2$ are both solutions of \eqref{LL1}, then $(D_{m}-z)(u_1-u_2)=0$ in $\O$,   and $P_-t_\S (u_1-u_2)= 0$   on  $\S$.  Thus,  $(u_1-u_2)\in\mathrm{dom}(\Hm)$ holds by Remark \ref{Sobolevremar},  and since $\Hm$ is self-adjoint by Theorem \ref{MITLipshictz} it follows that $u_1=u_2$, which proves the uniqueness. Next, observe that the function 
$$v_g= \mathcal{E}_{\O}(P_-g)- (\Hm -z)^{-1}(D_m-z)\mathcal{E}_{\O}(P_-g)$$
 is a solution to \eqref{LL1}. Indeed, we have  $\mathcal{E}_{\O}(P_-g)\in\mathit{H}^1(\O)^4$ and thus $v_g\in\mathit{H}^1(\O)^4$, moreover,  we clearly have that  $P_-t_{\S}v_g=g$ and $(D_m-z)v_g=0$.  Since we already know that the solution to \eqref{LL1} is unique, it follows that  $v_g$ is independent of the extension operator $\mathcal{E}_{\O}$, and hence there is a unique solution in $\mathit{H}^1(\O)^4$ to the elliptic problem  \eqref{LL1}.

Let us show the assertion $(\mathrm{i})$. Let $\psi\in P_-\mathit{H}^{1/2}(\S)^4$ and $f\in\mathit{L}^{2}(\O)^4$, then using the Green's formula and the fact that $P_+ (-i\alpha\cdot n)= (-i\alpha\cdot n)P_-=-\beta P_-$ we get that 
 \begin{align*}
\langle E^{\O}_{m}(z)\psi, f \rangle_{\mathit{L}^{2}(\O)^4}&=\langle E^{\O}_{m}(z)\psi, (\Hm-\overline{z})(\Hm-\overline{z})^{-1}f \rangle_{\mathit{L}^{2}(\O)^4}\\
&=\langle E^{\O}_{m}(z)\psi, (D_m-\overline{z})(\Hm-\overline{z})^{-1}f \rangle_{\mathit{L}^{2}(\O)^4}\\
&=\langle (D_m-z)E^{\O}_{m}(z)\psi, (\Hm-\overline{z})^{-1}f \rangle_{\mathit{L}^{2}(\O)^4}\\
&+\langle (-i\alpha\cdot n)t_{\S}E^{\O}_{m}(z)\psi, t_{\S}(\Hm-\overline{z})^{-1}f \rangle_{\mathit{L}^{2}(\S)^4}\\
&=\langle (-i\alpha\cdot n)P_-t_{\S}E^{\O}_{m}(z)\psi, P_+t_{\S}(\Hm-\overline{z})^{-1}f \rangle_{\mathit{L}^{2}(\S)^4}\\
&=\langle \psi, -\beta P_+t_{\S}(\Hm-\overline{z})^{-1}f \rangle_{\mathit{L}^{2}(\S)^4}
\end{align*}
 which entails that $-\beta P_+t_{\S}(\Hm-\overline{z})^{-1}$ is the adjoint of $E^{\O}_{m}(z)$ and proves $(\mathrm{i})$.

Now we are going to show the assertion  $(\mathrm{ii})$. So, let $K$  be a compact set of $\cc$, and note that for all $m> \sup\{ |\mathrm{Re}(z)|: z\in K\}$ it holds that $K\subset\rho(D_{m})\subset\rho(\Hm)$. Hence, $v:=E^{\O}_{m}(z)g $ is well defined for any $z\in K$ and $g\in P_-\mathit{H}^{1/2}(\S)^4$.  Then a straightforward application of the Green's formula yields that
\begin{align}\label{inequa4}
\begin{split}
0=\left|\left| (D_{m}-z) v \right|\right|^2_{\mathit{L}^2(\O)^4}=& \left|\left| (i\alpha\cdot\nabla-z) v \right|\right|^2_{\mathit{L}^2(\O)^4}+m^2\left|\left|v \right|\right|^2_{\mathit{L}^2(\O)^4}\\
&+m\left(\langle -i(\alpha\cdot n)t_{\S}v,\beta t_{\S}v \rangle_{\mathit{L}^{2}(\S)^4}-2 \mathrm{Re}(z)\langle v,\beta v \rangle_{\mathit{L}^2(\O)^4} \right).
\end{split}
\end{align}
Observe that 
\begin{align*}
\langle -i(\alpha\cdot n)t_{\S}v,\beta t_{\S}v \rangle_{\mathit{L}^{2}(\S)^4}= \langle (P_+ -P_- )t_{\S}v, t_{\S}v \rangle_{\mathit{L}^{2}(\S)^4}= \left|\left| P_+ t_{\S}v \right|\right|^2_{\mathit{L}^{2}(\S)^4}- \left|\left| P_- t_{\S}v \right|\right|^2_{\mathit{H}^{1/2}(\S)^4}.
\end{align*}
Since $ P_- t_{\S}v=g$ and  $ P_+ t_{\S}v=\mathscr{A}_{m}(g)$ hold by definition,  and  that
$$-\mathrm{Re}(z)\langle v,\beta v \rangle_{\mathit{L}^2(\O)^4} \geqslant - |\mathrm{Re}(z)| \left|\left|v \right|\right|^2_{\mathit{L}^2(\O)^4}$$
 holds by Cauchy-Schwarz inequality, it follows from \eqref{inequa4} that
 \begin{align*}
 \left|\left| g \right|\right|^2_{\mathit{L}^{2}(\S)^4} \geqslant  m\left|\left|v \right|\right|^2_{\mathit{L}^2(\O)^4}  -2|\mathrm{Re}(z)| \left|\left|v \right|\right|^2_{\mathit{L}^2(\O)^4} +\left|\left| \mathscr{A}_{m}(g)  \right|\right|^2_{\mathit{L}^{2}(\S)^4}.
\end{align*}
Thus, if we take  $m_0\geqslant 4 \sup\{ |\mathrm{Re}(z)|: z\in K\}$, then  
 \begin{align*}
  \left|\left| \mathscr{A}_{m}(g)  \right|\right|^2_{\mathit{L}^{2}(\S)^4}+  \frac{m}{2}\left|\left|v \right|\right|^2_{\mathit{L}^2(\O)^4}\leqslant  \left|\left| g \right|\right|^2_{\mathit{L}^{2}(\S)^4}
  \end{align*}
  holds for any $m\geqslant m_0$, which prove the desired estimate for $ E^{\O}_{m}(z)$.
  
  Let us now show the assertion $(\mathrm{iii})$, so let $z\in\rho(D_m)$ and recall that $ \Phi^{\O}_{z,m}(\Lambda^{z}_{m})^{-1}: \mathit{H}^{1/2}(\S)^4 \rightarrow \mathit{H}^{1}(\O)^4$ is well defined and bounded by  Lemma \ref{prop of C} . Since  $\phi^z_{m}$ is a fundamental solution of $(D_{m}-z)$,  it holds that 
\begin{align*}
(D_{m}-z)\Phi^{\O}_{z,m}(\Lambda^{z}_{m})^{-1}[g]=0 \,\text{ in } \mathit{L}^{2}(\O)^4, \quad\forall g\in\mathit{H}^{1/2}(\S)^4.
\end{align*}
Now, observe that if  $g\in P_-\mathit{H}^{1/2}(\S)^4$, then a direct application of the identity \eqref{traceof Phi} yields that
 \begin{align*}
 t_{\S}\Phi^{\O}_{z,m}(\Lambda^{z}_{m})^{-1}[g]= \left( -\frac{i}{2}(\alpha\cdot n) + \mathscr{C}_{z,m}\right)(\Lambda^{z}_{m})^{-1}[g]= g-P_+\beta(\Lambda^{z}_{m})^{-1}[g].
\end{align*}
Consequently, we get 
 \begin{align*}
P_- t_{\S}\Phi^{\O}_{z,m}(\Lambda^{z}_{m})^{-1}[g]= g \,\,\text{ and }\,\, P_+ t_{\S}\Phi^{\O}_{z,m}(\Lambda^{z}_{m})^{-1}[g]=-P_+\beta(\Lambda^{z}_{m})^{-1}[g],
\end{align*}
which means that $ E^{\O}_{m}(z)[g]= \Phi^{\O}_{z,m}(\Lambda^{z}_{m})^{-1}[g]$ is the unique solution to the boundary value problem \eqref{LL1}, and proves the identity $\mathscr{A}_{m}=-P_{+}\beta(\Lambda^z_{m})^{-1}P_{-}$.

We are going to prove assertion $(\mathrm{iv})$, so fix $z,\xi \in\rho(\Hm)$ and let $g\in P_-\mathit{H}^{1/2}(\S)^4$.  Then, by definition of $E^{\O}_{m}(z)$ we have that 
 \begin{align*}
 (D_m -\xi)E^{\O}_{m}(\xi)g&= (D_m -\xi)(\mathit{I}_{4}+(\xi- z)(\Hm -\xi)^{-1})  E^{\O}_{m}(z)g,\\
 &= (D_m -z)E^{\O}_{m}(z)g -(\xi- z)E^{\O}_{m}(z) g +(\xi- z) (D_m -\xi)(\Hm -\xi)^{-1}  E^{\O}_{m}(z)g,\\
 &=  (\xi- z) E^{\O}_{m}(z)g -(\xi- z) E^{\O}_{m}(z)g=0. 
 \end{align*}
Since $(\Hm -\xi)^{-1} E^{\O}_{m}(z)g\in\mathrm{dom}(\Hm)$, and hence $ P_-t_{\S}(\Hm -\xi)^{-1} E^{\O}_{m}(z)g=0$, it follows that  $P_-t_{\S}E^{\O}_{m}(\xi)g=P_-t_{\S}E^{\O}_{m}(z)g=g$, which prove the identity \eqref{properiesofPS1}. Now,  \eqref{properiesofPS2} follows by applying $P_+t_{\S}$ to the representation \eqref{properiesofPS1} and using assertion $(\mathrm{i})$. 

It remains to prove item $(\mathrm{v})$. We first consider the case $z\in\rho(D_m)$, then the claim for $z\in\rho(\Hm)\setminus \rho(D_m)$ follows by the representation formula \eqref{properiesofPS1} . Fix $z\in\rho(D_m)$ and recall that the operators $\mathscr{C}_{z,m}$ and $\Lambda^{z}_{m}$ are bounded invertible in $\mathit{H}^{1/2}(\S)^4$ by Lemma \ref{prop of C}$(\mathrm{ii})$-$(\mathrm{iii})$ and \eqref{PRI}. Since  $\mathscr{C}^{\ast}_{z,m}=\mathscr{C}_{\overline{z},m}$, by duality it follows that $\Lambda^{z}_{m}$ admits a bounded and everywhere defined inverse in $\mathit{H}^{-1/2}(\S)^4$.  This together with Lemma \ref{prop of C}$(\mathrm{i})$ and item $(\mathrm{iii})$ of this proposition show that $E^{\O}_{m}(z)$ admits a continuous extension from $P_-\mathit{H}^{-1/2}(\S)^4$ to $\mathit{H}(\alpha,\O)$. This completes the proof of the proposition.
 \qed
\\ 
\begin{remark}\label{unifoBo} The proof above gives more, namely that for all $m_0>0$,  $K\subset\rho(D_{m_0})$ a compact set and $z\in K$,  there is $m_1\gg 1$ such that 
 \begin{align*}
  \sup_{m\geqslant m_1}\left|\left| \mathscr{A}_{m}  \right|\right|_{P_-\mathit{H}^{1/2}(\S)^4\longrightarrow P_+\mathit{L}^{2}(\S)^4} \lesssim 1.
    \end{align*}
\end{remark}
\begin{remark}\label{ReguLip} Thanks to Theorem \ref{MITLipshictz} and Remark \ref{remLip},  if  $\O$ is a Lipschitz domain, then  $E^{\O}_{m}(z)$ is the unique solution in $\mathit{H}^{1/2}(\O)^4$ to the system \eqref{LL1} for datum in $\mathit{L}^2(\S)^4$.  Moreover, the PS operator $\mathscr{A}_{m}=-P_{+}\beta(\Lambda^z_{m})^{-1}P_{-}$ is well-defined and bounded as an operator from $P_-\mathit{L}^{2}(\S)^4$ to $P_+\mathit{L}^{2}(\S)^4$. 
\end{remark}
	In the rest of this section, we will only address the case $z\in\rho(D_m)$ and we show that the Poincar\'e-Steklov operator $\mathscr{A}_m$ from Definition \ref{PSdef}  is a homogeneous pseudodifferential operators of order $0$ and capture its principal symbol in local coordinates. To this end, we first study the pseudodifferential properties of the Cauchy operator $\mathscr{C}_{z,m}$. Once this is done,  we use the explicit formula of $\mathscr{A}_m$ given by \eqref{def2PS} and the symbol calculus  to obtain  the principal symbol of $\mathscr{A}_m$.

 Recall the definition of $\phi^{z}_m$ from \eqref{defsolo}, and observe that 
$$\phi^{z}_m(x-y)= k^{z}(x-y)+ w(x-y),$$
 where 
\begin{align*}\label{}
k^{z}(x-y)&=\frac{e^{i\sqrt{z^2-m^2}|x-y|}}{4\pi|x-y|}\left(z +m\beta+ \sqrt{z^2-m^2}\alpha\cdot\frac{x-y}{|x-y|}\right) +i\frac{e^{i\sqrt{z^2-m^2}|x-y|}-1}{4\pi|x-y|^3}\alpha\cdot(x-y),\\
w(x-y)&=\frac{i}{4\pi|x-y|^3}\alpha\cdot(x-y).
\end{align*}  
Using this, it follows that 
\begin{align}\label{decomposeCauchy}
\begin{split}
\mathscr{C}_{z,m}[f](x)=&   \lim\limits_{\rho\searrow 0}\int_{|x-y|>\rho}w(x-y)f(y)\mathrm{d\sigma}(y)+\int_{\S}k^{z}(x-y)f(y)\mathrm{d\sigma}(y)\\
=& W[f](x)+ K[f](x).
 \end{split}
\end{align}
As $|k^{z}(x-y)|=\mathcal{O}(|x-y|^{-1})$ when $|x-y|\rightarrow0$,  using the standard layer potential techniques 
(see, e.g. \cite[Chap. 3, Sec. 4]{T2} and \cite[Chap. 7, Sec. 11]{T1}) it is not hard to prove that the integral operator $K$ gives rise to a pseudodifferential operator of order $-1$, i.e. $K\in  \mathit{Op} \sS^{-1}(\S)$. Thus,  we can (formally) write
\begin{align}\label{SymCauchy}
\mathscr{C}_{z,m}= W  \quad \mathrm{mod} \, \mathit{Op} \sS^{-1}(\S),
\end{align}
which means that the operator $W$ encodes the main contribution in the pseudodifferential character of $\mathscr{C}_{z,m}$. So  we only need to  focus on the study of the pseudodifferential  properties of $W$.  The following theorem makes this heuristic more rigorous. Its proof follows similar arguments as in \cite{AKM,Miy,MR}.

\begin{theorem}\label{Th SymClas} Let $\mathscr{C}_{z,m}$ be as  \eqref{CauchyOpe}, $W$ as in \eqref{decomposeCauchy} and  $\mathscr{A}_m$ as in  Definition \ref{PSdef}.  Then $\mathscr{C}_{z,m}$, $W$ and $\mathscr{A}_m$ are homogeneous pseudodifferential operators of order $0$, and we have 
\begin{align*}
\mathscr{C}_{z,m}&= \frac{1}{2}\alpha\cdot\frac{\nabla_\S}{\sqrt{-\Delta_\S}}\quad  \mathrm{mod} \, \mathit{Op} \sS^{-1}(\S),\\
\mathscr{A}_m&=\frac{1}{\sqrt{-\Delta_\S}}S\cdot(\nabla_\S\wedge n)P_-  \quad  \mathrm{mod} \, \mathit{Op} \sS^{-1}(\S)= 
\frac{D_\Sigma}{\sqrt{-\Delta_\S}}  P_- \quad  \mathrm{mod} \, \mathit{Op} \sS^{-1}(\S).
\end{align*}
\end{theorem}

\textbf{Proof.} We first deal with the operator $W$. So, let $\psi_k:\S\rightarrow\rr$, $k=1,2$,  be a $\mathit{C}^{\infty}$-smooth function. Clearly, if $\mathrm{supp}(\psi_2)\cap\mathrm{supp}(\psi_1)=\emptyset$, then  $\psi_2 W\psi_1$ gives rise to a bounded operator from $\mathit{H}^{-j}(\S)^4$ into $\mathit{H}^{j}(\S)^4$, for all $j\geqslant0$. \\
Now, fix a local chart $(U,V,\var)$ as in Subsection \ref{SymPseuOpSur} and  recall the definition of first fundamental form $\mathrm{I}$ and  the metric tensor $G(\tx)$. That is,  for all $x\in U$ we have $x=\varphi(\tx)=(\tx,\chi(\tx))$ with $\tx\in V$, and where  the graph of $\chi:V\rightarrow \rr$  coincides with $U$. Notice that if we assume that $\psi_k$ is  compactly supported with $\mathrm{supp}(\psi_k)\subset U$, then, in this setting, the operator  $\psi_2 W\psi_1$ has the form 
\begin{align}\label{calculW1}
\begin{split}
\psi_2 W[\psi_1f](x) =&\psi_2(x)\mathrm{p.v} \int_{V}i\alpha\cdot\frac{\varphi(\tx)-\varphi(\ty)}{4\pi|\varphi(\tx)-\varphi(\ty)|^3}\psi_1(\varphi(\ty))f(\varphi(\ty))\sqrt{g(\ty)}\mathrm{d}\ty\\
=&\psi_2(x)\sqrt{g(\tx)}\mathrm{p.v} \int_{V}i\alpha\cdot\frac{\varphi(\tx)-\varphi(\ty)}{4\pi|\varphi(\tx)-\varphi(\ty)|^3}\psi_1(\varphi(\ty)f(\varphi(\ty))\mathrm{d}\ty\\
&+ \psi_2(x)\int_{V}i\alpha\cdot\frac{\varphi(\tx)-\varphi(\ty)}{4\pi|\varphi(\tx)-\varphi(\ty)|^3}f(\varphi(\ty))\left( \sqrt{g(\ty)}-\sqrt{g(\tx)} \right) \mathrm{d}\ty,
\end{split}
 \end{align}
 where $g$ is the determinant of the metric tensor $G$. Since $g(\cdot)$ is smooth, it follows that 
 $$|\sqrt{g(\ty)}-\sqrt{g(\tx)}| \lesssim |\tx -\ty|.$$
 Therefore, the last integral operator  on the right-hand side of \eqref{calculW1} has a non singular kernel and does not require to write it as an integral operator in the principal value sense.  Thus, a simple computation using Taylor's formula shows that
\begin{align*}
|x-y |^2&=|\varphi(\tx)-\varphi(\ty)|^2
= \langle \tx-\ty, G(\tx)(\tx-\ty)\rangle( 1+\mathcal{O}|\tx-\ty|),
\end{align*}
where the definition of $\mathrm{I}$  was used in the last equality.  It follows from the above computations that 
\begin{align*}
|x-y |^{-3}&= \frac{1}{ \langle \tx-\ty, G(\tx)(\tx-\ty)\rangle^{3/2}} + k_1(\tx,\ty),
\end{align*}
where the kernel $k_1$ satisfies $|k_1(\tx,\ty)|=\mathcal{O}(|\tx-\ty|^{-2})$, when $|\tx-\ty|\rightarrow0$. Consequently,  we get that  
\begin{equation*}\label{} 
  \frac{x_j-y_j}{|x-y |^{3}}= \left\{
  \begin{aligned}
&\frac{\tx_j-\ty_j}{\langle \tx-\ty, G(\tx)(\tx-\ty)\rangle^{3/2}}+(\tx_j-\ty_j)k_1(\tx,\ty),\quad \text{for }j=1,2,\\
&\frac{\langle \nabla\chi,\tx-\ty\rangle}{\langle \tx-\ty, G(\tx)(\tx-\ty)\rangle^{3/2}} + k_2(\tx,\ty) ,\quad \text{for }j=3,
\end{aligned}
  \right.
\end{equation*}
with $|k_2(\tx,\ty)|=\mathcal{O}(|\tx-\ty|^{-1})$, when $|\tx-\ty|\rightarrow0$. Note that this implies
\begin{align*}
\alpha\cdot\left(\frac{x-y}{|x-y |^{3}}\right)=\alpha\cdot\frac{(\tx-\ty,\langle \nabla\chi,\tx-\ty\rangle)}{\langle \tx-\ty, G(\tx)(\tx-\ty)\rangle^{3/2}} + \mathcal{O}(|\tx-\ty|^{-1})\mathit{I}_4.
\end{align*}
Combining the above computations and \eqref{calculW1}, we deduce that 
 \begin{align}\label{calculW}
\begin{split}
\psi_2 W[\psi_1f](x)=\psi_2(x)\sqrt{g(\tx)}\mathrm{ p.v} \int_{V}i\alpha\frac{(\tx-\ty,\langle \nabla\chi,\tx-\ty\rangle)}{\langle \tx-\ty, G(\tx)(\tx-\ty)\rangle^{3/2}}f(\varphi(\ty))\mathrm{d}\ty + \psi_2(x)L[\psi_1f](x),
\end{split}
 \end{align}
 where  $L$ is an integral operator with a kernel $l(x,y)$  satisfying 
 \begin{align*}
 |l(x,y)|=\mathcal{O}(|x-y|^{-1}) \,\, \text{ when }\, |x-y|\rightarrow0.
\end{align*}
Thus,  similar arguments as the ones in \cite[Chap. 7, Sec. 11]{T1}  yield that $L$ is a pseudodifferential operator of order $-1$. Now, for $h\in\mathit{L}^2(\rr^2)$ and $k=1,2$, observe that if we set
\begin{align*}
R_k[h](\tx)= \frac{ i\sqrt{g(\tx)}}{4\pi}\int_{\rr^2}r_k(\tx,\tx-\ty)h(\ty)\mathrm{d}(\ty),
\end{align*}
where 
\begin{align*}
 r_k(\tx,\tx-\ty)=\frac{\tx_k-\ty_k}{ \langle \tx-\ty, G(\tx)(\tx-\ty)\rangle^{3/2}}, \quad \tx\neq\ty.
\end{align*}
Then the standard formula connecting a pseudodifferential operator and its symbol yields
\begin{align*}
R_k[h](\tx)=\frac{ i\sqrt{g(\tx)}}{(2\pi)^2}\int_{\rr^2}\int_{\rr^2}e^{i\langle \tx-\ty,\xi\rangle}q_k(\tx,\xi)h(\ty)\mathrm{d}\xi\mathrm{d}\ty,
\end{align*}
where 
\begin{align*}
q_k(\tx,\xi)= \frac{i\sqrt{g(\tx)}}{2}\int_{\rr^2}e^{-i\langle \omega,\xi\rangle}r_k(\tx,\omega)\mathrm{d}\xi.
\end{align*}
Recall the definition of $Q$ from \eqref{defPmatrix} and set $\omega=Q(\tx)\tau$. Also recall that 
\begin{align}\label{FourierdeX/X}
 \int_{\rr^2}e^{-i\langle \omega,\xi\rangle}\frac{\omega_k}{|\omega|^3}\mathrm{d}\omega= -i\frac{\xi_k}{|\xi|},\quad k=1,2.
\end{align}
Thus, the above change of variables together with the properties  \eqref{lesrelations} and \eqref{FourierdeX/X} yield that 
\begin{align*}
q_k(\tx,\xi)&=  \frac{i}{2}\int_{\rr^2}e^{-i\langle \tau,Q^{\mathrm{t}}(\tx)\xi\rangle}\frac{(Q^{\mathrm{t}}(\tx)\tau)_k}{|\tau|^3}\mathrm{d}\tau=\frac{(G^{-1}(\tx)\xi)_k}{2\langle G^{-1}(\tx) \xi,\xi\rangle^{1/2}}=\frac{g_{k1}\xi_1+g_{k2}\xi_2}{2\langle G^{-1}(\tx) \xi,\xi\rangle^{1/2}},
\end{align*}
which means that $q_k(\tx,\xi)$  is  homogeneous of degree $0$ in $\xi$. Therefore,  $R_k$ is a homogeneous pseudodifferential operators of degree $0$. 
From the above observation and \eqref{calculW} if follows that 
\begin{align*}
\psi_2 W\psi_1= \psi_2\alpha\cdot\left(R_1,R_2,  \partial_1\chi(\tx)R_1+ \partial_2\chi(\tx)R_2\right)\psi_1+ \psi_2L\psi_1.
\end{align*}
Since $L$ is a pseudodifferential operator of order -1, we deduce that $W$ is a homogeneous pseudodifferential operators of order $0$, and exploiting (\ref{symbolLG}), we obtain that
\begin{align}\label{SymW}
W&= \frac{1}{2}\alpha\cdot\frac{\nabla_\S}{\sqrt{-\Delta_\S}}\quad\mathrm{mod}\, \mathit{Op} \sS^{-1}(\S).
\end{align}
Thanks to \eqref{SymCauchy} and  \eqref{SymW},  we deduce that the Cauchy operator $ \mathscr{C}_{z,m}$ has  the same principal symbol as the operator $ W$. 
 
 Now we are going to deal with the operator $\mathscr{A}_m$. Note that  we have 
\begin{align}\label{formule de A_m}
 \frac{1}{2}\left(\beta+  \alpha\cdot\frac{\nabla_\S}{\sqrt{-\Delta_\S}}\right)^2=\mathit{I}_4,
\end{align}
and as $\mathscr{A}_m$ is given by the formula 
\begin{align*}
\mathscr{A}_m= -P_+\beta\left(\frac{1}{2}\beta+ \mathscr{C}_{z,m}\right)^{-1}P_-,
\end{align*}
using \eqref{formule de A_m} and the standard mollification arguments, it follows from  the product formula for calculus of pseudodifferential operators  that,  in local coordinates, the symbol of $\mathscr{A}_m$ denoted by $q_{\mathscr{A}_m}$ has the form
\begin{align*}
q_{\mathscr{A}_m}(\tx,\xi)= -P_+\beta\left(\beta+\alpha\cdot\left(\frac{\xi_G}{ \langle G^{-1}\xi,\xi\rangle^{1/2}}\right)
\right)P_- +p(\tx,\xi),
\end{align*}
where $p \in\, \sS^{-1}(\S)$ and $\xi_G$ defined in (\ref{symbolLG}) is the principal symbol of $\nabla_\Sigma$. Therefore, we get 
\begin{align*}
q_{\mathscr{A}_m}(\tx,\xi)= -P_+ \, \beta \, \alpha\cdot \xi_G  \,  \langle G^{-1}\xi,\xi\rangle^{-1/2} \, P_- +p(\tx,\xi).
\end{align*}
Hence, using the fact that $P_\pm$ are projectors, a simple computation shows 
\begin{align*}
q_{\mathscr{A}_m}(\tx,\xi)= -  i\alpha\cdot n^\varphi(\tx) \, \alpha\cdot \xi_G  \,  \langle G^{-1}\xi,\xi\rangle^{-1/2} \,  P_- +p(\tx,\xi).
\end{align*}
Finally, from results of Section \ref{sDirac} we deduce
\begin{align*}
q_{\mathscr{A}_m}(\tx,\xi)=- \alpha \cdot n^\varphi(\tx) \, \alpha\cdot \xi_G  \,  \langle G^{-1}\xi,\xi\rangle^{-1/2} P_- + p(\tx,\xi) = 
S\cdot\left(\frac{\xi_G \wedge n^\varphi(\tx) }{ \langle G^{-1}\xi,\xi\rangle}\right)P_- + p(\tx,\xi),
\end{align*}
and 
\begin{align*}
\mathscr{A}_m= \frac{D_\Sigma }{\sqrt{-\Delta_\S}} P_-  \quad  \mathrm{mod} \, \mathit{Op} \sS^{-1}(\S)=\frac{1}{\sqrt{-\Delta_\S}}S\cdot(\nabla_\S\wedge n )P_- \quad\mathrm{mod}\, \mathit{Op} \sS^{-1}(\S).
\end{align*}
It justifies that $ \mathscr{A}_m$ is a homogeneous pseudodifferential operators of order $0$ and completes the proof of the theorem. 
\qed

\section{Approximation of the Poincar\'e-Steklov operators for large masses}\label{SEC5}

Although the technique used in the last section allows us to treat the layer potential operator $\mathscr{A}_m$ as pseudodifferential operator and to derive its principal symbol. However, it does not allow us to capture the dependence on $m$. The main goal of this section is to study  the Poincar\'e-Steklov operator, $\mathscr{A}_m$,  as a  $m$-dependent pseudodifferential operator when $m$ is large enough. For this purpose, we consider  $h=1/m$ as a semiclassical parameter (for $m\gg1$) and use the system \eqref{LL1} instead of the layer potential formula of $\mathscr{A}_m$.  Roughly speaking, we will look for a local approximate formula for the solution of \eqref{LL1}. Once this is done, we use the regularization property of the resolvent of the MIT bag operator to catch the semiclassical principal symbol of $\mathscr{A}_m$.

Throughout this section, we assume that $m>1$, $z\in\rho(\Hm)$ and that $\O$ is smooth with a compact boundary $\S:=\partial\O$. Next, we introduce the semiclassical parameter $h=m^{-1}\in(0,1]$, and we set $\mathscr{A}^h:=\mathscr{A}_m$. 
Then, the following theorem is the main result of this section, it ensures that $\mathscr{A}^h$ is a $h$-pseudodifferential operator of order $0$ and gives its semiclassical principal symbol.
\begin{theorem}\label{pseudo de Am} 
Let $h\in(0,1]$ and $z\in\rho(\Hm)$, and let $\mathscr{A}^{h}$ be as above. 
Then for any $N \in \nn$, there exists a $h$-pseudodifferential operator of order $0$, $\mathscr{A}_N^{h} \in \mathit{Op}^h \sS^{0}(\S)$ such that for $h$ sufficiently small, and any $0 \leq l \leq N + \frac12$
\begin{equation*}
\| \mathscr{A}^{h} -  \mathscr{A}_N^{h}  \|_{H^{\frac12}(\S) \to H^{N+\frac32-l}(\S)}  =  O(h^{N+ \frac12 + l}),
 \end{equation*} 	
and 
 \begin{equation*}
 \mathscr{A}_N^{h} = \frac{h D_\Sigma}{\sqrt{-h^2\Delta_{\S}+I}+I}\, P_- \quad\mathrm{mod}\,  \,  h \mathit{Op}^h \sS^{-1}(\S).
 \end{equation*}
\end{theorem}

Let us consider $\mathbb{A}=\{ (U_{\varphi_j}, V_{\varphi_ j},\var_j): j\in\{1,\cdots, N\} \}$ an atlas of $\S$ and $(U_{\varphi}, V_\varphi,\var) \in \mathbb{A}$.  As in Section \ref{SEC3} we consider the case where $U_{\varphi}$ is the graph of a smooth function $\chi$, and we assume that $\O$ corresponds locally to the side $x_3 > \chi(x_1,x_2)$. Then,  for
\begin{align*} 
U_{\varphi} =  & \{ (x_1,x_2, \chi(x_1,x_2)) ; \, (x_1,x_2) \in V_\varphi \}; \quad \varphi((x_1,x_2, \chi(x_1,x_2)) = (x_1,x_2) \\
\mathcal{V}_{\varphi,\varepsilon}:= & \{ (y_1,y_2, y_3 + \chi(y_1,y_2)) ; \, (y_1,y_2,y_3) \in V_\varphi \times (0, \varepsilon) \} \subset \O,
	\end{align*} 
with $\varepsilon$ sufficiently small, we have the following homeomorphism:
\begin{align*} 
\phi :  \mathcal{V}_{\varphi,\varepsilon} & \longrightarrow  V_\varphi \times (0, \varepsilon)\\
 (x_1,x_2, x_3)  & \mapsto (x_1,x_2, x_3 - \chi(x_1,x_2)).
	\end{align*} 
Then  the  pull-back 
\begin{align*} 
\phi^*:  C^\infty(V_\varphi \times (0, \varepsilon)) & \longrightarrow   C^\infty(\mathcal{V}_{\varphi,\varepsilon} ) \\
v  & \mapsto \phi^*v:= v  \circ \phi 
	\end{align*} 
transforms the differential operator $D_m$ restricted on $\mathcal{V}_{\varphi, \varepsilon}$ into the following operator on $V_\varphi \times (0, \varepsilon)$:
\begin{align*} 
\widetilde{D}^\varphi_m:=  (\phi^{-1} )^* D_m (\phi )^* & = 
 -i\left( \alpha_1 \partial_{y_1} +   \alpha_2 \partial_{y_2} - ( \alpha_1 \partial_{x_1} \chi  +\alpha_2 \partial_{x_2} \chi - \alpha_3)\partial_{y_3}\right) + m \beta \\
 & =  -i( \alpha_1 \partial_{y_1}  +   \alpha_2 \partial_{y_2}) + \sqrt{1+ | \nabla \chi |^2}(i \alpha \cdot n^\varphi)(\ty )\partial_{y_3}  + m \beta,
 \end{align*} 
where $\ty=(y_1,y_2)$ and $n^\varphi=(\varphi^{-1})^*n$ is the pull-back of the outward pointing normal to $\O$ restricted on $V_\varphi$:
\begin{equation*}
n^\varphi(\ty) = \frac{1}{\sqrt{1+ | \nabla \chi |^2}}
 \begin{pmatrix}
\partial_{x_1} \chi  \\   \partial_{x_2} \chi   \\ - 1 
\end{pmatrix} (y_1,y_2).
\end{equation*}
	\begin{figure}[h]
\includegraphics[scale=0.2]{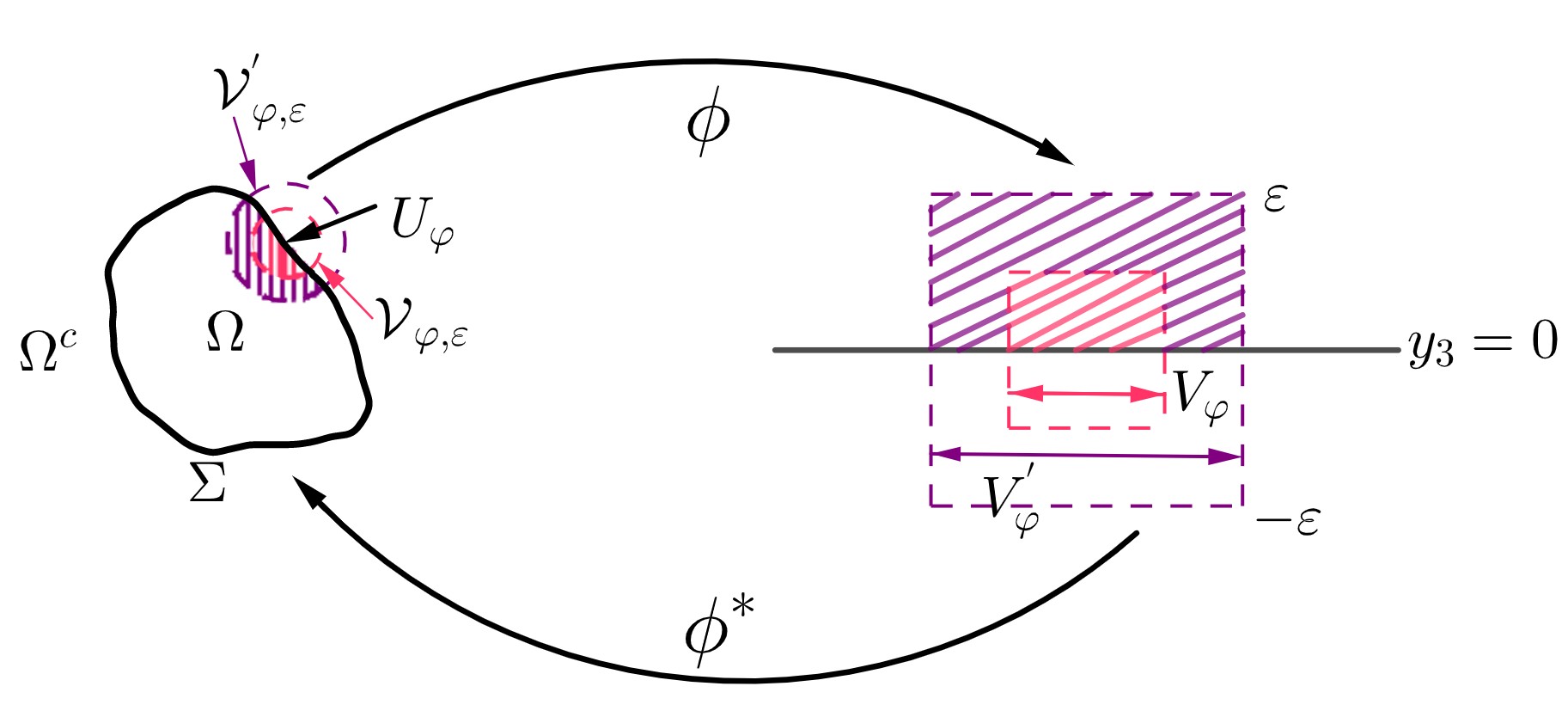} 
    \label{figure}
    \caption{Change of coordinates}
    \end{figure}
For the projectors $P_\pm$, we have:
\begin{align*}P_\pm^\varphi : =(\varphi^{-1})^*P_\pm (\varphi)^*= \frac{1}{2}\Big(\mathit{I}_4\mp i\beta \, \alpha \cdot n^\varphi(\ty)\Big).
\end{align*}
Thus, in the variable $y \in V_\varphi \times (0, \varepsilon)$, the equation \eqref{LL1}  
becomes:
\begin{equation}
	\left\{
\begin{aligned}\label{T11}	
(\widetilde{D}^\varphi_m  -z )u	=0, & \quad  \text{ in } V_\varphi \times (0, \varepsilon) ,\\
	\G^\varphi_- u= g^\varphi= g \circ \varphi^{-1}, & \quad    \text{ on } V_\varphi \times \{0\}, 
	\end{aligned}
	\right.
	\end{equation}
where $ \G^\varphi_\pm = P^\varphi_{\pm}t_{\{ y_3 = 0 \}}$. 

By isolating the derivative with respect to $y_3$, and using that $(i\alpha \cdot n^\varphi)^{-1} = - i \alpha \cdot n^\varphi$,  the system \eqref{T11} becomes:
\begin{equation}
	\left\{
\begin{aligned}\label{T1b}	
\partial_{y_3} u	= \frac{i \alpha \cdot n^\varphi (\ty)}{\sqrt{1+ | \nabla \chi (\ty) |^2}} 
	\Big( -i \alpha_1 \partial_{y_1}  - i \alpha_2 \partial_{y_2}   + m \beta -z \Big) u
 , & \quad  \text{ in } V_\varphi \times (0, \varepsilon) ,\\
	\G^\varphi_- u= g^\varphi , & \quad    \text{ on } V_\varphi \times \{0\}.
	\end{aligned}
	\right.
	\end{equation}
Let us now introduce the matrices-valued symbols 
\begin{equation}\label{defL0L1}
L_0(\ty, \xi) := \frac{i \alpha \cdot n^\varphi (\ty)}{\sqrt{1+ | \nabla \chi (\ty) |^2}} 
	\Big( \alpha \cdot \xi    +  \beta \Big) ; \qquad L_1(\ty):= \frac{-i z \alpha \cdot n^\varphi (\ty)}{\sqrt{1+ | \nabla \chi (\ty) |^2}} ,
\end{equation}
with $\xi = (\xi_1, \xi_2)$ identified with $(\xi_1, \xi_2,0)$.
Then \eqref{T1b} is equivalent to
\begin{equation}
	\left\{
	\begin{aligned}\label{T1bh}	
h \partial_{y_3} u	= L_0(\ty, h D_{\ty})  u + h L_1(\ty)  u 
 , & \quad  \text{ in } V_\varphi \times (0, \varepsilon) ,\\
	\G^\varphi_- u= g^\varphi , & \quad    \text{ on } V_\varphi \times \{0\}.
	\end{aligned}
	\right.
	\end{equation}
%
%
Before constructing an approximate solution of the system \eqref{T1bh}, let us give some properties of $L_0$.
\subsection{Algebric properties of $L_0$}
The following lemma will be used in the sequel, it gathers some useful properties which allow us to simplify the expression of $L_0(\ty,\xi)$. We omit the proof since it is an easy consequence of the anticommutation relations of the Dirac's matrices and the formulas \eqref{AlphaSigma}.

	\begin{lemma}\label{properties3.1} Let  $n^\varphi$ and $\xi$ be as above, and let $S$ be as in \eqref{spinangular}. Then, for any  $z \in \cc$ and any $\tau\in\rr^3$ such that  $\tau \perp n^\varphi$, the following identities hold:
	\begin{align*}
	 \left(S \cdot \tau   - i m \beta (\alpha \cdot n^\varphi (\ty))  \right)^2 = \left(|\tau|^2  + m^2  \right) I_4.
	 \end{align*}
	 \begin{align*}
 P_\pm^{\varphi} (S\cdot\tau)=(S\cdot\tau)P^{\varphi}_\mp\quad\text{ and }\quad P^{\varphi}_\pm(i\alpha\cdot n^{\varphi})= (i\alpha\cdot n^{\varphi})P^{\varphi}_\mp .
	  \end{align*}
		\end{lemma}
		The next proposition gathers the main properties of the operator $L_0$.
	\begin{proposition}\label{expo de q} Let $L_0(\ty,\xi)$ be as in \eqref{defL0L1}, then we have 
        \begin{align*}
	L_0(\ty,\xi) &= \frac{1 }{\sqrt{1+ | \nabla \chi (\ty)|^2}} 
	\Big( i \xi \cdot n^\varphi(\ty)  + S \cdot (n^\varphi(\ty)  \wedge \xi)   - i  \beta (\alpha \cdot n^\varphi(\ty) )    \Big),\\
	 &=i \xi \cdot \tn^\varphi (\ty) + \frac{\lambda (\ty,\xi) }{\sqrt{1+ | \nabla \chi (\ty)|^2}}\Pi_+ (\ty,\xi) - \frac{\lambda (\ty,\xi) }{\sqrt{1+ | \nabla \chi (\ty)|^2}}\Pi_-  (\ty,\xi)
	\end{align*}
	where 
	\begin{align}\label{deflambdapro}
	\begin{split}
	\lambda(\ty, \xi):=&\sqrt{ |{n}^\varphi(\ty) \wedge\xi|^2  + 1 }= \sqrt{ \langle G(\ty)^{-1}\xi,\xi\rangle  + 1 } ,\\
	\tn^{\varphi}(\ty) :=& \frac{1 }{\sqrt{1+ | \nabla \chi |^2}} n^{\varphi}(\ty),\\
	 \Pi_\pm  ( \ty,\xi) :=& \frac{1}2 \left(  I_4 \pm    \frac{S \cdot (n^\varphi(\ty) \wedge \xi)   - i \beta (\alpha \cdot n^\varphi (\ty)) }{\lambda(\ty, \xi)}  \right),
	 \end{split}
	\end{align} 
	with $G$ the induced metric defined in Section \ref{SymPseuOpSur}. 
	
	In particular, the symbol $L_0(\ty, \xi)$ is elliptic in ${\sS}^1$ and it admits two eigenvalues $\rho_\pm (\cdot, \cdot) \in {\sS}^1$ of multiplicity $2$ which are given by 
	\begin{align}\label{eigenlam}
	\rho_{\pm}(\ty,\xi)= \frac{i   n^\varphi(\ty)\cdot \xi \pm  \lambda(\ty,\xi) }{\sqrt{1+ | \nabla \chi |^2}} ,
	\end{align}
	and for which there exists $c>0$ such that 
	\begin{equation}\label{ellip}
\pm \Re \rho_\pm (\ty, \xi) > c \langle \xi \rangle,
\end{equation}
uniformly with respect to $\ty$. Moreover, $\Pi_\pm(\ty , \xi)$ are the projections onto $\mathrm{Kr}(L_0(\ty,\xi) -  \rho_\pm (\ty, \xi)  I_4 )$, belong to the symbol class ${\sS}^0$ and satisfy:
\begin{align}\label{PPiP}
\begin{split}
P^\varphi_\pm \,  \Pi_\pm ( \ty,\xi)\,  P^\varphi_\pm = k_{+}^\varphi ( \ty,\xi)  P^\varphi_\pm \quad\text{and }\, P^\varphi_\pm \,  \Pi_\mp ( \ty,\xi)\,  P^\varphi_\mp=\mp\Theta^{\varphi}(\ty, \xi)P^\varphi_\mp,
\end{split}
\end{align}
	with 
\begin{align} \label{P-Pi-2}
\begin{split}
k_\pm^\varphi ( \ty,\xi) = \frac{1}{2}\left(1\pm \frac{1}{\lambda(\ty, \xi)}\right)  ,\quad \Theta^{\varphi}(\ty, \xi)=\frac{1}{2\lambda(\ty, \xi)}\left(S\cdot (n^\varphi(\ty)\wedge\xi) \right).
\end{split}
\end{align}
That is, $k^{\varphi}_+$ is a positive function of ${\sS}^0$,  $(k^{\varphi}_+)^{-1} \in {\sS}^0$  and  $\Theta^{\varphi}\in {\sS}^0$. 
\end{proposition}	
\begin{remark}\label{isoPPi}
Thanks to the property \eqref{PPiP} a $4 \times 4$-matrix $A$ is uniquely determined by $P^{\varphi}_- A$ and $\Pi_+ A$ and we have:
$$A = P^{\varphi}_- A + P^{\varphi}_+ A= P^{\varphi}_- A + \frac1{k^{\varphi}_+} P^{\varphi}_+ \Pi_+ P^{\varphi}_+ A = \Big( I - \frac{P^{\varphi}_+ \Pi_+}{k^{\varphi}_+}\Big) P^{\varphi}_- A + \frac{P^{\varphi}_+}{k^{\varphi}_+}  \Pi_+  A.$$
\end{remark}

	\textbf{Proof of Proposition \ref{expo de q}.} By definition it is clear that $L_0(\ty, \xi)$ belongs to the symbol class ${\sS}^1$, $\Pi_\pm(\ty , \xi), \Theta^{\varphi} \in {\sS}^0$,   $k^{\varphi}_+$ a positive function of ${\sS}^0$ and $(k^{\varphi}_+)^{-1} \in {\sS}^0$. Now, by \eqref{AlphaSigma} we obtain that 
	  \begin{align*}
	L_0(\ty,\xi) &= \frac{1 }{\sqrt{1+ | \nabla \chi (\ty)|^2}} 
	\Big( i \xi \cdot n^\varphi(\ty)  + S \cdot (n^\varphi(\ty)  \wedge \xi)   - i  \beta (\alpha \cdot n^\varphi(\ty) )    \Big),
	\end{align*}
and since $(n^\varphi  \wedge \xi)\perp n^\varphi$, Lemma \ref{properties3.1} yields that  
	$$ \left(  S \cdot (n^\varphi(\ty)  \wedge \xi)   - i  \beta (\alpha \cdot n^\varphi(\ty) )   \right)^2=  |{n}^\varphi\wedge\xi|^2  + 1 = (\lambda(\ty, \xi))^2, $$ 
	with  $\lambda$ as  in \eqref{deflambdapro}. From this we deduce that $L_0(\ty,\xi)$ has two eigenvalues $\rho_\pm $ which are given by \eqref{eigenlam} and $\Pi_\pm(\ty , \xi)$ are the corresponding projectors onto $\mathrm{Kr}(L_0(\ty,\xi) -  \rho_\pm (\ty, \xi)  I_4 )$. 
Next, using  (\ref{symbolLD})
we get for some  $c>0$ independent of $\ty$ that 
	\begin{equation*}
\pm \Re \rho_\pm (\ty, \xi)=  \frac{\sqrt{ |{n}^\varphi\wedge\xi|^2  + 1 } }{\sqrt{1+ | \nabla \chi |^2}} =  \frac{\sqrt{ \langle G(\ty)^{-1}\xi,\xi\rangle  + 1}  }{\sqrt{1+ | \nabla \chi |^2}}   \geqslant c (1+ |\xi|),
\end{equation*}
which gives \eqref{ellip} and shows that $ \rho_\pm$ are elliptic in ${\sS}^1$. Consequently, we also get that $L_0(\ty, \xi)$ is elliptic in ${\sS}^1$.  Now,  using Lemma \ref{properties3.1} and  the properties  \eqref{spinangular},  a simple computation shows that
\begin{align*}
	P^{\varphi}_+   \Pi_{\pm} &= k_\pm^\varphi P^{\varphi}_+ \pm\frac{1}{2\lambda}\left(S\cdot (n^\varphi(\ty)\wedge\xi) \right) \,  P^{\varphi}_-,\\
	P^{\varphi}_-  \Pi_\pm&= k_\mp^\varphi P^{\varphi}_- \pm\frac{1}{2\lambda}\left(S\cdot (n^\varphi(\ty)\wedge\xi)\right) \,  P^{\varphi}_+,
\end{align*}
with $k_\pm^\varphi $ as in \eqref{P-Pi-2}. Hence, \eqref{PPiP} directly follows from the above formulas.\qed\\

\subsection{Semiclassical parametrix for the boundary problem}\label{ssSclBPb}
In this section,  we construct the approximate solution of the system \eqref{Dj} mentioned in the introduction. For simplicity of notation, in the sequel  we will use $y$ and $P_{\pm}$ instead of $\ty$ and $P_{\pm}^{\varphi}$, respectively.

We  are going to construct a local approximate solution of the following first order system:
\begin{equation*}
	\left\{
	\begin{aligned}
	h \partial_{\tau} u^h & = L_0(y, h D_y)  u^h + h L_1(y) u^h ,  &  \quad \text{ in }  &  \rr^2 \times (0, + \infty),\\
	P_- u^h  & = f,  &  \text{ on } &  \rr^2 \times \{0\}.
	\end{aligned}
	\right.
	\end{equation*}
To be precise, we will look for a solution $u^h$  in the following form:
\begin{equation}\label{uh}
u^h(y,\tau)= 
\mathit{Op}^h (A^h(\cdot, \cdot, \tau) ) f = 
 \int_{\rr^2} A^h(y, h \xi , \tau) e^{i y \cdot \xi} \hat{f}(\xi) \text{d} \xi ,
\end{equation}
with $A^h(\cdot , \cdot , \tau) \in {\sS}^0$ for any $\tau>0$  constructed inductively in the form:
$$A^h(y,  \xi , \tau) \sim \sum_{j \geq 0} h^j A_j(y, \xi , \tau).$$

%

The action of $h \partial_{\tau} - L_0(y, h D_y)  -  h L_1(y) $ on $A^h(y, h D_y, \tau) f$ is given by $T^h(y, h D_y, \tau) f$, with
$$T^h(y, \xi , \tau) = h (\partial_\tau A) (y, \xi , \tau) - L_0(y,\xi) A (y, \xi , \tau) + h \Big(  L_1(y) A (y, \xi , \tau) - \partial_\xi L_0(y, \xi) \cdot \partial_y  A (y, \xi , \tau) \Big).$$
Then we look for $A_0$ satisfying:
\begin{equation}\label{SystS0}
	\left\{
	\begin{aligned}
	h \partial_{\tau} A_0 (y, \xi , \tau)  & =   L_0(y, \xi) A_0 (y, \xi , \tau) , \\
	P_- (y) A_0 (y, \xi , \tau)   & =   P_- (y) ,
	\end{aligned}
	\right.
	\end{equation}
and for $j \geq 1$,
\begin{equation}\label{SystSj}
	\left\{
	\begin{aligned}
	h \partial_{\tau} A_j (y, \xi , \tau)  & =   L_0(y, \xi) A_j (y, \xi , \tau) +  L_1(y) A_{j-1} (y, \xi , \tau) - \partial_\xi L_0(y, \xi) \cdot \partial_y  A_{j-1} (y, \xi , \tau), \\
	P_- (y) A_j (y, \xi , \tau)   & =  0 ,
	\end{aligned}
	\right.
	\end{equation}
	
	Let us introduce a class of parametrized symbols, in which we will construct the family $A_j $:
	
\begin{equation*}
{\cP}_{h}^{m}:= \{ b(\cdot, \cdot, \tau ) \in {\sS}^{m} ; \; \forall (k,l) \in \nn ^2 , \, \tau^k \partial_\tau^l b(\cdot, \cdot, \tau ) \in h^{k-l} {\sS}^{m-k+l}  \}; \quad m \in \zz.
\end{equation*} 	

\begin{proposition}\label{SolSystS0}
There exists $A_0 \in {\cP}_{h}^{0}$ solution of \eqref{SystS0} given by:
$$A_0 (y, \xi , \tau) =  \frac{\Pi_-(y,\xi) P_-(y)}{k^{\varphi}_+(y,\xi)} \, e^{h^{-1} \tau \rho_-(y,\xi)}.$$

\end{proposition}
\textbf{Proof.} 
The solutions of the differential system $h \partial_{\tau} A_0  =   L_0  A_0$ are $A_0 (y, \xi , \tau) =   e^{h^{-1} \tau L_0(y,\xi)} A_0 (y, \xi , 0)$. By definition of $\rho_\pm$ and $\Pi_\pm$, we have:
\begin{equation}\label{expL0}
e^{h^{-1} \tau L_0(y,\xi)}  = e^{h^{-1} \tau \rho_-(y,\xi)} \Pi_-(y,\xi) + e^{h^{-1} \tau \rho_+(y,\xi)} \Pi_+(y,\xi).
\end{equation} 
It follows from \eqref{ellip} that $A_0 $ belongs to ${\sS}^0$  for any $\tau >0$ if and only if $ \Pi_+(y,\xi) A_0 (y, \xi , 0) =0$. Moreover, the boundary condition $P_- A_0  =   P_- $ implies $P_- (y) A_0 (y, \xi , 0)  =   P_- (y)$.  Thus, thanks to Remark \ref{isoPPi}, we deduce that 
$$A_0 (y, \xi , 0) =  P_-(y) - \frac{P_+ \Pi_+ P_-}{k^{\varphi}_+} (y,\xi) = P_-(y) + \frac{P_+ \Pi_- P_-}{k^{\varphi}_+} (y,\xi) = \frac{\Pi_- P_-}{k^{\varphi}_+} (y,\xi).$$
The properties of $\rho_-$, $\Pi_-$, $P_-$ and $k_+$ given in Proposition \ref{expo de q}, imply that $(k^{\varphi}_+)^{-1}\Pi_- P_- \in \sS^0$ and that $e^{h^{-1} \tau \rho_-(y,\xi)} \in \cP^0_h$. This concludes the proof of Proposition \ref{SolSystS0}.\qed\\

For the other terms $A_j$, $j \geq 1$, we have:

\begin{proposition}\label{SolSystSj}
Let $A_0$ be defined by Proposition \ref{SolSystS0}. Then for any $j\geq 1$, there exists $A_j \in h^j {\cP}_{h}^{-j}$ solution of \eqref{SystSj} which has the form:
\begin{equation}\label{eqAj}
 A_j (y, \xi , \tau) =  e^{h^{-1} \tau \rho_-(y,\xi)} \, \sum_{k=0}^{2j} (h^{-1} \tau  \langle \xi \rangle   )^k B_{j,k} (y, \xi),
 \end{equation} 	
with $B_{j,k}  \in h^j \sS^{-j}$.
\end{proposition}
\textbf{Proof.} 
Since $A_0$ has already the claimed form by Proposition \ref{SolSystS0},  so for  $A_j$ with $j\geq 1$, it is sufficient to prove the induction step. Thus,  assume that there exists  $A_j \in h^j {\cP}_{h}^{-j}$ solution of \eqref{SystSj} satisfying the above property and let us prove that the same holds for $A_{j+1}$. In order to be a solution of the differential system $h \partial_{\tau} A_{j+1}  =   L_0  A_{j+1} +L_1 A_j - \partial_\xi L_0 \cdot \partial_y  A_{j}$, for $A_{j+1}$ we have:
 \begin{equation}\label{Agene}
 A_{j+1}  =   e^{h^{-1} \tau L_0} {A_{j+1}}_{| \tau=0}   + e^{h^{-1} \tau L_0} \int_0^\tau e^{-h^{-1} s L_0} ( L_1 A_j - \partial_\xi L_0 \cdot \partial_y  A_{j} )\text{d}s,
 \end{equation} 	
 where $L_1 A_j$ has still the form \eqref{eqAj}, and we have 
 $$  \partial_y  A_{j}  =  e^{h^{-1} \tau \rho_-} \, \Big( h^{-1} \tau \partial_y \rho_- + \partial_y \Big) \sum_{k=0}^{2j} (h^{-1} \tau  \langle \xi \rangle   )^k B_{j,k}. $$
 Thus, thanks to the  properties $ \rho_-$ and $B_{j,k}$, the quantity $( L_1 A_j - \partial_\xi L_0 \cdot \partial_y  A_{j} )(y, \xi, s)$ has the form:
 \begin{equation}\label{Agene1}
  e^{h^{-1} s \rho_-(y,\xi)} \, \sum_{k=0}^{2j+1} (h^{-1} \tau  \langle \xi \rangle   )^k \tilde{B}_{j,k} (y, \xi) 
  \end{equation} 	
 with $\tilde{B}_{j,k} \in h^j \sS^{-j}$.
 So, by using the decomposition \eqref{expL0}, for the second term of the r.h.s. of \eqref{Agene} we have:
 \begin{equation}\label{2ndAgene}
 e^{h^{-1} \tau L_0} \int_0^\tau e^{-h^{-1} s L_0} ( L_1 A_j - \partial_\xi L_0 \cdot \partial_y  A_{j} )\text{d}s
 = e^{h^{-1} \tau \rho_-}  \Pi_-  I_-^j (\tau) 
 +
  e^{h^{-1} \tau \rho_+}  \Pi_+  I_+^j (\tau) 
 \end{equation} 	
with  
$$
I_\pm^j (\tau) = \int_0^\tau  e^{h^{-1} s( \rho_- - \rho_\pm)} \sum_{k=0}^{2j+1} (h^{-1} s  \langle \xi \rangle   )^k \tilde{B}_{j,k} \text{d}s,
$$
For $I_-^j $, the exponential term is equal to $1$ and by integration of $s^k$, we obtain:
\begin{equation}\label{I-}
I_-^j (\tau) =   \sum_{k=0}^{2j+1} (h^{-1} \tau  \langle \xi \rangle   )^{k+1} \frac{ h  \langle \xi \rangle^{-1} }{k+1} \tilde{B}_{j,k} .
\end{equation} 	
For $I_+^j $, let us introduce $P_k$ the polynomial of degree $k$ such that 
$$ \int_0^\tau e^{\lambda s} s^k \text{d}s = \frac{1}{\lambda^{k+1}} ( e^{\tau \lambda} P_k(\tau \lambda) - P_k(0) ),$$
for any $\lambda \in \cc^*$. With this notation in hand, we easily see that the term $e^{ \tau^h \rho_+}  \Pi_+  I_+^j (\tau)$ has the following form:
\begin{equation}\label{I+}
e^{ \tau^h \rho_+}  \Pi_+  I_+^j (\tau)   = \Pi_+ \,  \sum_{k=0}^{2j+1} \frac{h \langle \xi \rangle^k}{(\rho_--\rho_+)^{k+1}}  \tilde{B}_{j,k}  \Big(   e^{ \tau^h \rho_-} P_k(  \tau^h ( \rho_- - \rho_+)  ) - e^{ \tau^h \rho_+} P_k(  0 )   \Big),
\end{equation} 	
where $\tau^h := h^{-1} \tau$. Thus, combining \eqref{I-} and \eqref{I+} with \eqref{Agene}, \eqref{2ndAgene} and \eqref{expL0}, yields that 
 \begin{equation*}
 A_{j+1}  =   e^{h^{-1} \tau \rho_+} \Big( \Pi_+ {A_{j+1}}_{| \tau=0} -  \widetilde{{B}^+_{j+1}}\Big)  + e^{h^{-1} \tau \rho_-} \Big( \Pi_- {A_{j+1}}_{| \tau=0}  + \sum_{k=0}^{2(j+1)} (h^{-1} \tau  \langle \xi \rangle   )^{k}  \widetilde{{B}^-_{j+1,k}} \Big) ,
  \end{equation*} 	
where 
$$\widetilde{{B}^+_{j+1}}  = \Pi_+  \sum_{k=0}^{2j+1} \frac{h \langle \xi \rangle^k}{(\rho_--\rho_+)^{k+1}} P_k(0) \tilde{B}_{j,k} \,  \in \, h^{j+1} \sS^{-j-1},$$
and $\widetilde{{B}^-_{j+1,k}}  \in  h^{j+1} \sS^{-j-1} $ as a linear combination of products of $\Pi_- \in \sS^0$, $h \langle \xi \rangle^{-1}$ (or $h \langle \xi \rangle^k (\rho_--\rho_+)^{-k-1}$) belonging to $h \sS^{-1}$,  and of $\tilde{B}_{j,k} \in h^j \sS^{-j}$. 

Now, in order to have $A_{j+1} \in \sS^0$,  we let the contribution of the exponentially growing term vanish by choosing 
 \begin{equation*}
 \Pi_+ {A_{j+1}}(y, \xi, 0) = \widetilde{{B}^+_{j+1}} (y, \xi) . 
   \end{equation*} 
 Then, thanks to Remark \ref{isoPPi}, the boundary condition $P_-(y) {A_{j+1}}(y, \xi, 0)=0$ gives  
 \begin{equation*}
{A_{j+1}}(y, \xi, 0) = \frac{P_+ \Pi_+ }{k^{\varphi}_+}  \widetilde{{B}^+_{j+1}} (y, \xi) . 
   \end{equation*} 
 Finally, we have 
 $${A_{j+1}}(y, \xi, \tau ) = e^{h^{-1} \tau \rho_-(y, \xi)} \Big(  \frac{\Pi_- P_+ \Pi_+ }{k^{\varphi}_+}  \widetilde{{B}^+_{j+1}} (y, \xi)  + \sum_{k=0}^{2(j+1)} (h^{-1} \tau  \langle \xi \rangle   )^{k}  \widetilde{{B}^-_{j+1,k}}(y, \xi) \Big),$$
and  Proposition \ref{SolSystSj} is proven with 
$$ {B}_{j+1,0} = \frac{\Pi_- P_+ \Pi_+ }{k^{\varphi}_+}  \widetilde{{B}^+_{j+1}} + \widetilde{{B}^-_{j+1,0}},$$
 and for $k \geq 1$, ${{B}_{j+1,k}} = \widetilde{{B}^-_{j+1,k}}$.\qed\\

\begin{remark} The computation of each term $B_{j,0}$ can be done recursively, but this leads to complicated calculations. For example $B_{1,0}$ has the following form
\begin{align*}
B_{1, 0}(y,\xi)= - h\Pi_{+}a_{0}\left(  \frac{(z  + i \alpha\cdot\partial_{y} )}{2\l}  - \dfrac{i\alpha\cdot\partial_{y}\rho_{-}}{4\l^2}  \right) \Pi_- A_0(y,\xi),
\end{align*}
with $a_0(\ty)= i\alpha\cdot\tn^\varphi(\ty)$.
\end{remark}

Thanks to the relation \eqref{uh}, to any $A^h \in {\cP}_{h}^{0}$ we associate a bounded operator from $\mathit{L}^2(\rr^2)$ into $\mathit{L}^2(\rr^2\times (0,+ \infty))$. The boundedness in the variable $y \in \rr^2$ is a consequence of the Calderon-Vaillancourt theorem (see (\ref{Calderon-Vaillancourt})), and in the variable $\tau \in (0,+ \infty)$ it is essentially the multiplication by an $\mathit{L}^\infty$-function. Moreover, for $A_j$ of the form \eqref{eqAj}, we have the following mapping property which captures the Sobolev space regularity.
\begin{proposition}\label{SobolevMap}
Let $A_j$, $j\geq 0$, be of the form \eqref{eqAj}. Then, for any $s \geq -j-\frac12$, the operator ${\mathcal A}_j$ defined by 
\begin{align*}
{\mathcal A}_j : f \longmapsto  ({\mathcal A}_j f) (y,y_3) = \int_{\rr^2} A_j(y, h \xi , y_3) e^{i y \cdot \xi} \hat{f}(\xi) \mathrm{d} \xi
\end{align*}
gives rise to a bounded operator from $\mathit{H}^s(\rr^2)$ into $\mathit{H}^{s+j+\frac12}(\rr^2\times (0,+ \infty))$. Moreover,  for any $l \in [0, j+\frac12]$ we have:
\begin{equation}\label{SobolevEsti}
\| {\mathcal A}_j\|_{\mathit{H}^s \to \mathit{H}^{s+j+\frac12-l}}  =  O(h^{l-s}).
 \end{equation} 	

\end{proposition}
\textbf{Proof.} 
First, let us prove the result for $s = k -j-\frac12$, $k \in \nn$,  between the semiclassical Sobolev spaces 
\begin{align*}
\mathit{H}_{\text{scl}}^s(\rr^2)&:= \langle h D_y\rangle^{-s} \mathit{L}^2(\rr^2) \\
 \mathit{H}^{k}_{\text{scl}}(\rr^2\times (0,+ \infty))&: = \{ u \in L^2  ; \;  \langle h D_y\rangle ^{k_1} (h \partial_{y_3})^{k_2}  u \in \mathit{L}^2 {\text{ for }} (k_1,k_2) \in \nn^2, k_1+k_2 = k \},
 \end{align*}
where $\langle h D_y\rangle = \sqrt{-h^2 \Delta_{\rr^2} +I}$. 
Then, for $f \in  \mathit{H}^s(\rr^2)^{4}$, we have:
\begin{align}
\begin{split}
\| {\mathcal A}_j f \|^2_{\mathit{H}^{k}_{\text{scl}}(\rr^2\times (0,+ \infty))} = & \sum_{k_1+k_2 = k} \|  \langle h D_y\rangle^{k_1} (h \partial_{y_3})^{k_2} {\mathcal A}_j f \|^2_{\mathit{L}^2(\rr^2\times (0,+ \infty))} \\
\label{eqAj1}
=& \sum_{k_1+k_2 = k}  \int_0^{+\infty} \|\langle h D_y\rangle^{k_1} (h \partial_{y_3})^{k_2}( {\mathcal A}_j f) (\cdot , y_3) \|^2_{\mathit{L}^2(\rr^2)} \text{d}y_3.
\end{split}
\end{align}
Thanks to the ellipticity property \eqref{ellip}, for $A_j$ given by Proposition \ref{SolSystSj} we have:
$$(h \partial_{y_3})^{k_2} A_j(y, \xi, y_3) = h^j  b_j(y, \xi; y_3) e^{-h^{-1} y_3 \frac{c}2 \langle \xi \rangle}  \langle \xi \rangle^{k_2-j},$$
with $b_j$ satisfying, for any $(\alpha, \beta) \in {\nn}^2 \times {\nn}^2$ there exists  $C_{\alpha,\beta} >0$ such that:
$$  | \partial_y^\alpha \partial_{\xi}^\beta b_j(y, \xi ; y_3) | \leq C_{\alpha, \beta}, \quad \forall (y, \xi; y_3) \in \rr^2\times \rr^2\times (0,+ \infty).$$
Consequently, thanks ot the Calder\'on-Vaillancourt theorem (see \eqref{Calderon-Vaillancourt}), we can write:
\begin{align*}\langle h D_y\rangle^{k_1} (h \partial_{y_3})^{k_2} {\mathcal A}_j = h^j {\mathcal B}_j (y_3) \langle h D_y\rangle^{k_1+k_2-j} e^{-h^{-1} y_3 \frac{c}{2} \langle h D_y\rangle},
\end{align*}
with $({\mathcal B}_j (y_3))_{y_3>0}$ a family of bounded operators on $\mathit{L}^2(\rr^2)$, and uniformly bounded with respect to $y_3>0$.
Then, for $f \in  \mathit{H}^s(\rr^2)^{4}$, we have:
$$ \|\langle h D_y\rangle^{k_1} (h \partial_{y_3})^{k_2}( {\mathcal A}_j f) (\cdot , y_3) \|^2_{\mathit{L}^2(\rr^2)} \lesssim
h^j  \|  \langle h D_y\rangle^{k_1+k_2-j} e^{-h^{-1} y_3 \frac{c}2 \langle h D_y\rangle}  f    \|^2_{\mathit{L}^2(\rr^2)},$$
and from \eqref{eqAj1}  we deduce that
\begin{align*}
\| {\mathcal A}_j f  \|^2_{\mathit{H}^{k}_{\text{scl}}(\rr^2\times (0,+ \infty))} \lesssim  h^{2j+1}   \|  \langle h D_y\rangle^{k-j-\frac12}   f    \|^2_{\mathit{L}^2(\rr^2)} =  h^{2j+1}   \|    f    \|^2_{\mathit{H}_{\text{scl}}^{k-j-\frac12}(\rr^2)},
\end{align*}
where we used that for any $l \in \nn$, $f \in \mathit{H}_{\text{scl}}^{l-\frac12}(\rr^2)$, 
\begin{align*}
\|  \langle h D_y\rangle^{l} e^{-h^{-1} y_3 \frac{c}2 \langle h D_y\rangle}  f    \|^2_{\mathit{L}^2(\rr^2)} &= \langle e^{-h^{-1} y_3 c \langle h D_y\rangle}  \langle h D_y\rangle^{l}  f \, , \, \langle h D_y\rangle^{l}  f  \rangle_{\mathit{L}^2} \\
&= - \frac{h}{c} \frac{\partial}{\partial_{ y_3}} \langle e^{-h^{-1} y_3 c \langle h D_y\rangle}  \langle h D_y\rangle^{l-1}  f \, , \, \langle h D_y\rangle^{l}  f  \rangle_{\mathit{L}^2} .\end{align*}
By interpolation arguments we thus deduce that for any $j \in \nn$, $s \geq -j-\frac12$, it holds that
\begin{equation*}
\| {\mathcal A}_j\|_{\mathit{H}_{\text{scl}}^s \to \mathit{H}_{\text{scl}}^{s+j+\frac12}}  =  O(h^{j+ \frac12}),
 \end{equation*} 
 proving the estimate  \eqref{SobolevEsti} and completing the proof of the proposition. \qed


\begin{proposition}\label{ApproxSystSC}
Let $f \in H^s(\rr^2)$ and $A_j$, $j\geq 0$, be as in Propositions \ref{SolSystS0} and \ref{SolSystSj}. Then for any $N \geq -s-\frac12$, the function
$u_N^h = \sum_{j=0}^N h^j {\mathcal A}_j f$ satisfies:
\begin{equation}\label{ApproxSystSC}
	\left\{
	\begin{aligned}
	h \partial_{\tau} u_N^h & - L_0(y, h D_y)  u_N^h - h L_1(y) u_N^h =  h^{N+1} {\mathcal R}_N^h f ,  &  \quad \text{ in }  &  \rr^2 \times (0, + \infty),\\
	P_- u_N^h  & = f,  &  \text{ on } &  \rr^2 \times \{0\},
	\end{aligned}
	\right.
	\end{equation}
with 
$${\mathcal R}_N^h: f \longmapsto  \int_{\rr^2} \Big( L_1 A_N - \partial_\xi L_0  \cdot \partial_y  A_N   \Big) (y, h \xi , \tau) e^{i y \cdot \xi} \hat{f}(\xi) \mathrm{d} \xi,   $$
 a bounded operator from $H^s(\rr^2)$ into $H^{s+N+\frac12}(\rr^2\times (0,+ \infty))$ satisfying  for any $l \in [0, N+\frac12]$:
\begin{equation}\label{EstiRN}
\| {\mathcal R}^h_N\|_{\mathit{H}^s \to \mathit{H}^{s+N+\frac12-l}}  =  O(h^{l-s}).
 \end{equation} 	
\end{proposition}
\textbf{Proof.} 
By construction of the sequence $(A_j)_{j \in \{0,\cdots,N-1\} }$ we have the system \eqref{ApproxSystSC} with  ${\mathcal R}^h_N = \mathit{Op}^h (r_N^h(\cdot, \cdot, \tau) )$,
$$ r_N^h(y, \xi , \tau) = \Big( L_1 A_N - \partial_\xi L_0  \cdot \partial_y  A_N   \Big) (y, \xi , \tau) ,$$
(see the beginning of Section \ref{ssSclBPb}). As in the proof of  Proposition \ref{SolSystSj},
$r_N^h$ has the form \eqref{Agene1} (with $j=N$).  Then,  as in the proof of Proposition \ref{SobolevMap} we obtain the estimate \eqref{EstiRN}.\qed\\


\subsection{Proof of Theorem \ref{pseudo de Am}}

In this section, we apply the above construction in order to prove Theorem \ref{pseudo de Am}.

Let $g\in P_- \mathit{H}^{1/2}(\partial\O)^4$,  $(U_{\varphi}, V_{\varphi},\var)$ a chart of the atlas $\mathbb{A}$ and $\psi_1,\psi_2\in\mathit{C}^{\infty}_{0}(U_{\varphi})$. Then $f:= (\varphi^{-1})^* (\psi_2 g) $ is a function of $\mathit{H}^{1/2}(V_{\varphi})^4$ which can be extended by $0$ to a function of $\mathit{H}^{1/2}(\rr^2)^4$. Then for $h = 1/m$ and any $N \in \nn$, the previous construction provides a function $u_N^h \in \mathit{H}^{1}(\rr^2 \times (0, + \infty))^4$ satisfying
\begin{equation*}
	\left\{
\begin{aligned}	
(\widetilde{D}^\varphi_m  -z ) u_N^h	=  & h^{N+1} {\mathcal R}_N^h f ,  \quad  &\text{ in } \rr^2 \times (0, \varepsilon) ,\\
	\G_- u_N^h= & f,  \quad \quad  & \text{ on } \rr^2 \times \{0\}, 
	\end{aligned}
	\right.
	\end{equation*}	
with $u_N^h = \sum_{j=0}^N h^j {\mathcal A}_j f$ (see Proposition \ref{SobolevMap}) and ${\mathcal R}^h_N f \in H^{N+1}(\rr^2 \times (0, \varepsilon))$ with norm in $H^{N+1-l}$,   $l \in [0, N+ \frac12]$, bounded by  $O(h^{l-\frac12})$. Consequently, $v_N^h := \phi^* u_N^h$, defined on $\mathcal{V}_{\varphi,\varepsilon}$, satisfies:
\begin{equation*}
	\left\{
\begin{aligned}\label{T1}	
({D}_m  -z ) v_N^h	= & h^{N+1} \phi^* ({\mathcal R}_N^h f ),  \quad  &\text{ in } \mathcal{V}_{\varphi,\varepsilon},\\
	\G_- v_N^h= & \psi_2 g,  \quad   & \text{ on } U_{\varphi}.
	\end{aligned}
	\right.
	\end{equation*}	
 Now, let $E_{m}^{\O}(z)[\psi_2 g] \in \mathit{H}^1(\O)^4$ be as in Definition  \ref{PSdef}. Since  $\G_- v_N^h =\G_- E_{m}^{\O}(z)[\psi_2 g]=\psi_2 g$,  then the following equality holds in $\mathcal{V}_{\varphi,\varepsilon}$: 
$$v_N^h - E_{m}^{\O}(z)[\psi_2 g] =  h^{N+1}  (\Hm-z)^{-1} \phi^* \Big({\mathcal R}_N^h (\varphi^{-1})^* (\psi_2 g )\Big).$$
From this, we deduce that
$$\psi_1 \mathscr{A}_m \psi_2 (g) := \psi_1  \G_+E_{m}^{\O}(z)[\psi_2 g] = \psi_1  \G_+  v_N^h - h^{N+1}   \psi_1  \G_+ (H_{\text{MIT}}-z)^{-1} \phi^* \Big({\mathcal R}_N^h (\varphi^{-1})^* (\psi_2 g )\Big).$$
Since $\phi\downharpoonright_{ U_{\varphi}} = \varphi$, for any $u \in \mathit{H}^1( V_{\varphi} \times (0,\varepsilon))^4$, we have that
\begin{align*}
 \G_+  \phi^*( u) = \varphi^* ( P_+ u\downharpoonright_{V_{\varphi} \times \{0\}}), \quad \psi_1 \G_+  v_N^h = \psi_1 \varphi^*  \mathit{Op}^h(a_N^h) (\varphi^{-1})^* \psi_2 g,
\end{align*} 
with
\begin{equation*}
a_N^h(\ty, \xi)= \sum_{j=0}^N h^j P_+  A_j (y, \xi, 0) = \sum_{j=0}^N h^j P_+  B_{j,0}(y, \xi), 
\end{equation*}
where  $B_{j,0} \in  h^j \sS^{-j} $    are introduced in Proposition \ref{SolSystSj}.  Thus, from Proposition \ref{SolSystS0}, in local coordinates,  the principal semiclassical symbol of $\mathscr{A}_m$ is given by
$$P_+  B_{0,0}(y, \xi)= P_+ A_0 (y, \xi , 0) =   \frac{P_+ \Pi_- P_-}{k^{\varphi}_+} (y,\xi).$$
Thanks to the property  \eqref{PPiP}  it is equal to
$$- \Theta^{\varphi}P_- (y,\xi) = \frac{S\cdot (\xi \wedge n^\varphi(y))}{\sqrt{\langle G(y)^{-1}\xi,\xi\rangle  + 1 } +1 }P_- (y,\xi) .$$

We conclude the proof of Theorem \ref{pseudo de Am} from results of Section \ref{sDirac} and by proving the following Lemma which is a consequence of the above considerations, the regularity estimates from Theorem \ref{MITLipshictz}-$(\mathrm{iii})$,  Theorem \ref{the de regu}-$(\mathrm{i})$ and Proposition \ref{definitionPSOP}.
\begin{lemma}\label{esti12Am}
Let $\psi_1,\psi_2\in\mathit{C}^\infty (\S)$ such that  $\mathrm{supp}(\psi_1)\cap \mathrm{supp}(\psi_2)=\emptyset$. Then, for $m_0>0$ sufficiently large, $m\geqslant m_0$,  and for any $(k,N)\in\nn^{\ast}\times \nn^{\ast}$  it holds that
 \begin{align*}
 \left\|  \psi_1\mathscr{A}_m\psi_2 \right\|_{P_-\mathit{H}^{1/2}(\S)^4\rightarrow P_+\mathit{H}^{k}(\S)^4}= \mathcal{O}(m^{-N}).
\end{align*}
 \end{lemma}
 \textbf{Proof.} Let  $\psi_1,\psi_2\in\mathit{C}^\infty (\S)$ with disjoint supports. Thanks to Theorem \ref{MITLipshictz}-$(\mathrm{iii})$ and  Theorem \ref{the de regu}-$(\mathrm{i})$, to prove the lemma it suffices to show that for any $(N_1,N_2) \in \nn^2$, there exists $C_{N_1,N_2}$ such that 
\begin{align}\label{Estimation pr m}
\begin{split}
\| (\psi_1  \mathscr{A}_m    \psi_2)g \|_{P_+ H^{N_2+\frac12}(\S)^4}
 \leq& \frac{C_{N_1,N_2}}{\sqrt{m}} \left(\Pi_{k=0}^{N_2-1} \|  (\Hm-z)^{-1} \|_{\mathit{H}^{k}(\O)^4 \rightarrow \mathit{H}^{k+1}(\O)^4 }\right)\\
 &\times  \|  (\Hm-z)^{-1} \|_{L^2(\O)^4 \rightarrow  \mathit{L}^2(\O)^4 }^{N_1}\| g \|_{P_-\mathit{H}^{1/2}(\S)^4}.
 \end{split}
 \end{align}
 For this,  let us introduce $\Phi_1 \in \mathit{C}_0^\infty (\overline{\O})$ such that $\Phi_1=1$ near $\mathrm{supp}(\psi_1)$ and $\Phi_1=0$ near $\mathrm{supp}(\psi_2)$. Thus for $g\in P_- \mathit{H}^{1/2}(\partial\O)^4$ and $E_{m}^{\O}(z)[\psi_2 g] \in \mathit{H}^1(\O)$  as in Definition  \ref{PSdef}, the function $u_{1,2}:=\Phi_1 E_{m}^{\O}(z)[\psi_2 g]$ satisfies:
 \begin{equation*}
	\left\{
\begin{aligned}
({D}_m  -z ) u_{1,2}	= & [ {D}_0 \, , \, \Phi_1] E_{m}^{\O}(z)[\psi_2 g] ,  \quad  &\text{ in } \O,\\
	\G_- u_{1,2} = & \Phi_1\downharpoonright_{ \S} \psi_2 g =0,  \quad    &\text{ on } \S.
	\end{aligned}
	\right.
	\end{equation*}	
Then, $u_{1,2} = (\Hm-z)^{-1}   [ {D}_0 \, , \, \Phi_1] E_{m}^{\O}(z)[\psi_2 g]$,  and for any $\widetilde{\Phi_1}  \in \mathit{C}_0^\infty (\overline{\O})$ equals to $1$  near $\mathrm{supp}(\psi_1)$ we have:
$$ \psi_1  \mathscr{A}_m    \psi_2  (g) = \psi_1  \G_+  \widetilde{\Phi_1} (\Hm-z)^{-1}   [ {D}_0 \, , \, \Phi_1] E_{m}^{\O}(z)[\psi_2 g].$$
Moreover, by choosing $\widetilde{\Phi_1} $ such that $\widetilde{\Phi_1}  \prec \Phi_1$, that is  $\Phi_1 = 1$ on $\mathrm{supp}(\widetilde{\Phi_1})$, both functions $\widetilde{\Phi_1} $ and $[ {D}_0 \, , \, \Phi_1]$ have disjoint supports, and we can then apply the following telescopic formula:
\begin{align*}
\widetilde{\Phi_1} (\Hm-z)^{-1}  (1- \chi_1) =& \widetilde{\Phi_1} (\Hm-z)^{-1} [D_0 , \chi_N] \cdots  (\Hm-z)^{-1} [D_0 , \chi_2]\\
& (\Hm-z)^{-1}  (1- \chi_1) ,
\end{align*}
for $(\chi_i)_{1 \leq i \leq N}$ a family of compactly supported smooth functions such that $\widetilde{\Phi_1} \prec \chi_N \prec \chi_{N-1} \prec \cdots   \prec \chi_1 \prec \Phi_1$. Since $ [ {D}_0 \, , \, \Phi_1] = (1- \chi_1)  [ {D}_0 \, , \, \Phi_1] $, the above telescopic formula allows us to write $\psi_1  \mathscr{A}_m    \psi_2  (g) $ as a product of $N$ cutoff resolvents of $\Hm$. Now, by Proposition \ref{definitionPSOP} we have 
\begin{align*}
 \left|\left| E_{m}^{\O}(z)[\psi_2 g] \right|\right|_{\mathit{L}^2(\O)^4}&\lesssim \frac{1}{\sqrt{m}}\left|\left| g \right|\right|_{\mathit{L}^{2}(\S)^4}. 
 \end{align*}
 Thus, using the continuity of $\G_+$ from $\mathit{H}^{N_2+1}(\O)$ to $\mathit{H}^{N_2+\frac12}(\S)$, we then get the estimation \eqref{Estimation pr m}  for $N = N_1+N_2$, finishing the proof of the lemma. \qed\\

\begin{remark}\label{remarkpseu} Note that for any $m>0$ and $z\in\rho(\Hm)$,  the parametrix we have constructed for $\mathscr{A}_m$ is valid from the classical pseudodifferentiel point of view. Actually, Lemma \ref{esti12Am} is the only result where the assumption that $m$ is big enough has been assumed, and it is exclusively required to ensures that away from the diagonal the operator $\mathscr{A}_m$ is negligible in $1/m$. In the same vein, if $m$ is fixed then the proof of Lemma \ref{esti12Am} still ensures that away from the diagonal $\mathscr{A}_m$ is regularizing. Consequently, we deduce that for any $m>0$ and $z\in\rho(\Hm)$, the operator $\mathscr{A}_m$ is a homogeneous pseudodifferential operator of order $0$, and that 
\begin{align*}
\mathscr{A}_m&=\frac{D_\Sigma}{\sqrt{-\Delta_\S}} P_- \quad  \mathrm{mod} \, \mathit{Op} \sS^{-1}(\S),
\end{align*}
which is in accordance with Theorem \ref{Th SymClas}.
\end{remark}
\begin{remark} If $\O$ is the upper half-plane $\{ (x_1,x_2,x_3) \in \rr^3; \; x_3 >0 \}$, we easily obtain that $\mathscr{A}_m$ is a Fourier multiplier with symbol
$$a_m(\xi)=  -\frac{i\alpha_{3}(\alpha\cdot \xi  -z)}{\sqrt{|\xi|^2 + m^2}+m}P_{-}    .$$
\end{remark}
 \section{Resolvent convergence to the MIT bag model}\label{SEC6}
In the whole section,  $\O\subset\rr^3$ denotes a bounded smooth domain, we set 
  \begin{align*}
\O_{i}=\O,  \quad \O_{e}=\rr^3\setminus\overline{\O}\quad\text{ and }\,\, \S=\partial\O,
	\end{align*}
and we let $n$ be the outward (with respect to $\O_{i}$) unit normal vector field on $\S$.
	
Fix $m>0$ and let  $M>0$. Consider the perturbed Dirac operator  	
\begin{align*}
H_M\var= (D_m +M\beta 1_{\O_{e}})\var,  \quad \forall \var\in \mathrm{dom}(H_M):=  \mathit{H}^1(\rr^3)^4,
\end{align*} 
where $1_{\O_{e}}$ is the characteristic function of $\O_e$. Using Kato-Rellich theorem and Weyl's theorem, it is easy to see that $(H_M, \mathrm{dom}(H_M))$ is self-adjoint and that 
\begin{align*}
	\mathrm{Sp}_{\mathrm{ess}}(H_{M}) =(-\infty,-(m+M)]\cup [m+M,+\infty),\\
	\mathrm{Sp}(H_{M})\cap(-(m+M),m+M) \text{ is purely discrete}.
\end{align*} 
Now, let $\Hm$ be the MIT bag operator acting  on $\mathit{L}^2(\O_i)^4$, that is 
\begin{align*}
\Hm v=D_mv\quad \forall v\in\mathrm{dom}(\Hm):=  \left\{ v \in\mathit{H}^1(\Omega_{i} )^4 : P_-t_{\S} v=0 \text{ on }\S \right\},
\end{align*} 
where $t_{\S}$ and $P_{\pm}$ are the trace operator and the orthogonal projection from Subsection \ref{SubNot}.

The aim of this section is to use the properties of the Poincar\'e-Steklov operators carried out in the previous sections to study  the resolvent of $H_M$ when $M$ is large enough. Namely, we give a Krein-type resolvent formula in terms of the resolvent of  $\Hm$, and we show that the convergence of $H_M$ toward $\Hm$ holds in the norm resolvent sense  with a convergence rate of $\mathcal{O}(1/M)$, which improves the result of \cite{BCLS}. 

 Before stating the main results of this section, we need to introduce some notations and definitions.  First, we introduce the following Dirac auxiliary operator
\begin{align*}
	\widetilde{H}_{M}u&=D_{m+M}u, \quad \forall u\in \mathrm{dom}(\widetilde{H}_{M}):=  \left\{ u\in\mathit{H}^1(\O_{e})^4 : P_+ t_{\S}u=0\text{ on } \S \right\}.
\end{align*} 
Notice that $\widetilde{H}_{M}$ is the MIT bag operator on $\O_e$ (the boundary condition is with $P_+$ because the normal $n$ is incoming for $\O_e$).
Since $\O_e$ is unbounded, Theorem  \ref{MITLipshictz} together with Remark \ref{RemMITLipshictz} imply that $(\widetilde{H}_{M},\mathrm{dom}(\widetilde{H}_{M}))$ is self-adjoint and that 
\begin{align*}
\mathrm{Sp}(\widetilde{H}_{M})=\mathrm{Sp}_{\mathrm{ess}}(\widetilde{H}_{M}) =(-\infty,-(m+M)]\cup [m+M,+\infty).
\end{align*}
In particular, $\rho(H_{M}) \subset \rho(  \widetilde{H}_{M})  )$. 
Let  $ z\in\rho(\Hm)\cap \rho(  \widetilde{H}_{M})  )$, 
$g\in P_-\mathit{H}^{1/2}(\S)^4$ and  $h\in P_+\mathit{H}^{1/2}(\S)^4$.  We denote by $E^{\O_{i}}_{m}(z): P_-\mathit{H}^{1/2}(\S)^4\rightarrow \mathit{H}^{1}(\O_{i})^4$ the unique solution of the boundary value problem:
	\begin{equation}\label{L1b}
	\left\{
	\begin{aligned}
	(D_m-z)v&=0, \quad \text{ in } \O_{i},\\
	P_-t_{\S}v&= g,  \quad \text{ in } \S.
	\end{aligned}
	\right.
	\end{equation}
Similarly,    we denote by $E^{\O_{e}}_{m+M}(z): P_+\mathit{H}^{1/2}(\S)^4\rightarrow \mathit{H}^{1}(\O_{e})^4$   the unique solution of the boundary value problem:  
\begin{equation}\label{L2b}
\left\{
\begin{aligned}
(D_{m+M}-z)u&=0, \quad \text{ in } {\O_{e}},\\
P_+t_\S u&= h,  \quad \text{ in } \S.
\end{aligned}
\right.
\end{equation}
Define the Poincar\'e-Steklov operators associated to the above problems by 
\begin{align*}
\mathscr{A}_{m}^{i}= P_+ t_{\S} E^{\O_{i}}_{m}(z)P_{-} \quad \text{and}\quad \mathscr{A}_{m+M}^{e}=P_- t_{\S} E^{\O_{e}}_{m+M}(z)P_{+}.
\end{align*}

 \begin{notation} In the sequel we shall denote by $R_{M}(z)$, $\widetilde{R}_{M}(z)$ and  $R_{\text{MIT}}(z)$ the resolvent of $H_{M}$, $\widetilde{H}_{M}$ and $\Hm$, respectively. We also use the notations:
\begin{itemize}
\item  $\G_\pm=P_\pm t_{\S}$ and  $\G=\G_+r_{\O_{i}}+\G_-r_{\O_{e}}$.
\item $E_M(z)=e_{\O_{i}}E^{\O_{i}}_{m}(z)P_{-} + e_{\O_{e}}E^{\O_{e}}_{m+M}(z)P_{+}$.
\item $ \widetilde{ R}_{\text{MIT}}(z)= e_{\O_{i}}R_{\text{MIT}}(z)r_{\O_{i}}+ e_{\O_{e}}\widetilde{R}_{M}(z)r_{\O_{e}}$.
\end{itemize}
\end{notation}
 With these notations in hand, we can state the main results of this section. The following theorem is the main tool to show the large coupling convergence with a rate of convergence of $\mathcal{O}(1/M)$.
\begin{theorem}\label{Formulederes} There is $M_0>0$ such that for all $M>M_0$ and all $ z\in\rho(\Hm)\cap\rho(H_{M})$,  the operator $\Psi_M(z):=\left( I-\mathscr{A}^{i}_{m}-\mathscr{A}^{e}_{m+M} \right)$ is bounded invertible in $\mathit{H}^{1/2}(\S)^4$, and the inverse is given by 
 \begin{align*}
  \Psi^{-1}_M(z)=  \left( \mathit{I}_4-\mathscr{A}_{m}^{i}\mathscr{A}_{m+M}^{e}-\mathscr{A}_{m+M}^{e}\mathscr{A}_{m}^{i} \right)^{-1}\left( I+\mathscr{A}_{m}^{i}+ \mathscr{A}_{m+M}^{e}  \right), 
  \end{align*} 
and the following resolvent formula holds:
 \begin{align}\label{fullreso}
R_M(z)= \widetilde{ R}_{\text{MIT}}(z)+E_M(z)\Psi^{-1}_M(z)\G\widetilde{ R}_{\mathrm{MIT}}(z).
\end{align} 
\end{theorem}
	\begin{remark}
	By Proposition \ref{definitionPSOP} $(\mathrm{ii})$  we have that 
	\begin{align*}
	\left(  E^{\O_i}_{m}(z)\right)^{\ast}=-\beta \G_+R_\mathrm{MIT} (\overline{z}) \quad\text{and}\quad  \left(  E^{\O_e}_{m+M}(z)\right)^{\ast}=-\beta \G_-\widetilde{R}_{M} (\overline{z}),
	\end{align*}
 for any $ z\in\rho(\Hm)\cap\rho(H_{M})$. Thus, the resolvent formula \eqref{fullreso} can be written in the form 
 \begin{align*}
R_M(z)= \widetilde{ R}_{\mathrm{MIT}}(z)- (\beta\G\widetilde{ R}_{\mathrm{MIT}}(\overline{z}))^{\ast} {\Psi}_M^{-1}(z)\G\widetilde{ R}_{\mathrm{MIT}}(z).
\end{align*} 
	\end{remark}
	
Before going through the proof of Theorem \ref{Formulederes}  we first establish a regularity result that will play a crucial role in the rest of this section. It concerns the dependence on the parameter $M$ of the norm of an auxiliary operator which involves the composition of the operators $\mathscr{A}_{m}^{i}$ and $\mathscr{A}_{m+M}^{e}$.  In the proof we use the symbols $\hat{u} $ and $\mathcal{F}[u]$ to denote the Fourier transform of $u$. 
\begin{proposition}\label{PseudESt}  Let $\mathscr{A}_{m}^{i}$ and $\mathscr{A}_{m+M}^{e}$ be as above. Then,  there is $M_0>0$ such that for every  $M>M_0$ and all $ z\in\rho(\Hm)\cap\rho(H_M)$  the following hold true:    
 \begin{itemize}
  \item[(i)] For any $s\in\rr$ the operator $\Xi_M(z):\mathit{H}^{s}(\S)^4\longrightarrow \mathit{H}^{s}(\S)^4$ defined by 
\begin{align}\label{Opein1}
\Xi_M(z)= \left( \mathit{I}_4-\mathscr{A}_{m}^{i}\mathscr{A}_{m+M}^{e}-\mathscr{A}_{m+M}^{e}\mathscr{A}_{m}^{i} \right)^{-1},
\end{align}
is everywhere defined and uniformly bounded with respect to $M$.
  \item[(ii)] The Poincar\'e-Steklov operator, $\mathscr{A}_{m+M}^{e}$, satisfies the estimate
   \begin{align*}
 \left|\left|  \mathscr{A}_{m+M}^{e}  \right|\right|_{P_+\mathit{H}^{s+1}(\S)^4 \rightarrow  P_-\mathit{H}^{s}(\S)^4 } \lesssim h,\quad \forall s\in\rr.
 \end{align*}
\end{itemize}
\end{proposition}

\textbf{Proof.} ($\mathrm{i}$) Set $\tau := (m+M)$, then the result essentially follows from the fact that $\Xi_{M}(z)$ is a $1/\tau$-pseudodifferential operator of order $0$. Indeed,  fix $ z\in\rho(\Hm)\cap\rho(H_{M})$ and set $h= \t^{-1}$. Then, from Theorem \ref{Th SymClas} and Remark \ref{remarkpseu} we know that $\mathscr{A}_{m}^{i}$ is a homogeneous pseudodifferential operator of order $0$. Thus $\mathscr{A}_{m}^{i}$ can also be viewed as a $h$-pseudodifferential operators of order $0$. That is,  $\mathscr{A}_{m}^{i}\in \mathit{Op}^h \sS^{0}(\S)$, and in local coordinates, its semiclassical principal symbol is given by 
	$$p_{h, \mathscr{A}^{i}_{m}}(x,\xi)=\frac{S\cdot(\xi\wedge n(x))P_-}{|\xi\wedge n(x)|} ,$$
	where we identify $\xi \in \rr^2$ with $\overline{\xi} = (\xi_1, \xi_2,  0)^t \in \rr^3$, and for $x= \varphi(\tx) \in \Sigma$, $n(x)$ stands for  $n^\varphi(\tx)$.
	 Similarly, thanks to Theorem \ref{pseudo de Am},  we also know that for $h_0$ sufficiently small (and hence $M_0$ big enough) and all $h<h_0$, $\mathscr{A}_{m+M}^{e}$ is a $h$-pseudodifferential operator and that 
\begin{align*}
\mathscr{A}_{m+M}^{e}\in \mathit{Op}^h \sS^{0}(\S) ,\quad p_{h, \mathscr{A}^{e}_{m+M}}(x,\xi)=-\frac{S\cdot(\xi\wedge n(x))P_+}{\sqrt{|\xi\wedge n(x)|^2+1}+1} .
\end{align*}
Therefore, the symbol calculus yields for all $h<h_0$ that $\left( \mathit{I}_4-\mathscr{A}_{m}^{i}\mathscr{A}_{m+M}^{e}-\mathscr{A}_{m+M}^{e}\mathscr{A}_{m}^{i} \right)$ is a $1/\tau$-pseudodifferential operator of order $0$.
Now, a simple computation using  Lemma \ref{properties3.1} yields that 
\begin{align*}
\frac{S\cdot(\xi\wedge n(x))P_\pm S\cdot(\xi\wedge n(x))P_\mp}{|\xi\wedge n(x)|( \sqrt{|\xi\wedge n(x)|^2+1}+1)}= \frac{|\xi\wedge n(x)|P_\mp}{\sqrt{|\xi\wedge n(x)|^2+1}+1}.
\end{align*}
Thus 
\begin{align*}
  \mathit{I}_4- p_{h, \mathscr{A}^{i}_{m}}(x,\xi)p_{h, \mathscr{A}^{e}_{m+M}}(x,\xi) &-p_{h, \mathscr{A}^{e}_{m+M}}(x,\xi)p_{h, \mathscr{A}^{i}_{m}}(x,\xi)=  \mathit{I}_4 + \frac{|\xi\wedge n(x)|}{\sqrt{|\xi\wedge n(x)|^2+1}+1}\\
  &=  \frac{\sqrt{|\xi\wedge n(x)|^2+1}+1 +|\xi\wedge n(x)| }{\sqrt{|\xi\wedge n(x)|^2+1}+1} \gtrsim 1.
    \end{align*}
From this, we deduce that $\left( \mathit{I}_4-\mathscr{A}_{m}^{i}\mathscr{A}_{m+M}^{e}-\mathscr{A}_{m+M}^{e}\mathscr{A}_{m}^{i} \right)$ is elliptic in $ \mathit{Op}^h \sS^{0}(\S)$. Thus,  $\Xi_M(z)\in \mathit{Op}^h \sS^{0}(\S)$, and in local coordinates,  its semiclassical principal symbol is given by 
\begin{align*}
 p_{h, \Xi_M(z)}(x,\xi)=   \frac{\sqrt{|\xi\wedge n(x)|^2+1}+1}{\sqrt{|\xi\wedge n(x)|^2+1}+1 +|\xi\wedge n(x)| }.
    \end{align*}
As $\Xi_M(z)$ is a  $h$-pseudodifferential operators of order $0$,  it follows that  $ \Xi_M(z) :\mathit{H}^{s}(\S)^4\rightarrow \mathit{H}^{s}(\S)^4$ is well-defined and uniformly bounded with respect to $M$, for any $s\in\rr$,  proving the statement ($\mathrm{i}$) of the theorem. 

   The proof of the statement  ($\mathrm{ii}$) follows the standard arguments of the proof of the boundedness of classical pseudodifferential operators. Indeed, for $\t=m+M$, as a classical pseudodifferential operator, $\mathscr{A}^{e}_{\t}$ is given by:
   $$\mathscr{A}^{e}_{\t} = \frac{1}{\t} \mathit{Op}( a_{0,\t} ) + \frac{1}{\t} \mathit{Op}( a_{1,\t} ),$$
   with $a_{1,\t}$ uniformly bounded with respect to $\t \gg1$ in $ \sS^{0}$ and 
   $$  a_{0,\t} (x,\xi)   = -\frac{S\cdot(\xi\wedge n(x))P_+}{\sqrt{\t^{-2} |\xi\wedge n(x)|^2+1}+1}   $$
  which is  uniformly bounded in $ \sS^{1}$. Then $(\mathrm{ii})$ is a consequence of the Calder\'on-Vaillancourt theorem (see \eqref{Calderon-Vaillancourt}).

 \qed\\
 
We can now give the proof of Theorem \ref{Formulederes}.
	
	\textbf{Proof of Theorem \ref{Formulederes}.} Let $M_0$ be as in Proposition \ref{PseudESt} and $M>M_0$,  fix $ z\in\rho(\Hm)\cap\rho(H_{M})$ and let $f\in\mathit{L}^{2}(\rr^3)^4$. We set 
	\begin{align*}
	v=r_{\O_{i}}R_{M}(z)f\quad \text{ and} \quad u= r_{\O_{e}}R_{M}(z)f.
	\end{align*}
	Then $u$ and $v$ satisfy the following system 
	\begin{equation*}
	\left\{
	\begin{aligned}
	(D_m- z) v &=f \quad&\text{ in }\Omega_{i},\\
	(D_{m+M}-z)u &=f \quad&\text{in }\O_{e}, \\
	P_-t_{\S}v&= P_-t_{\S}u\quad& \text{ on } \S,\\ 
	P_+t_{\S}v&= P_+t_{\S}u\quad& \text{ on } \S.
	\end{aligned}
	\right.
	\end{equation*}  
	Since $E^{\O_{i}}_{m}(z)$ (resp. $ E^{\O_{e}}_{m+M}(z)$) gives the unique solution to the boundary value problem \eqref{L1b} (resp. \eqref{L2b}), and 
	\begin{align*}
	\G_-R_{\text{MIT}}(z)r_{\O_{i}}f= 0 \quad\text{and}\quad \G_+\widetilde{ R}_{M}(z)r_{\O_{e}}f= 0,
	\end{align*} 
	 if we let  
	\begin{align*}
	\varphi=\G_-u\quad \text{ and} \quad \psi= \G_+v,
	\end{align*} 
	then it is easy to check that 
	\begin{equation}\label{S2} 
	\left\{
	\begin{aligned}
	v &=R_{\text{MIT}}(z)r_{\O_{i}}f +E^{\O_{i}}_{m}(z)\varphi,\\
	u &=\widetilde{ R}_{M}(z)r_{\O_{e}}f +E^{\O_{e}}_{m+M}(z)\psi.
	\end{aligned}
	\right.
	\end{equation}
	 Hence, to get an explicit formula for $R_M(z)$ it remains to find the unknowns $\var$ and $\psi$. For this,  note that from \eqref{S2} we have
	\begin{equation}\label{S3} 
	\left\{
	\begin{aligned}
	\psi &= \G_+r_{\O_{i}}R_M(z)f =\G_+R_{\text{MIT}}(z)r_{\O_{i}}f +\G_+E^{\O_{i}}_{m}(z)[\varphi],\\
	\var&=  \G_-r_{\O_{e}}R_M(z)f =\G_-\widetilde{ R}_{M}(z)r_{\O_{e}}f +\G_-E^{\O_{e}}_{m+M}(z)[\psi].  
	\end{aligned}
	\right.
	\end{equation} 
Substituting the values of $\psi$ and $\var$ (from \eqref{S3}) into the system \eqref{S2},  we obtain 
	\begin{align}\label{preformula}
	\begin{split}
	R_M(z)=&e_{\O_{i}}R_{\text{MIT}}(z)r_{\O_{i}}+ e_{\O_{e}}\widetilde{R}_{M}(z)r_{\O_{e}}\\
	&+ \left(e_{\O_{i}} E^{\O_{i}}_{m}(z) \G_-r_{\O_{e}}+  e_{\O_{e}} E^{\O_{e}}_{m+M}(z) \G_+r_{\O_{i}}\right) R_M(z)  \\
	=& \widetilde{ R}_{\text{MIT}}(z)+E_M(z)\G R_{M}(z).
	\end{split}
	\end{align} 
	Note that, by definition of the Poincar\'e-Steklov operators, \eqref{S3} is equivalent to 
\begin{equation}\label{S3prime} 
\left\{
\begin{aligned}
\psi &= \G_+R_{\text{MIT}}(z)r_{\O_{i}}f +\mathscr{A}^{i}_{m}(\varphi),\\
\var&=  \G_-\widetilde{ R}_{M}(z)r_{\O_{e}}f +\mathscr{A}_{m+M}^{e}(\psi).  
\end{aligned}
\right.
\end{equation} 
Thus,  applying  $\G$ to the identity \eqref{preformula} yields that 
	\begin{align*}
	\G\widetilde{ R}_{\text{MIT}}(z)=   \left( I-\mathscr{A}^{i}_{m}-\mathscr{A}^{e}_{m+M} \right)\G R_M(z) =\Psi_M(z)\G R_M(z) .
	\end{align*} 
Now, we apply $( I+\mathscr{A}^{i}_{m}+\mathscr{A}_{m+M}^{e})$ to the last identity and we get 
\begin{align*}
\left( I+\mathscr{A}_{m}^{i}+ \mathscr{A}_{m+M}^{e}  \right) \G\widetilde{ R}_{\text{MIT}}(z) = \left( I-\mathscr{A}_{m}^{i}\mathscr{A}_{m+M}^{e}-\mathscr{A}_{m+M}^{e}\mathscr{A}_{m}^{i} \right)\G R_M(z)=:(\Xi_M(z))^{-1} \G R_M(z).
\end{align*} 
where $\Xi_M(z) $ is given by \eqref{Opein1}. Then, thanks to Proposition \ref{PseudESt} we know that  for $M>M_0$ the operator $(\Xi_M(z))^{-1} $ is bounded invertible  from $\mathit{H}^{1/2}(\S)^4$ into itself, which actually means that $\Psi_M$ is  bounded invertible  from $\mathit{H}^{1/2}(\S)^4$ into itself,  and that
$$ \Psi_M^{-1}= \Xi_M(z)\left( I+\mathscr{A}_{m}^{i}+ \mathscr{A}_{m+M}^{e}  \right).$$ 
From which it follows that  
$$ \G R_M(z)=  \Psi_M^{-1}(z)\G\widetilde{ R}_{\text{MIT}}(z).$$ 
Substituting this into formula \eqref{preformula} yields that
\begin{align*}
	R_M(z)= \widetilde{ R}_{\text{MIT}}(z)+E_M(z) \Psi_M^{-1}(z)\G\widetilde{ R}_{\text{MIT}}(z),
\end{align*} 
which achieves the proof of the theorem.\qed\\

As an immediate consequence of Theorem  \ref{Formulederes} and Proposition \ref{PseudESt} we have:
\begin{corollary}\label{une autre est} There is $M_0>0$ such that for every  $M>M_0$ and all $ z\in\rho(\Hm)\cap\rho(H_M)$,   the operators $\Xi^{\pm}_{M}(z):P_{\pm}\mathit{H}^{s}(\S)^4\rightarrow P_{\pm}\mathit{H}^{s}(\S)^4$ defined by 
\begin{align*}
\Xi^{+}_{M}(z)=\left(I-\mathscr{A}^{i}_{m}\mathscr{A}^{e}_{m+M}\right)^{-1}\quad\text{ andÂ }\quad \Xi^{-}_{M}(z) &=\left(I-\mathscr{A}^{e}_{m+M}\mathscr{A}^{i}_{m}\right)^{-1},
\end{align*}
are everywhere defined and bounded for any $s\in\rr$,  and it holds that 
 \begin{align*}
 \left|\left| \Xi^{\pm}_{M}(z) \right|\right|_{P_\pm\mathit{H}^{s+r}(\S)^4 \rightarrow  P_\pm\mathit{H}^{s}(\S)^4 } \lesssim M^{-r}, \quad \forall r>0. 
 \end{align*}
Moreover, if $v\in\mathit{H}^{1}(\rr^3)^4$ solves $(D_m +M\beta 1_{\O_{e}} -z)v=e_{\O_{i}}f$, for some  $f\in\mathit{L}^{2}(\O_{i})^4$. Then,  $r_{\O_{i}}v$ satisfies the following boundary value problem 
\begin{equation}\label{sysprinc} 
\left\{
\begin{aligned}
(D_m- z) r_{\O_{i}}v &=f \quad&\text{ in }\,\Omega_{i},\\
\G_-v&= \Xi^{-}_{M}(z)\mathscr{A}^{e}_{m+M}\G_+R_{\text{MIT}}(z)f  \,\,&\text{ on }\, \S,\\
\G_+v&= \G_+R_{\text{MIT}}(z)f  +\mathscr{A}^{i}_{m}\G_-v \,\,&\text{ on }\, \S.
\end{aligned}
\right.
\end{equation}  
\end{corollary}
\textbf{Proof.} We first note that $\Xi^{\pm}_{M}(z)=P_\pm  \Xi_M(z) P_\pm $. Thus, the first statement follows immediately from Proposition \ref{PseudESt} . Now, let $f\in\mathit{L}^{2}(\O_{i})^4$,  and suppose that $v\in\mathit{H}^{1}(\rr^3)^4$ solves $(D_m +M\beta 1_{\O_{e}} -z)v=e_{\O_{i}}f$. Thus  $(D_m- z) r_{\O_{i}}v =f $ in  $\Omega_{i}$, and if we set 
\begin{align*}
	\varphi=P_-t_{\S}v\quad \text{ and} \quad \psi= P_+t_{\S}v,
\end{align*} 
then, from  \eqref{S3prime} we easily get 
$$ \var=    \Xi^{-}_{M}(z)\mathscr{A}^{e}_{m+M}\G_+R_{\text{MIT}}(z)f \quad\text{and}\quad  \psi= \G_+R_{\text{MIT}}(z)f  +\mathscr{A}^{i}_{m}\varphi, $$ 
which means that  $r_{\O_{i}}v$ satisfies  \eqref{sysprinc}, and this completes the proof of the corollary. \qed\\

\begin{remark}\label{matrixform}Notice that from \eqref{S3prime} we have that
\begin{align*}
 \begin{pmatrix}
\G_+r_{\O_{i}}R_M(z)f\\
 \G_-r_{\O_{e}}R_M(z)f 
\end{pmatrix}= \begin{pmatrix}
\Xi^{+}_{M}(z)\\
\Xi^{-}_{M}(z)
\end{pmatrix}\begin{pmatrix}
\mathit{I}_4& \mathscr{A}^{i}_{m}\\
 \mathscr{A}^{e}_{m+M}&	\mathit{I}_4
	\end{pmatrix}\begin{pmatrix}
\G_+R_{\text{MIT}}(z)r_{\O_{i}}f \\
\G_-\widetilde{ R}_{M}(z)r_{\O_{e}}f
\end{pmatrix}.
\end{align*}
With this observation, we remark that the resolvent formula \eqref{fullreso} can also be written in the following matrix form
\begin{align*}
 \begin{pmatrix}
r_{\O_{i}}R_M(z)\\
r_{\O_{e}}R_M(z)
\end{pmatrix}= \begin{pmatrix}
R_{\text{MIT}}(z)r_{\O_{i}}\\
\widetilde{R}_{M}(z)r_{\O_{e}}
\end{pmatrix}+  \begin{pmatrix}
E^{\O_{i}}_{m}(z)\Xi^{-}_{M}(z)\\
E^{\O_{e}}_{m+M}(z)\Xi^{+}_{M}(z)
\end{pmatrix}\begin{pmatrix}
\mathscr{A}^{e}_{m+M}& \mathit{I}_4\\
\mathit{I}_4& \mathscr{A}^{i}_{m}
 \end{pmatrix}\begin{pmatrix}
\G_+R_{\text{MIT}}(z)r_{\O_{i}} \\
\G_-\widetilde{ R}_{M}(z)r_{\O_{e}}
\end{pmatrix}.
\end{align*}
\\
\end{remark}
An inspection of the proof of Theorem \ref{Formulederes} shows that, for any $M>0$,  $ z\in\rho(\Hm)\cap\rho(H_{M})$ and  $f\in\mathit{L}^{2}(\rr^3)^4$, one has 
\begin{align}\label{relating1}
	\G\widetilde{ R}_{\text{MIT}}(z)f=\Psi_M(z)\G R_M(z)f. 
\end{align} 
When $f$ runs through the whole space $\mathit{L}^{2}(\rr^3)^4$, then the values of $\G\widetilde{ R}_{\text{MIT}}(z)f$ and $\G R_M(z)f$ cover the whole space $\mathit{H}^{1/2}(\S)^4$, which means that $\mathrm{Rn}(\Psi_M(z))= \mathit{H}^{1/2}(\S)^4$. Hence, if one proves that $\mathrm{Kr}(\Psi_M(z))= \{0\}$, then $\Psi_M(z)$ would be boundedly invertible in $ \mathit{H}^{1/2}(\S)^4$, and thus \eqref{fullreso} holds without restriction on $M>0$. The following theorem provides a Birman-Schwinger-type principle relating $\mathrm{Kr}(H_{M}-z)$ with $\mathrm{Kr}(\Psi_M(z))$ and allows us to recover  the resolvent formula \eqref{fullreso} for any $M>0$.

 \begin{theorem}\label{Formulederes1} Let $M>0$ and let $\Psi_M$ be as in Theorem \ref{Formulederes}. Then, the following hold:
 \begin{itemize}
  \item[(i)] For any $a\in( -(m+M),m+M)\cap\rho(\Hm)$ we have $a\in \mathrm{Sp_{p}}(H_{M})\Leftrightarrow0\in \mathrm{Sp_{p}}(\Psi_{M}(a))$, and  it holds that 
\begin{align*}
\mathrm{Kr}(H_{M}-a)=\{ E_M(a)g: g\in \mathrm{Kr}(\Psi_M(a))\}.
\end{align*} 
In particular, $\mathrm{dim}\mathrm{Kr}(H_{M}-a)= \mathrm{dim}\mathrm{Kr}(\Psi_M(a))$ holds for all $a\in( -(m+M),m+M)\cap\rho(\Hm)$.
 \item[(ii)]  The operator $\Psi_M(z)$ is boundedly invertible in $\mathit{H}^{1/2}(\S)^4$  for all $ z\in\rho(\Hm)\cap\rho(H_{M})$, and the following resolvent formula holds:
 \begin{align}\label{fullreso2}
R_M(z)= \widetilde{ R}_{\text{MIT}}(z)+E_M(z)\Psi^{-1}_M(z)\G\widetilde{ R}_{\mathrm{MIT}}(z). 
\end{align} 
\end{itemize}
\end{theorem}
\textbf{Proof.}  ($\mathrm{i}$)  Let us first prove the implication $(\Longrightarrow)$.  Let $a\in (-(m+M),m+M)\cap\rho(\Hm)$ be such that $(H_M-a)\var=0$ for some $0\neq\var\in\mathit{H}^1(\rr^3)^4$. Set $\var_+=\var_{|\O_i}$ and $\var_-=\var_{|\O_e}$. Then, it is clear that $\var_+$ solves the system \eqref{L1b} with $g=\G_- \var$, and $\var_-$  solves the system \eqref{L2b} with $h=\G_+\var$. Thus, $\var_+= E_{m}^{\O_{i}}(a)\G_- \var$  and $\var_-= E_{m+M}^{\O_{e}}(a)\G_+ \var$. Hence,  $\var= E_M(a) t_{\S}\varphi$ and  $\G_\pm\varphi\neq0$, as otherwise $\varphi$ would be zero. Using this and the definition of the Poincar\'e-Steklov operators, we obtain that 
\begin{align*}
(\mathit{I}_{4}+ \mathscr{A}^{i}_{m}) \G_-\varphi=:t_{\S}\varphi_+  =t_{\S}\varphi=  t_{\S}\varphi_-:=  (\mathit{I}_{4} +\mathscr{A}_{m+M}^{e})\G_+\varphi,
\end{align*}
and since $t_{\S}\varphi\neq 0$ it follows that 
$$\Psi_{M}(a)t_{\S}\varphi=(\mathit{I}_{4}- \mathscr{A}^{i}_{m}-\mathscr{A}_{m+M}^{e})t_{\S}\varphi=0,$$
which means that $0\in \mathrm{Sp_{p}}(\Psi_{M}(a))$ and proves the inclusion $\mathrm{Kr}(H_{M}-a)\subset\{ E_M(a)g: g\in \mathrm{Kr}(\Psi_M(a))\}$. 

Now we turn to the proof of the implication $(\Longleftarrow)$.  Let $a\in (-(m+M),m+M)\cap\rho(\Hm)$  and assume that $0$ is an eigenvalue of $\Psi_{M}(a)$. Then,  there is $g\in \mathit{H}^{1/2}(\S)^{4}\setminus \lbrace 0\rbrace$ such that $\Psi_{M}(a)g=0$ on $\S$. Note that this is equivalent to
\begin{align}\label{EQQV}
(P_- + \mathscr{A}^{i}_{m})g=(P_+ + \mathscr{A}_{m+M}^{e})g .
\end{align}
Since  $a\in (-(m+M),m+M) \cap \rho(\Hm)$, the operators  $E_{m}^{\O_{i}}(a): P_-\mathit{H}^{1/2}(\S)^{4}\rightarrow \mathit{H}^1(\O_i)^4$ and $E_{m+M}^{\O_{e}}(a):P_+\mathit{H}^{1/2}(\S)^{4}\rightarrow \mathit{H}^1(\O_e)^4$ are well-defined and bounded.  Thus, if we let $\var= E_{M}(a)g=(E_{m}^{\O_{i}}(a)P_-g,E_{m+M}^{\O_{e}}(a)P_+g)$, then $\var\neq0$ and we have that   $(D_m-a)\var=0$ in $\O_i$,  and that $(D_{m+M}-a)\var=0$ in $\O_e$. Hence, it remains to show that $\var\in\mathit{H}^1(\rr^3)^4$. For this, observe that by \eqref{EQQV} we have
\begin{align*}
t_{\S}E_{m}^{\O_{i}}(a)P_-g=(P_- + \mathscr{A}^{i}_{m})g= (P_+ +\mathscr{A}_{m+M}^{e})g = t_{\S}E_{m+M}^{\O_{e}}(a)P_+g.
\end{align*}
Thanks to the boundedness properties of $E_{m}^{\O_{i}}(a)$ and $E_{m+M}^{\O_{e}}(a)$, it follows from the above computations that $\var=E_M(a)g\in\mathit{H}^1(\rr^3)^4\setminus\{0\}$ and satisfies the equation $(H_M-a)\var=0$. Therefore, $a\in \mathrm{Sp_{p}}(H_M)$ and the inclusion $\{ E_M(a)g: g\in \mathrm{Kr}(\Psi_M(a))\}\subset\mathrm{Kr}(H_{M}-a)$ holds, which completes the proof of ($\mathrm{i}$).

($\mathrm{ii}$) Let $ z\in\rho(\Hm)\cap\rho(H_{M})$ and note that the self-adjointness of  $H_M$ together with assertion ($\mathrm{i}$) imply that $\mathrm{Kr}(\Psi_M(z))= \{0\}$, as otherwise $\mathrm{Kr}(H_{M}-z)\neq\{0\}$. Since $\mathrm{Rn}(\Psi_M(z))= \mathit{H}^{1/2}(\S)^4$ holds for all $ z\in\rho(\Hm)\cap\rho(H_{M})$, it follows that $\Psi_M(z)$ admits a bounded and everywhere defined inverse in $\mathit{H}^{1/2}(\S)^4$. Therefore, \eqref{relating1} yields that $ \G R_M(z)=  \Psi_M^{-1}(z)\G\widetilde{ R}_{\text{MIT}}(z)$, and the resolvent formula \eqref{fullreso2} follows from this and \eqref{preformula}.\qed\\

\begin{remark}Note the different nature of Theorems  \ref{Formulederes} and \ref{Formulederes1}, since the second one ensures the invertibility of $\Psi_M$ and yields the resolvent formula \eqref{fullreso2} without assumption, while the first one is based on a largeness assumption that allows us (thanks to  the semiclassical properties of the PS operators)  to obtain the explicit formula of the operator $(\Psi_M)^{-1}$. Besides,  note that in Theorem \ref{Formulederes1} we do not know a priori whether $(\Psi_M)^{-1}$ is uniformly bounded when $M$ is large, and hence \eqref{fullreso2} is not suitable for studying the large coupling convergence. 
\end{remark}

 In the next proposition we prove the norm convergence of $R_M(z)$  toward  $R_{\text{MIT}}(z)$ and estimate the rate of convergence.
	\begin{proposition}\label{Resolvent convergence}  For any compact set  $K\subset\rho(\Hm)$  there is $M_0>0$ such that for all $M>M_0$: $K\subset\rho(H_{M})$, and for all  $z\in K$ the resolvent $R_M$ admits an asymptotic expansion in $\mathcal{L}(\mathit{L}^2(\rr^3)^4)$ of the form:
	\begin{align}\label{asymptotic}
	R_M(z)=e_{\O_{i}}R_{\text{MIT}}(z)r_{\O_{i}}+ \frac{1}{M}\left( K_M(z) +  L_M(z)\right), 
	\end{align} 
where $K_M(z),\,L_M(z): \mathit{L}^2(\rr^3)^4\longrightarrow \mathit{L}^2(\rr^3)^4$ are uniformly bounded with respect to $M$ and satisfy
$$r_{\O_{i}}K_M(z) e_{\O_{i}}  =0=r_{\O_{e}}K_M(z) e_{\O_{e}}.$$
In particular, it holds that	
\begin{align}\label{RC1}
\left|\left| R_M(z)- e_{\O_{i}}R_{\text{MIT}}(z)r_{\O_{i}}\right|\right|_{\mathit{L}^2(\rr^3)^4\rightarrow \mathit{L}^2(\rr^3)^4}&=\mathcal{O}\left( \frac{1}{M}\right).
\end{align} 
\end{proposition}
Before giving the proof, we need the following estimates.
\begin{lemma}\label{Estimations des op}  Let $K\subset \cc$ be a compact set. Then, there is $M_0>0$ such that for all $M>M_0$:  $K\subset\rho(\widetilde{H}_{M})$ and for every $z\in K$ the following estimates hold:
\begin{align*}
\begin{split}
 \left|\left| \widetilde{ R}_{M}(z)f\right|\right|_{\mathit{L}^2(\O_{e})^4}+ \frac{1}{\sqrt{M}} \left|\left| \G_-\widetilde{ R}_{M}(z) f\right|\right|_{\mathit{L}^2(\S)^4} &\lesssim \frac{1}{M} \left|\left| f \right|\right|_{\mathit{L}^2(\O_{e})^4},\quad \forall f\in \mathit{L}^2(\O_{e})^4,\\
 \left|\left| \G_-\widetilde{ R}_{M}(z) f\right|\right|_{\mathit{H}^{-1/2}(\S)^4}&\lesssim \frac{1}{M} \left|\left| f \right|\right|_{\mathit{L}^2(\O_{e})^4}, \quad \forall f\in \mathit{L}^2(\O_{e})^4,\\
   \left|\left| E^{\O_{e}}_{m+M}(z)\psi \right|\right|_{\mathit{L}^2(\O_{e})^4}&\lesssim \frac{1}{\sqrt{M}}\left|\left| \psi \right|\right|_{\mathit{L}^{2}(\S)^4},\quad\forall \psi\in P_+\mathit{L}^{2}(\S)^4,\\
     \left|\left| E^{\O_{e}}_{m+M}(z)\psi \right|\right|_{\mathit{L}^2(\O_{e})^4}&\lesssim \frac{1}{M}\left|\left| \psi \right|\right|_{\mathit{H}^{1/2}(\S)^4},\quad\forall \psi\in P_+\mathit{H}^{1/2}(\S)^4. 
 \end{split}
 \end{align*}
\end{lemma}	
\textbf{Proof.} Fix a compact set $K\subset\cc$,  and note that for $M_1> \sup_{z\in K}\{  |\mathrm{Re}(z)| -m\} $ it holds that $K\subset\rho(D_{m+M_1})$, and hence $K\subset\rho(\widetilde{H}_{M})$  for all $M>M_1$. 
We next show the claimed estimates for $\widetilde{ R}_{M}(z)$ and $\G_-\widetilde{ R}_{M}(z)$. For this, let $z\in K$ and assume that $M>M_1$. Let $ \var\in \mathrm{dom}(\widetilde{H}_{M})$, then a straightforward application of the Green's formula yields that
\begin{align*}
\begin{split}
 \|  \widetilde{H}_{M}\var \|^2_{\mathit{L}^2(\O_{e})^4}=&  \| (\alpha\cdot\nabla)\var \|^2_{\mathit{L}^2(\O_e)^4}+(m+M)^2  \left|\left| \var \right|\right|^2_{\mathit{L}^2(\O_e)^4}+(m+M)  \left|\left| P_- t_{\S}\var \right|\right|^2_{\mathit{L}^2(\S)^4}.
\end{split}
\end{align*}
Using this and the Cauchy-Schwarz inequality we obtain that 
\begin{align*}
 \|  (\widetilde{H}_{M}-z)\var \|^2_{\mathit{L}^2(\O_{e})^4}=&  \|  \widetilde{H}_{M}\var \|^2_{\mathit{L}^2(\O_{e})^4} +|z|^2 \| \var \|^2_{\mathit{L}^2(\O_e)^4}  -2\mathrm{Re}(z) \langle  \widetilde{H}_{M}\var, \var \rangle_{\mathit{L}^2(\O_e)^4}\\
 \geqslant&  \|  \widetilde{H}_{M}\var \|^2_{\mathit{L}^2(\O_{e})^4} +|z|^2 \| \var \|^2_{\mathit{L}^2(\O_e)^4} -\frac{1}{2} \|  \widetilde{H}_{M}\var \|^2_{\mathit{L}^2(\O_{e})^4} -2|\mathrm{Re}(z)|^2 \| \var \|^2_{\mathit{L}^2(\O_e)^4}\\
 \geqslant& \left(\frac{(m+M)^2}{2}+|\mathrm{Im}(z)|^2 -|\mathrm{Re}(z)|^2\right)   \left|\left| \var \right|\right|^2_{\mathit{L}^2(\O_e)^4}+\frac{M}{2}  \left|\left| P_- t_{\S}\var \right|\right|^2_{\mathit{L}^2(\S)^4}.
\end{align*}
Therefore,   taking $\widetilde{ R}_{M}(z)f=\var$ and  $M\geqslant M_2\geqslant \sup_{z\in K}\{  \sqrt{|\mathrm{Re}(z)|^2-|\mathrm{Im}(z)|^2} -m\}$  we obtain the inequality 
\begin{align*}
 \left|\left| \widetilde{ R}_{M}(z)f\right|\right|_{\mathit{L}^2(\O_{e})^4} + \frac{1}{\sqrt{M}} \left|\left| \G_-\widetilde{ R}_{M}(z) f\right|\right|_{\mathit{L}^2(\S)^4}\lesssim \frac{1}{M} \left|\left| f \right|\right|_{\mathit{L}^2(\O_{e})^4}.
  \end{align*} 
Since $ \G_-$ is bounded from $\mathit{L}^2(\O_{e})^4$ into $\mathit{H}^{-1/2}(\S)^4$, it follows from the above inequality that 
\begin{align*}
 \left|\left| \G_-\widetilde{ R}_{M}(z) f\right|\right|_{\mathit{H}^{-1/2}(\S)^4}\lesssim  \left|\left| \G_-\right|\right|_{\mathit{L}^2(\O_{e})^4 \rightarrow\mathit{H}^{-1/2}(\S)^4} \left|\left|\widetilde{ R}_{M}(z) f\right|\right|_{\mathit{H}^{-1/2}(\S)^4}\lesssim   \frac{1}{M} \left|\left| f \right|\right|_{\mathit{L}^2(\O_{e})^4},
 \end{align*}
 for any $ f\in \mathit{L}^2(\O_{e})^4$, which gives the second inequality.
 
 Let us now turn to the proof of the claimed estimates for $E^{\O_{e}}_{m+M}(z)$.  Let $\psi\in P_+\mathit{L}^{2}(\S)^4$, then from the proof of Proposition \ref{definitionPSOP} we have 
 \begin{align*}
 \left|\left| \psi \right|\right|^2_{\mathit{L}^{2}(\S)^4} \geqslant  (m+M)\left|\left| E^{\O_{e}}_{m+M}(z)\psi \right|\right|^2_{\mathit{L}^2(\O_{e})^4}  -2|\mathrm{Re}(z)| \left|\left|E^{\O_{e}}_{m+M}(z)\psi \right|\right|^2_{\mathit{L}^2(\O_{e})^4}.
\end{align*}
 Thus, for any $M\geqslant M_3 \geqslant  \sup_{z\in K}\{ 4 |\mathrm{Re}(z)| -m\}$,  we get that 
 $$M  \left|\left|E^{\O_{e}}_{m+M}(z)\psi \right|\right|^2_{\mathit{L}^2(\O_{e})^4}\leqslant 2  \left|\left| \psi \right|\right|^2_{\mathit{L}^{2}(\S)^4},$$
  and this proves the first estimate for $ E^{\O_{e}}_{m+M}(z)$.  Finally, the last inequality is a consequence of the first one and Proposition \ref{definitionPSOP}. Indeed,  from Proposition \ref{definitionPSOP} $(\mathrm{ii})$ we know that $\beta\G_-\widetilde{ R}_{M}(\overline{z})$ is the adjoint of the operator $E^{\O_{e}}_{m+M}(z):P_+\mathit{H}^{1/2}(\S)^4\longrightarrow \mathit{L}^2(\O_{e})^4$.  Using this and the estimate fulfilled by $\G_-\widetilde{ R}_{M}(\overline{z})$  we obtain that 
 \begin{align*}
\left| \langle f,E^{\O_{e}}_{m+M}(z)\psi \rangle_{\mathit{L}^{2}(\O_{e})^4}\right| &=\left| \langle \G_-\widetilde{ R}_{M}(\overline{z})f,\beta\psi \rangle_{\mathit{H}^{-1/2}(\S)^4, \mathit{H}^{1/2}(\S)^4}\right|\\
&\leqslant  \left|\left| \G_-\widetilde{ R}_{M}(z) f\right|\right|_{\mathit{H}^{-1/2}(\S)^4}\left|\left| \psi \right|\right|_{\mathit{H}^{1/2}(\S)^4}\\
& \lesssim \frac{1}{M} \left|\left|  f\right|\right|_{\mathit{L}^{2}(\O_{e})^4}\left|\left| \psi \right|\right|_{\mathit{H}^{1/2}(\S)^4}.
\end{align*}
Since this is true for all $ f\in \mathit{L}^2(\O_{e})^4$, by duality arguments it follows that 
\begin{align*}
  \left|\left| E^{\O_{e}}_{m+M}(z)\psi \right|\right|_{\mathit{L}^2(\O_{e})^4}\lesssim \frac{1}{M}\left|\left| \psi \right|\right|_{\mathit{H}^{1/2}(\S)^4},\quad\forall \psi\in P_+\mathit{H}^{1/2}(\S)^4,
 \end{align*}
 which proves the last inequality. Hence, the lemma follows by taking $M_0= \max \{M_1,M_2, M_3\}$. \qed\\
\\

\textbf{Proof of Proposition \ref{Resolvent convergence}.} We first show \eqref{RC1} for some $M^{\prime}_0>0$ and any $z\in\cc\setminus \rr$. So, let us fix such a $z$ and let $f\in\mathit{L}^{2}(\rr^3)^4$. Then, it is clear that $z\in\rho(\Hm)\cap\rho(H_{M})$, and  from Theorem \ref{Formulederes} and Remark \ref{matrixform} we know that there is $M^{\prime}_0>0$ such that for all $M>M^{\prime}_0$ it holds that  
\begin{align*}
\left|\left| (R_M(z)- e_{\O_{i}}R_{\text{MIT}}(z)r_{\O_{i}})f\right|\right|_{\mathit{L}^2(\rr^3)^4}&\leqslant \left|\left|E^{\O_{i}}_{m}(z)\Xi^{-}_{M}(z)\mathscr{A}^{e}_{m+M}\G_+R_{\text{MIT}}(z)r_{\O_{i}} f\right|\right|_{\mathit{L}^2(\O_{i})^4} \\
&+ \left|\left|E^{\O_{i}}_{m}(z)\Xi^{-}_{M}(z)\G_-\widetilde{ R}_{M}(z)r_{\O_{e}}f\right|\right|_{\mathit{L}^2(\O_{i})^4}\\
&+ \left|\left|E^{\O_{e}}_{m+M}(z)\Xi^{+}_{M}(z)\G_+R_{\text{MIT}}(z)r_{\O_{i}} f\right|\right|_{\mathit{L}^2(\O_{e})^4} \\
&+ \left|\left|E^{\O_{e}}_{m+M}(z)\Xi^{+}_{M}(z)\mathscr{A}^{i}_{m}\G_-\widetilde{ R}_{M}(z)r_{\O_{e}}f\right|\right|_{\mathit{L}^2(\O_{e})^4}\\
 &+\left|\left|\widetilde{R}_{M}(z)r_{\O_{e}} f\right|\right|_{\mathit{L}^2(\O_{e})^4}=: J_1+J_2+J_3+J_4+J_5.
\end{align*} 
From  Lemma \ref{Estimations des op} we immediately get that $J_5 \lesssim M^{-1} \left|\left| f \right|\right|$. Next, notice that   $\G_+R_{\text{MIT}}(z): \mathit{L}^{2}(\O_i)^4\rightarrow \mathit{H}^{1/2}(\S)^4$, $\mathscr{A}^{i}_{m}: \mathit{H}^{1/2}(\S)^4\rightarrow \mathit{H}^{1/2}(\S)^4$ and  $E^{\O_{i}}_{m}(z):\mathit{H}^{-1/2}(\S)^4\rightarrow H(\alpha,\O_i)\subset  \mathit{L}^{2}(\O_i)^4$  (where $H(\alpha,\O_i)$ is defined by \eqref{Sobolev Dirac}) are bounded operators and do not depend on $M$.   Moreover, thanks to Corollary \ref{une autre est} we know that for all $s\in\rr$ there is $C>0$ independent of $M$ such that 
$$ \left|\left|\Xi^{\pm}_{M}(z)\right|\right|_{P_\pm\mathit{H}^{s}(\S)^4\rightarrow P_\pm\mathit{H}^{s}(\S)^4} \leqslant C.$$
Using this and the above observation, for $j\in\{1,2,3,4\}$,  we can estimate  $J_k$ as follows
\begin{align*}
J_1 &\lesssim   \left|\left| E^{\O_{i}}_{m}(z)\Xi^{-}_{M}(z) \right|\right|_{\mathit{H}^{-1/2}(\S)^4\rightarrow\mathit{L}^{2}(\O_i)^4}      \left|\left| \G_-\widetilde{ R}_{M}(z)r_{\O_{e}} f\right|\right|_{\mathit{H}^{-1/2}(\S)^4}, \\
J_2 &\lesssim    \left|\left|E^{\O_{e}}_{m+M}(z) \right|\right|_{\mathit{H}^{1/2}(\S)^4\rightarrow\mathit{L}^{2}(\O_e)^4}      \left|\left| \Xi^{+}_{M}(z)\G_+R_{\text{MIT}}(z)r_{\O_{i}} f\right|\right|_{\mathit{H}^{1/2}(\S)^4},\\
J_3 &\lesssim    \left|\left| E^{\O_{i}}_{m}(z)\Xi^{-}_{M}(z) \right|\right|_{P_-\mathit{H}^{-1/2}(\S)^4\rightarrow\mathit{L}^{2}(\O_i)^4}   \left|\left| \mathscr{A}^{e}_{m+M}\right|\right|_{\mathit{H}^{1/2}(\S)^4\rightarrow\mathit{H}^{-1/2}(\S)^4}    \left|\left| \G_+R_{\text{MIT}}(z)r_{\O_{i}} f\right|\right|_{\mathit{H}^{1/2}(\S)^4},\\
J_4 &\lesssim    \left|\left|E^{\O_{e}}_{m+M}(z) \right|\right|_{\mathit{L}^{2}(\S)^4\rightarrow\mathit{L}^{2}(\O_e)^4}   \left|\left|\Xi^{+}_{M}(z)\mathscr{A}^{i}_{m}\right|\right|_{\mathit{L}^{2}(\S)^4\rightarrow\mathit{L}^{2}(\S)^4}    \left|\left| \G_-\widetilde{ R}_{M}(z)r_{\O_{e}} f\right|\right|_{\mathit{L}^{2}(\S)^4}. 
\end{align*} 
Therefore, Theorem \ref{PseudESt}-($\mathrm{ii}$) together with Lemma \ref{Estimations des op} yield that 
\begin{align*}
J_k \lesssim \frac{1}{M} \left|\left| f \right|\right|_{\mathit{L}^{2}(\rr^3)^4}, \quad\text{for any } j\in\{1,2,3,4\}.
\end{align*}
Thus, we obtain the estimate 
 \begin{align}\label{LLLL}
\left|\left|  ( R_M(z)- e_{\O_{i}}R_{\mathrm{MIT}}(z)r_{\O_{i}})f \right|\right|_{\mathit{L}^{2}(\rr^3)^4}\leqslant \frac{C}{M} \left|\left| f \right|\right|_{\mathit{L}^{2}(\rr^3)^4}.
\end{align} 
 Moreover, the asymptotic expansion \eqref{asymptotic} holds with 
\begin{align*}
 L_M(z)=& M( e_{\O_{e}}\widetilde{R}_{M}(z)r_{\O_{e}} + e_{\O_{i}}  E^{\O_{i}}_{m}(z)\Xi^{-}_{M}(z)\mathscr{A}^{e}_{m+M}\G_+R_{\text{MIT}}(z)r_{\O_{i}} \\
  &+ e_{\O_{e}} E^{\O_{e}}_{m+M}(z)\Xi^{+}_{M}(z)\mathscr{A}^{i}_{m}\G_-\widetilde{ R}_{M}(z)r_{\O_{e}}),
\end{align*}  
and
\begin{align*}
	K_M(z)= M\left( e_{\O_{i}} E^{\O_{i}}_{m}(z)\Xi^{-}_{M}(z)\G_-\widetilde{ R}_{M}(z)r_{\O_{e}} +  e_{\O_{e}} E^{\O_{e}}_{m+M}(z)\Xi^{+}_{M}(z)\G_+R_{\text{MIT}}(z)r_{\O_{i}}\right),
\end{align*} 
and we clearly see that $r_{\O_{i}}K_M(z) e_{\O_{i}}  =0=r_{\O_{e}}K_M(z) e_{\O_{e}}$. 

Finally, since \eqref{LLLL} holds true for every $z\in\cc\setminus\rr$, for any fixed compact subset $K\subset\rho(\Hm)$, one can show by arguments similar to those in the proof of  \cite[Lemma A.1]{BCLS} that there is $M_0>M^{\prime}_0$ such that $K\subset \rho(H_M)$. Therefore, the proposition follows with the same arguments as before.\qed

  \subsection{Comments and  further remarks}
  In this part we discuss possible generalizations of our results and comment on the usefulness of the pseudodifferential properties of the Poincar\'e-Steklov operators.
\begin{itemize}
  \item[(1)]    First note that all the results in this article which are proved without the use of the (semi) classical properties of the Poincar\'e-Steklov
  operator are valid when $\S$ is just $\mathit{C}^{1,\o}$-smooth with $\o\in(1/2,1)$, and can also be generalized without difficulty to the case of local deformation of the plane $\rr^2\times\{0\}$ (see \cite{BB2} where the self-adjointness of $\Hm$ and the regularity properties of $\Phi^{\O}_{z,m}$, $\mathscr{C}_{z,m}$ and $\Lambda^{z}_{m}$ were shown for this case). We mention, however, that in the latter case the spectrum of the MIT bag operator is equal to that of the free Dirac operator, cf. \cite[Theorem 4.1]{BB2}. 
  \item[(2)] It should also be noted that there are several boundary conditions that lead to self-adjoint realizations of the Dirac operator on domains (see, e.g., \cite{AMSV,BHM,BB3}) and for which the associated PS operators can be analyzed in a similar way as for the MIT bag model. In particular, one can consider the PS operator $\mathscr{B}_m(z)$ associated with the self-adjoint Dirac operator 
 \begin{align*}
\widetilde{H}_{MIT}(m) v=D_mv\quad \forall v\in\mathrm{dom}(\widetilde{H}_{MIT}(m)):=  \left\{ v \in\mathit{H}^1(\Omega_{i} )^4 : P_+t_{\S} v=0 \text{ on }\S \right\}.
\end{align*} 
According to the previous considerations,  this operator can be viewed as an analogue of the Neumann-to-Dirichlet map for the Dirac operator.  Moreover,  the same arguments as in the proof of Theorem \ref{Th SymClas} show that 
\begin{align*}
\mathscr{B}(z)_m&= \frac{1}{\sqrt{-\Delta_\S}}S\cdot(\nabla_\S\wedge n)P_+\quad  \mathrm{mod} \, \mathit{Op} \sS^{-1}(\S) =  \frac{D_\Sigma}{\sqrt{-\Delta_\S}} P_+\quad  \mathrm{mod} \, \mathit{Op} \sS^{-1}(\S) ,
\end{align*}
 for all $z\in\rho(D_m)\cap\rho(\tilde{H}_{MIT}(m))$. 
  \item[(3)]   As already mentioned in the introduction, in \cite{BCLS} it was shown that (in the two-dimensional massless case) the norm resolvent convergence of $H_M$ to $\Hm$ holds with a convergence rate of $M^{-1/2}$. Their proof is based on two main ingredients: The first is a resolvent identity (see \cite[Lemma 2.2]{BCLS} for the exact formula),  and the second is the following inequality 
\begin{align}\label{GG} 
\left|\left| \G_- R_{M}(z) f\right|\right|_{\mathit{L}^{2}(\S)^4}\lesssim \frac{1}{\sqrt{M}} \left|\left| f \right|\right|_{\mathit{L}^2(\rr^3)^4},
\end{align} 
which is a consequence of the lower bound
  $$ \left|\left| \nabla\psi \right|\right|^2_{\mathit{L}^{2}(\O_e)^4} + M^2 \left|\left| \psi \right|\right|^2_{\mathit{L}^{2}(\O_e)^4} \geqslant (M-C) \left|\left| t_{\S}\psi \right|\right|^2_{\mathit{L}^{2}(\S)^4},$$ 
which holds for all $\psi\in\mathit{H}^1(\rr^3)^4$ and $M$ large enough (see \cite[Lemma 4]{SV} for the proof in the 2D-case and \cite[Proposition 2.1 (i)]{ALMR} for the 3D-case). Note that the resolvent formula \eqref{preformula} together with \eqref{GG} yield the same result. Indeed,  from  \eqref{S3} and  \eqref{GG} we easily get the inequality
\begin{align*}
\left|\left| \G_+ R_{M}(z) f\right|\right|_{\mathit{L}^{2}(\S)^4}\lesssim \left|\left| f \right|\right|_{\mathit{L}^2(\rr^3)^4}. 
\end{align*} 
This together with \eqref{preformula} and Lemma \ref{Estimations des op} yield
\begin{align*}
\left|\left| (R_M(z)- e_{\O_{i}}R_{\text{MIT}}(z)r_{\O_{i}})f\right|\right|_{\mathit{L}^2(\rr^3)^4}&\leqslant  \left|\left| E^{\O_{i}}_{m}(z)\G_-r_{\O_{e}}R_M(z)f \right|\right|_{\mathit{L}^2(\O_{i})^4}+\left|\left|\widetilde{R}_{M}(z)r_{\O_{e}} f\right|\right|_{\mathit{L}^2(\O_{e})^4}\\
&+ \left|\left| E^{\O_{e}}_{m+M}(z)\G_+r_{\O_{i}}R_M(z)f \right|\right|_{\mathit{L}^2(\O_{e})^4}\\
& \lesssim  \frac{1}{\sqrt{M}} \left|\left| f \right|\right|_{\mathit{L}^2(\rr^3)^4}.
\end{align*} 
\item[(4)] Finally, let us point out that a first order asymptotic expansion of the eigenvalues of $H_M$ in terms of the eigenvalues of $\Hm$ was established in \cite{ALMR} when $M\rightarrow\infty$. In their proof the authors used the min-max characterization and optimization techniques. Note that it is also possible to obtain such a result using the properties of the PS operator, the Krein formula from Theorem \ref{Formulederes} and the finite-dimensional perturbation theory (cf. Kato \cite{Kato} for example), see, e.g.,  \cite{BB1,BrCa02} for similar arguments.  Note also that the asymptotic expansion of the eigenvalues of $H_M$ depends only on the term $E^{\O_{i}}_{m}(z)\Xi^{-}_{M}(z)\mathscr{A}^{e}_{m+M}\G_+R_{\text{MIT}}(z)r_{\O_{i}}$. Indeed, let $\lambda_{MIT}$ be an eigenvalue of $\Hm$ with multiplicity $l$, and let $(f_1,\cdots,f_l)$ be an $L^2(\Omega_i)^4$-orthonormal basis of $\mathrm{Kr}(\Hm-\lambda_{MIT} \mathit{I}_4)$. Then, using the explicit resolvent formula from Remark \ref{matrixform} we see that  
  \begin{align*}
  \langle  R_{M}(z)e_{\O_i}f_k,  e_{\O_i}f_j\rangle_{\mathit{L}^{2}(\rr^3)^4}&=  \langle   E^{\O_{i}}_{m}(z)\Xi^{-}_{M}(z)\mathscr{A}^{e}_{m+M}\G_+R_{\text{MIT}}(z)f_k,f_j\rangle_{\mathit{L}^{2}(\O_i)^4}\\
  &=  \langle  \Xi^{-}_{M}(z)\mathscr{A}^{e}_{m+M}\G_+R_{\text{MIT}}(z)f_k,-\beta \G_+R_{\text{MIT}}(\overline{z})f_j\rangle_{\mathit{L}^{2}(\S)^4}\\
 &=  \frac{1}{(z-\lambda_{MIT})^2}\langle  \Xi^{-}_{M}(z)\mathscr{A}^{e}_{m+M}\G_+f_k,-\beta \G_+f_j\rangle_{\mathit{L}^{2}(\S)^4},
  \end{align*}
   which means that $E^{\O_{i}}_{m}(z)\Xi^{-}_{M}(z)\mathscr{A}^{e}_{m+M}\G_+R_{\text{MIT}}(z)r_{\O_{i}}$ is the only term that intervenes in the asymptotic expansion of the eigenvalues of $H_M$. Besides, recall that the principal symbol of $\Xi^{-}_{M}(z)\mathscr{A}^{e}_{m+M}$ is given by 
  \begin{align*}
 q_M(x,\xi)=   -\frac{S\cdot(\xi\wedge n(x))P_+}{\sqrt{|\xi\wedge n(x)|^2+(m+M)^2} +|\xi\wedge n(x)|+(m+M)},
    \end{align*}
    and for $M>0$ large enough one has
      \begin{align*}
 q_M(x,\xi)=  -\frac{1}{2M}S\cdot(\xi\wedge n(x))P_+ \sum_{l=1}^{\infty}\frac{1}{{M}^{l+1}}p_{l} (x,\xi)P_+ , \quad p_{l} \in\sS^{-l}. 
    \end{align*}
Using this, we formally deduce that for sufficiently large $M$, $H_M$ has exactly $l$ eigenvalues $(\lambda^M_k)_{1\leqslant k\leqslant l}$  counted according to their multiplicities (in $B(\lambda_{MIT},\eta)$ with $B(\lambda_{MIT},\eta)\cap \mathrm{Sp}(\Hm)=\lbrace\lambda_{MIT}\rbrace$) and these eigenvalues admit an asymptotic expansion of the form
   \begin{align}\label{devasy}
  \lambda^M_k = \lambda_{MIT}+\frac{1}{M}\mu_k +  \sum_{j=2}^{N}\frac{1}{M^j}\mu_k^j + \textit{O}\left(M^{-(N+1)}\right).
  \end{align}
   where $(\mu_k)_{1\leqslant k\leqslant l}$ are the eigenvalues of the matrix $\mathcal{M}$ with coefficients:
  $$m_{kj}=\frac{1}{2}\langle  \beta Op(S\cdot(\xi\wedge n(x)) )\G_+ f_k, \G_+f_j\rangle_{\mathit{L}^{2}(\S)^4}.$$
\end{itemize}


\section*{Acknowledgement} B. Benhellal was supported by the ERC-2014-ADG Project HADE Id. 669689 (European Research Council) and by the Spanish State Research Agency through BCAM Severo Ochoa excellence accreditation SEV- 2017-0718. The authors are grateful to Luis Vega for stimulating discussions on this work and thank Grigori Rozenblum  for pointing out the references on the R-matrix theory.


\end{document}